\documentclass[a4paper,11pt,leqno]{article}
\usepackage{amsfonts}
\usepackage{amsmath}
\usepackage{amssymb}
\numberwithin{equation}{section}

\makeatletter
\def\diagram{\leftwidth=\z@ \rightwidth=\z@ \topheight=\z@
\botheight=\z@ \setbox\@picbox\hbox\bgroup}
\def\enddiagram{\egroup\wd\@picbox\rightwidth\unitlength
\ht\@picbox\topheight\unitlength \dp\@picbox\botheight\unitlength
\hskip\leftwidth\unitlength\box\@picbox}
\def\bfig{\begin{diagram}}
\def\efig{\end{diagram}}
\newcount\wideness \newcount\leftwidth \newcount\rightwidth
\newcount\highness \newcount\topheight \newcount\botheight
\def\ratchet#1#2{\ifnum#1<#2 \global #1=#2 \fi}
\def\putbox(#1,#2)#3{%
\horsize{\wideness}{#3} \divide\wideness by 2 {\advance\wideness
by #1 \ratchet{\rightwidth}{\wideness}} {\advance\wideness by -#1
\ratchet{\leftwidth}{\wideness}} \vertsize{\highness}{#3}
\divide\highness by 2 {\advance\highness by #2
\ratchet{\topheight}{\highness}} {\advance\highness by -#2
\ratchet{\botheight}{\highness}} \put(#1,#2){\makebox(0,0){$#3$}}}
\def\putlbox(#1,#2)#3{%
\horsize{\wideness}{#3} {\advance\wideness by #1
\ratchet{\rightwidth}{\wideness}} {\ratchet{\leftwidth}{-#1}}
\vertsize{\highness}{#3} \divide\highness by 2 {\advance\highness
by #2 \ratchet{\topheight}{\highness}} {\advance\highness by -#2
\ratchet{\botheight}{\highness}}
\put(#1,#2){\makebox(0,0)[l]{$#3$}}}
\def\putrbox(#1,#2)#3{%
\horsize{\wideness}{#3} {\ratchet{\rightwidth}{#1}}
{\advance\wideness by -#1 \ratchet{\leftwidth}{\wideness}}
\vertsize{\highness}{#3} \divide\highness by 2 {\advance\highness
by #2 \ratchet{\topheight}{\highness}} {\advance\highness by -#2
\ratchet{\botheight}{\highness}}
\put(#1,#2){\makebox(0,0)[r]{$#3$}}}

\def\adjust[#1]{} 
\newcount \coefa
\newcount \coefb
\newcount \coefc
\newcount\tempcounta
\newcount\tempcountb
\newcount\tempcountc
\newcount\tempcountd
\newcount\xext
\newcount\yext
\newcount\xoff
\newcount\yoff
\newcount\gap%
\newcount\arrowtypea
\newcount\arrowtypeb
\newcount\arrowtypec
\newcount\arrowtyped
\newcount\arrowtypee
\newcount\height
\newcount\width
\newcount\xpos
\newcount\ypos
\newcount\run
\newcount\rise
\newcount\arrowlength
\newcount\halflength
\newcount\arrowtype
\newdimen\tempdimen
\newdimen\xlen
\newdimen\ylen
\newsavebox{\tempboxa}%
\newsavebox{\tempboxb}%
\newsavebox{\tempboxc}%
\newdimen\w@dth
\def\setw@dth#1#2{\setbox\z@\hbox{$#1$}\w@dth=\wd\z@
\setbox\@ne\hbox{$#2$}\ifnum\w@dth<\wd\@ne \w@dth=\wd\@ne \fi
\advance\w@dth by 1.2em}
\def\t@^#1_#2{\def\n@one{#1}\def\n@two{#2}\mathrel{\setw@dth{#1}{#2}
\mathop{\hbox to \w@dth{\rightarrowfill}}\limits
\ifx\n@one\empty\else ^{\box\z@}\fi \ifx\n@two\empty\else
_{\box\@ne}\fi}}
\def\t@@^#1{\@ifnextchar_ {\t@^{#1}}{\t@^{#1}_{}}}
\def\to{\@ifnextchar^ {\t@@}{\t@@^{}}}
\def\t@left^#1_#2{\def\n@one{#1}\def\n@two{#2}\mathrel{\setw@dth{#1}{#2}
\mathop{\hbox to \w@dth{\leftarrowfill}}\limits
\ifx\n@one\empty\else ^{\box\z@}\fi \ifx\n@two\empty\else
_{\box\@ne}\fi}}
\def\t@@left^#1{\@ifnextchar_ {\t@left^{#1}}{\t@left^{#1}_{}}}
\def\toleft{\@ifnextchar^ {\t@@left}{\t@@left^{}}}
\def\two@^#1_#2{\def\n@one{#1}\def\n@two{#2}\mathrel{\setw@dth{#1}{#2}
\mathop{\vcenter{\hbox to \w@dth{\rightarrowfill}\kern-1.7ex
                 \hbox to \w@dth{\rightarrowfill}}%
       }\limits
\ifx\n@one\empty\else ^{\box\z@}\fi \ifx\n@two\empty\else
_{\box\@ne}\fi}}
\def\tw@@^#1{\@ifnextchar_ {\two@^{#1}}{\two@^{#1}_{}}}
\def\two{\@ifnextchar^ {\tw@@}{\tw@@^{}}}
\def\tofr@^#1_#2{\def\n@one{#1}\def\n@two{#2}\mathrel{\setw@dth{#1}{#2}
\mathop{\vcenter{\hbox to \w@dth{\rightarrowfill}\kern-1.7ex
                 \hbox to \w@dth{\leftarrowfill}}%
       }\limits
\ifx\n@one\empty\else ^{\box\z@}\fi \ifx\n@two\empty\else
_{\box\@ne}\fi}}
\def\t@fr@^#1{\@ifnextchar_ {\tofr@^{#1}}{\tofr@^{#1}_{}}}
\def\tofro{\@ifnextchar^ {\t@fr@}{\t@fr@^{}}}

\def\mon{\mathop{\m@th\hbox to
      14.6\P@{\lasyb\char'51\hskip-2.1\P@$\arrext$\hss
$\mathord\rightarrow$}}\limits} 
\def\leftmono{\mathrel{\m@th\hbox to
14.6\P@{$\mathord\leftarrow$\hss$\arrext$\hskip-2.1\P@\lasyb\char'50%
}}\limits} 
\mathchardef\arrext="0200       

\def\settypes(#1,#2,#3){\arrowtypea#1 \arrowtypeb#2 \arrowtypec#3}
\def\settoheight#1#2{\setbox\@tempboxa\hbox{#2}#1\ht\@tempboxa\relax}%
\def\settodepth#1#2{\setbox\@tempboxa\hbox{#2}#1\dp\@tempboxa\relax}%
\def\settokens[#1`#2`#3`#4]{%
     \def\tokena{#1}\def\tokenb{#2}\def\tokenc{#3}\def\tokend{#4}}
\def\setsqparms[#1`#2`#3`#4;#5`#6]{%
\arrowtypea #1 \arrowtypeb #2 \arrowtypec #3 \arrowtyped #4
\width #5 \height #6 }
\def\setpos(#1,#2){\xpos=#1 \ypos#2}

\def\settriparms[#1`#2`#3;#4]{\settripairparms[#1`#2`#3`1`1;#4]}%
\def\settripairparms[#1`#2`#3`#4`#5;#6]{%
\arrowtypea #1 \arrowtypeb #2 \arrowtypec #3 \arrowtyped #4
\arrowtypee #5 \width #6 \height #6 }
\def\resetparms{\settripairparms[1`1`1`1`1;500]\width 500}
\resetparms
\def\mvector(#1,#2)#3{
\put(0,0){\vector(#1,#2){#3}}%
\put(0,0){\vector(#1,#2){26}}%
}
\def\evector(#1,#2)#3{{
\arrowlength #3
\put(0,0){\vector(#1,#2){\arrowlength}}%
\advance \arrowlength by-30
\put(0,0){\vector(#1,#2){\arrowlength}}%
}}
\def\horsize#1#2{%
\settowidth{\tempdimen}{$#2$}%
#1=\tempdimen \divide #1 by\unitlength }
\def\vertsize#1#2{%
\settoheight{\tempdimen}{$#2$}%
#1=\tempdimen
\settodepth{\tempdimen}{$#2$}%
\advance #1 by\tempdimen \divide #1 by\unitlength }
\def\putvector(#1,#2)(#3,#4)#5#6{{%
\ifnum3<\arrowtype \putdashvector(#1,#2)(#3,#4)#5\arrowtype \else
\ifnum\arrowtype<-3 \putdashvector(#1,#2)(#3,#4)#5\arrowtype \else
\xpos=#1 \ypos=#2 \run=#3 \rise=#4 \arrowlength=#5 \ifnum
\arrowtype<0
    \ifnum \run=0
        \advance \ypos by-\arrowlength
    \else
        \tempcounta \arrowlength
        \multiply \tempcounta by\rise
        \divide \tempcounta by\run
        \ifnum\run>0
            \advance \xpos by\arrowlength
            \advance \ypos by\tempcounta
        \else
            \advance \xpos by-\arrowlength
            \advance \ypos by-\tempcounta
        \fi
    \fi
    \multiply \arrowtype by-1
    \multiply \rise by-1
    \multiply \run by-1
\fi \ifcase \arrowtype
\or \put(\xpos,\ypos){\vector(\run,\rise){\arrowlength}}%
\or \put(\xpos,\ypos){\mvector(\run,\rise)\arrowlength}%
\or \put(\xpos,\ypos){\evector(\run,\rise){\arrowlength}}%
\fi\fi\fi }}
\def\putsplitvector(#1,#2)#3#4{
\xpos #1 \ypos #2 \arrowtype #4 \halflength #3 \arrowlength #3
\gap 140 \advance \halflength by-\gap \divide \halflength by2
\ifnum\arrowtype>0
   \ifcase \arrowtype
   \or \put(\xpos,\ypos){\line(0,-1){\halflength}}%
       \advance\ypos by-\halflength
       \advance\ypos by-\gap
       \put(\xpos,\ypos){\vector(0,-1){\halflength}}%
   \or \put(\xpos,\ypos){\line(0,-1)\halflength}%
       \put(\xpos,\ypos){\vector(0,-1)3}%
       \advance\ypos by-\halflength
       \advance\ypos by-\gap
       \put(\xpos,\ypos){\vector(0,-1){\halflength}}%
   \or \put(\xpos,\ypos){\line(0,-1)\halflength}%
       \advance\ypos by-\halflength
       \advance\ypos by-\gap
       \put(\xpos,\ypos){\evector(0,-1){\halflength}}%
   \fi
\else \arrowtype=-\arrowtype
   \ifcase\arrowtype
   \or \advance \ypos by-\arrowlength
       \put(\xpos,\ypos){\line(0,1){\halflength}}%
       \advance\ypos by\halflength
       \advance\ypos by\gap
       \put(\xpos,\ypos){\vector(0,1){\halflength}}%
   \or \advance \ypos by-\arrowlength
       \put(\xpos,\ypos){\line(0,1)\halflength}%
       \put(\xpos,\ypos){\vector(0,1)3}%
       \advance\ypos by\halflength
       \advance\ypos by\gap
       \put(\xpos,\ypos){\vector(0,1){\halflength}}%
   \or \advance \ypos by-\arrowlength
       \put(\xpos,\ypos){\line(0,1)\halflength}%
       \advance\ypos by\halflength
       \advance\ypos by\gap
       \put(\xpos,\ypos){\evector(0,1){\halflength}}%
   \fi
\fi }
\def\putmorphism(#1)(#2,#3)[#4`#5`#6]#7#8#9{{%
\run #2 \rise #3 \ifnum\rise=0
  \puthmorphism(#1)[#4`#5`#6]{#7}{#8}#9%
\else\ifnum\run=0
  \putvmorphism(#1)[#4`#5`#6]{#7}{#8}#9%
\else
\setpos(#1)%
\arrowlength #7 \arrowtype #8 \ifnum\run=0 \else\ifnum\rise=0
\else \ifnum\run>0
    \coefa=1
\else
   \coefa=-1
\fi \ifnum\arrowtype>0
   \coefb=0
   \coefc=-1
\else
   \coefb=\coefa
   \coefc=1
   \arrowtype=-\arrowtype
\fi \width=2 \multiply \width by\run \divide \width by\rise
\ifnum \width<0  \width=-\width\fi \advance\width by60 \if l#9
\width=-\width\fi
\putbox(\xpos,\ypos){#4}
{\multiply \coefa by\arrowlength
\advance\xpos by\coefa \multiply \coefa by\rise \divide \coefa
by\run \advance \ypos by\coefa
\putbox(\xpos,\ypos){#5} }%
{\multiply \coefa by\arrowlength
\divide \coefa by2 \advance \xpos by\coefa \advance \xpos by\width
\multiply \coefa by\rise \divide \coefa by\run \advance \ypos
by\coefa
\if l#9%
   \putrbox(\xpos,\ypos){#6}%
\else\if r#9%
   \putlbox(\xpos,\ypos){#6}%
\fi\fi }%
{\multiply \rise by-\coefc
\multiply \run by-\coefc \multiply \coefb by\arrowlength \advance
\xpos by\coefb \multiply \coefb by\rise \divide \coefb by\run
\advance \ypos by\coefb \multiply \coefc by70 \advance \ypos
by\coefc \multiply \coefc by\run \divide \coefc by\rise \advance
\xpos by\coefc \multiply \coefa by140 \multiply \coefa by\run
\divide \coefa by\rise \advance \arrowlength by\coefa
\ifcase\arrowtype
\or \put(\xpos,\ypos){\vector(\run,\rise){\arrowlength}}%
\or \put(\xpos,\ypos){\mvector(\run,\rise){\arrowlength}}%
\or \put(\xpos,\ypos){\evector(\run,\rise){\arrowlength}}%
\fi}\fi\fi\fi\fi}}

\newcount\numbdashes \newcount\lengthdash \newcount\increment
\def\howmanydashes{
\numbdashes=\arrowlength \lengthdash=40 \divide\numbdashes by
\lengthdash \lengthdash=\arrowlength \divide\lengthdash by
\numbdashes
\increment=\lengthdash \multiply\lengthdash by 3
\divide\lengthdash by 5 }
\def\putdashvector(#1)(#2,#3)#4#5{%
\ifnum#3=0 \putdashhvector(#1){#4}#5 \else \ifnum#2=0
\putdashvvector(#1){#4}#5\fi\fi}
\def\putdashhvector(#1,#2)#3#4{{%
\arrowlength=#3 \howmanydashes
\multiput(#1,#2)(\increment,0){\numbdashes}%
{\vrule height .4pt width \lengthdash\unitlength} \arrowtype=#4
\xpos=#1 \ifnum\arrowtype<0 \advance\arrowtype by 7 \fi
\ifcase\arrowtype \or \advance\xpos by 10
    \put(\xpos,#2){\vector(-1,0){\lengthdash}}
    \advance\xpos by 40
    \put(\xpos,#2){\vector(-1,0){\lengthdash}}
\or \advance \xpos by 10
    \put(\xpos,#2){\vector(-1,0){\lengthdash}}
    \advance\xpos by  \arrowlength
    \advance\xpos by  -50
    \put(\xpos,#2){\vector(-1,0){\lengthdash}}
\or \advance\xpos by 10
    \put(\xpos,#2){\vector(-1,0){\lengthdash}}
\or \advance\xpos by \arrowlength
    \advance\xpos by -\lengthdash
    \put(\xpos,#2){\vector(1,0){\lengthdash}}
\or {\advance\xpos by 10
    \put(\xpos,#2){\vector(1,0){\lengthdash}}}
    \advance\xpos by \arrowlength
    \advance\xpos by -\lengthdash
    \put(\xpos,#2){\vector(1,0){\lengthdash}}
\or \advance\xpos by \arrowlength
    \advance\xpos by -\lengthdash
    \put(\xpos,#2){\vector(1,0){\lengthdash}}
    \advance\xpos by -40
    \put(\xpos,#2){\vector(1,0){\lengthdash}}
   \fi
}}
\def\putdashvvector(#1,#2)#3#4{{%
\arrowlength=#3 \howmanydashes \ypos=#2 \advance\ypos by
-\arrowlength
\multiput(#1,#2)(0,\increment){\numbdashes}%
    {\vrule width .4pt height \lengthdash\unitlength}
\arrowtype=#4 \ypos=#2 \ifnum\arrowtype<0 \advance\arrowtype by 7
\fi \ifcase\arrowtype \or \advance\ypos by \arrowlength
\advance\ypos by -40
    \put(#1,\ypos){\vector(0,1){\lengthdash}}
    \advance\ypos by -40
    \put(#1,\ypos){\vector(0,1){\lengthdash}}
\or \advance\ypos by 10
    \put(#1,\ypos){\vector(0,1){\lengthdash}}
    \advance\ypos by \arrowlength \advance\ypos by -40
    \put(#1,\ypos){\vector(0,1){\lengthdash}}
\or \advance\ypos by \arrowlength \advance\ypos by -40
    \put(#1,\ypos){\vector(0,1){\lengthdash}}
\or \advance\ypos by 10
    \put(#1,\ypos){\vector(0,-1){\lengthdash}}
\or \advance\ypos by 10
    \put(#1,\ypos){\vector(0,-1){\lengthdash}}
    \advance\ypos by \arrowlength \advance\ypos by -40
    \put(#1,\ypos){\vector(0,-1){\lengthdash}}
\or \advance\ypos by 10
    \put(#1,\ypos){\vector(0,-1){\lengthdash}}
    \advance\ypos by 40
    \put(#1,\ypos){\vector(0,-1){\lengthdash}}
\fi }}
\def\puthmorphism(#1,#2)[#3`#4`#5]#6#7#8{{%
\xpos #1 \ypos #2 \width #6 \arrowlength #6 \arrowtype=#7
\putbox(\xpos,\ypos){#3\vphantom{#4}}%
{\advance \xpos by\arrowlength
\putbox(\xpos,\ypos){\vphantom{#3}#4}}%
\horsize{\tempcounta}{#3}%
\horsize{\tempcountb}{#4}%
\divide \tempcounta by2 \divide \tempcountb by2 \advance
\tempcounta by30 \advance \tempcountb by30 \advance \xpos
by\tempcounta \advance \arrowlength by-\tempcounta \advance
\arrowlength by-\tempcountb
\putvector(\xpos,\ypos)(1,0)\arrowlength\arrowtype \divide
\arrowlength by2 \advance \xpos by\arrowlength
\vertsize{\tempcounta}{#5}%
\divide\tempcounta by2 \advance \tempcounta by20
\if a#8 %
   \advance \ypos by\tempcounta
   \putbox(\xpos,\ypos){#5}%
\else
   \advance \ypos by-\tempcounta
   \putbox(\xpos,\ypos){#5}%
\fi}}
\def\putvmorphism(#1,#2)[#3`#4`#5]#6#7#8{{%
\xpos #1 \ypos #2 \arrowlength #6 \arrowtype #7
\settowidth{\xlen}{$#5$}%
\putbox(\xpos,\ypos){#3}%
{\advance \ypos by-\arrowlength
\putbox(\xpos,\ypos){#4}}%
{\advance\arrowlength by-140 \advance \ypos by-70 \ifdim\xlen>0pt
   \if m#8%
      \putsplitvector(\xpos,\ypos)\arrowlength\arrowtype
   \else
   \putvector(\xpos,\ypos)(0,-1)\arrowlength\arrowtype
   \fi
\else
   \putvector(\xpos,\ypos)(0,-1)\arrowlength\arrowtype
\fi}%
\ifdim\xlen>0pt
   \divide \arrowlength by2
   \advance\ypos by-\arrowlength
   \if l#8%
      \advance \xpos by-40
      \putrbox(\xpos,\ypos){#5}%
   \else\if r#8%
      \advance \xpos by40
      \putlbox(\xpos,\ypos){#5}%
   \else
      \putbox(\xpos,\ypos){#5}%
   \fi\fi
\fi }}
\def\putsquarep<#1>(#2)[#3;#4`#5`#6`#7]{{%
\setsqparms[#1]%
\setpos(#2)%
\settokens[#3]%
\puthmorphism(\xpos,\ypos)[\tokenc`\tokend`{#7}]{\width}{\arrowtyped}b%
\advance\ypos by \height
\puthmorphism(\xpos,\ypos)[\tokena`\tokenb`{#4}]{\width}{\arrowtypea}a%
\putvmorphism(\xpos,\ypos)[``{#5}]{\height}{\arrowtypeb}l%
\advance\xpos by \width
\putvmorphism(\xpos,\ypos)[``{#6}]{\height}{\arrowtypec}r%
}}
\def\putsquare{\@ifnextchar <{\putsquarep}{\putsquarep%
   <\arrowtypea`\arrowtypeb`\arrowtypec`\arrowtyped;\width`\height>}}
\def\square{\@ifnextchar< {\squarep}{\squarep
   <\arrowtypea`\arrowtypeb`\arrowtypec`\arrowtyped;\width`\height>}}
\def\squarep<#1>[#2`#3`#4`#5;#6`#7`#8`#9]{{
\setsqparms[#1]
\diagram
\putsquarep<\arrowtypea`\arrowtypeb`\arrowtypec`
\arrowtyped;\width`\height>
(0,0)[#2`#3`#4`{#5};#6`#7`#8`{#9}]
\enddiagram
}}                                                 
\def\putptrianglep<#1>(#2,#3)[#4`#5`#6;#7`#8`#9]{{%
\settriparms[#1]%
\xpos=#2 \ypos=#3 \advance\ypos by \height
\puthmorphism(\xpos,\ypos)[#4`#5`{#7}]{\height}{\arrowtypea}a%
\putvmorphism(\xpos,\ypos)[`#6`{#8}]{\height}{\arrowtypeb}l%
\advance\xpos by\height
\putmorphism(\xpos,\ypos)(-1,-1)[``{#9}]{\height}{\arrowtypec}r%
}}
\def\putptriangle{\@ifnextchar <{\putptrianglep}{\putptrianglep
   <\arrowtypea`\arrowtypeb`\arrowtypec;\height>}}
\def\ptriangle{\@ifnextchar <{\ptrianglep}{\ptrianglep
   <\arrowtypea`\arrowtypeb`\arrowtypec;\height>}}
\def\ptrianglep<#1>[#2`#3`#4;#5`#6`#7]{{
\settriparms[#1]
\diagram
\putptrianglep<\arrowtypea`\arrowtypeb`
\arrowtypec;\height>
(0,0)[#2`#3`#4;#5`#6`{#7}]
\enddiagram
}}                                            
\def\putqtrianglep<#1>(#2,#3)[#4`#5`#6;#7`#8`#9]{{%
\settriparms[#1]%
\xpos=#2 \ypos=#3 \advance\ypos by\height
\puthmorphism(\xpos,\ypos)[#4`#5`{#7}]{\height}{\arrowtypea}a%
\putmorphism(\xpos,\ypos)(1,-1)[``{#8}]{\height}{\arrowtypeb}l%
\advance\xpos by\height
\putvmorphism(\xpos,\ypos)[`#6`{#9}]{\height}{\arrowtypec}r%
}}
\def\putqtriangle{\@ifnextchar <{\putqtrianglep}{\putqtrianglep
   <\arrowtypea`\arrowtypeb`\arrowtypec;\height>}}
\def\qtriangle{\@ifnextchar <{\qtrianglep}{\qtrianglep
   <\arrowtypea`\arrowtypeb`\arrowtypec;\height>}}
\def\qtrianglep<#1>[#2`#3`#4;#5`#6`#7]{{
\settriparms[#1]
\width=\height                                
\diagram
\putqtrianglep<\arrowtypea`\arrowtypeb`
\arrowtypec;\height>
(0,0)[#2`#3`#4;#5`#6`{#7}]
\enddiagram
}}
\def\putdtrianglep<#1>(#2,#3)[#4`#5`#6;#7`#8`#9]{{%
\settriparms[#1]%
\xpos=#2 \ypos=#3
\puthmorphism(\xpos,\ypos)[#5`#6`{#9}]{\height}{\arrowtypec}b%
\advance\xpos by \height \advance\ypos by\height
\putmorphism(\xpos,\ypos)(-1,-1)[``{#7}]{\height}{\arrowtypea}l%
\putvmorphism(\xpos,\ypos)[#4``{#8}]{\height}{\arrowtypeb}r%
}}
\def\putdtriangle{\@ifnextchar <{\putdtrianglep}{\putdtrianglep
   <\arrowtypea`\arrowtypeb`\arrowtypec;\height>}}
\def\dtriangle{\@ifnextchar <{\dtrianglep}{\dtrianglep
   <\arrowtypea`\arrowtypeb`\arrowtypec;\height>}}
\def\dtrianglep<#1>[#2`#3`#4;#5`#6`#7]{{
\settriparms[#1]
\width=\height                                
\diagram
\putdtrianglep<\arrowtypea`\arrowtypeb`
\arrowtypec;\height>
(0,0)[#2`#3`#4;#5`#6`{#7}]
\enddiagram
}}
\def\putbtrianglep<#1>(#2,#3)[#4`#5`#6;#7`#8`#9]{{%
\settriparms[#1]%
\xpos=#2 \ypos=#3
\puthmorphism(\xpos,\ypos)[#5`#6`{#9}]{\height}{\arrowtypec}b%
\advance\ypos by\height
\putmorphism(\xpos,\ypos)(1,-1)[``{#8}]{\height}{\arrowtypeb}r%
\putvmorphism(\xpos,\ypos)[#4``{#7}]{\height}{\arrowtypea}l%
}}
\def\putbtriangle{\@ifnextchar <{\putbtrianglep}{\putbtrianglep
   <\arrowtypea`\arrowtypeb`\arrowtypec;\height>}}
\def\btriangle{\@ifnextchar <{\btrianglep}{\btrianglep
   <\arrowtypea`\arrowtypeb`\arrowtypec;\height>}}
\def\btrianglep<#1>[#2`#3`#4;#5`#6`#7]{{
\settriparms[#1]
\width=\height                               
\diagram
\putbtrianglep<\arrowtypea`\arrowtypeb`
\arrowtypec;\height>
(0,0)[#2`#3`#4;#5`#6`{#7}]
\enddiagram
}}
\def\putAtrianglep<#1>(#2,#3)[#4`#5`#6;#7`#8`#9]{{%
\settriparms[#1]%
\xpos=#2 \ypos=#3 {\multiply \height by2
\puthmorphism(\xpos,\ypos)[#5`#6`{#9}]{\height}{\arrowtypec}b}%
\advance\xpos by\height \advance\ypos by\height
\putmorphism(\xpos,\ypos)(-1,-1)[#4``{#7}]{\height}{\arrowtypea}l%
\putmorphism(\xpos,\ypos)(1,-1)[``{#8}]{\height}{\arrowtypeb}r%
}}
\def\putAtriangle{\@ifnextchar <{\putAtrianglep}{\putAtrianglep
   <\arrowtypea`\arrowtypeb`\arrowtypec;\height>}}
\def\Atriangle{\@ifnextchar <{\Atrianglep}{\Atrianglep
   <\arrowtypea`\arrowtypeb`\arrowtypec;\height>}}
\def\Atrianglep<#1>[#2`#3`#4;#5`#6`#7]{{
\settriparms[#1]
\width=\height                                     
\diagram
\putAtrianglep<\arrowtypea`\arrowtypeb`
\arrowtypec;\height>
(0,0)[#2`#3`#4;#5`#6`{#7}]
\enddiagram
}}
\def\putAtrianglepairp<#1>(#2)[#3;#4`#5`#6`#7`#8]{{%
\settripairparms[#1]%
\setpos(#2)%
\settokens[#3]%
\puthmorphism(\xpos,\ypos)[\tokenb`\tokenc`{#7}]{\height}{\arrowtyped}b%
\advance\xpos by\height
\puthmorphism(\xpos,\ypos)[\phantom{\tokenc}`\tokend`{#8}]%
{\height}{\arrowtypee}b%
\advance\ypos by\height
\putmorphism(\xpos,\ypos)(-1,-1)[\tokena``{#4}]{\height}{\arrowtypea}l%
\putvmorphism(\xpos,\ypos)[``{#5}]{\height}{\arrowtypeb}m%
\putmorphism(\xpos,\ypos)(1,-1)[``{#6}]{\height}{\arrowtypec}r%
}}
\def\putAtrianglepair{\@ifnextchar <{\putAtrianglepairp}{\putAtrianglepairp%
   <\arrowtypea`\arrowtypeb`\arrowtypec`\arrowtyped`\arrowtypee;\height>}}
\def\Atrianglepair{\@ifnextchar <{\Atrianglepairp}{\Atrianglepairp%
   <\arrowtypea`\arrowtypeb`\arrowtypec`\arrowtyped`\arrowtypee;\height>}}
\def\Atrianglepairp<#1>[#2;#3`#4`#5`#6`#7]{{
\settripairparms[#1]
\settokens[#2]
\width=\height                                
\diagram
\putAtrianglepairp                            
<\arrowtypea`\arrowtypeb`\arrowtypec`
\arrowtyped`\arrowtypee;\height>
(0,0)[{#2};#3`#4`#5`#6`{#7}]
\enddiagram
}}
\def\putVtrianglep<#1>(#2,#3)[#4`#5`#6;#7`#8`#9]{{%
\settriparms[#1]%
\xpos=#2 \ypos=#3 \advance\ypos by\height {\multiply\height by2
\puthmorphism(\xpos,\ypos)[#4`#5`{#7}]{\height}{\arrowtypea}a}%
\putmorphism(\xpos,\ypos)(1,-1)[`#6`{#8}]{\height}{\arrowtypeb}l%
\advance\xpos by\height \advance\xpos by\height
\putmorphism(\xpos,\ypos)(-1,-1)[``{#9}]{\height}{\arrowtypec}r%
}}
\def\putVtriangle{\@ifnextchar <{\putVtrianglep}{\putVtrianglep
   <\arrowtypea`\arrowtypeb`\arrowtypec;\height>}}
\def\Vtriangle{\@ifnextchar <{\Vtrianglep}{\Vtrianglep
   <\arrowtypea`\arrowtypeb`\arrowtypec;\height>}}
\def\Vtrianglep<#1>[#2`#3`#4;#5`#6`#7]{{
\settriparms[#1]
\width=\height                                 
\diagram
\putVtrianglep<\arrowtypea`\arrowtypeb`
\arrowtypec;\height>
(0,0)[#2`#3`#4;#5`#6`{#7}]
\enddiagram
}}
\def\putVtrianglepairp<#1>(#2)[#3;#4`#5`#6`#7`#8]{{
\settripairparms[#1]%
\setpos(#2)%
\settokens[#3]%
\advance\ypos by\height
\putmorphism(\xpos,\ypos)(1,-1)[`\tokend`{#6}]{\height}{\arrowtypec}l%
\puthmorphism(\xpos,\ypos)[\tokena`\tokenb`{#4}]{\height}{\arrowtypea}a%
\advance\xpos by\height
\puthmorphism(\xpos,\ypos)[\phantom{\tokenb}`\tokenc`{#5}]%
{\height}{\arrowtypeb}a%
\putvmorphism(\xpos,\ypos)[``{#7}]{\height}{\arrowtyped}m%
\advance\xpos by\height
\putmorphism(\xpos,\ypos)(-1,-1)[``{#8}]{\height}{\arrowtypee}r%
}}
\def\putVtrianglepair{\@ifnextchar <{\putVtrianglepairp}{\putVtrianglepairp%
    <\arrowtypea`\arrowtypeb`\arrowtypec`\arrowtyped`\arrowtypee;\height>}}
\def\Vtrianglepair{\@ifnextchar <{\Vtrianglepairp}{\Vtrianglepairp%
    <\arrowtypea`\arrowtypeb`\arrowtypec`\arrowtyped`\arrowtypee;\height>}}
\def\Vtrianglepairp<#1>[#2;#3`#4`#5`#6`#7]{{
\settripairparms[#1]
\settokens[#2]
\diagram
\putVtrianglepairp                             
<\arrowtypea`\arrowtypeb`\arrowtypec`
\arrowtyped`\arrowtypee;\height>
(0,0)[{#2};#3`#4`#5`#6`{#7}]
\enddiagram
}}

\def\putCtrianglep<#1>(#2,#3)[#4`#5`#6;#7`#8`#9]{{%
\settriparms[#1]%
\xpos=#2 \ypos=#3 \advance\ypos by\height
\putmorphism(\xpos,\ypos)(1,-1)[``{#9}]{\height}{\arrowtypec}l%
\advance\xpos by\height \advance\ypos by\height
\putmorphism(\xpos,\ypos)(-1,-1)[#4`#5`{#7}]{\height}{\arrowtypea}l%
{\multiply\height by 2
\putvmorphism(\xpos,\ypos)[`#6`{#8}]{\height}{\arrowtypeb}r}%
}}
\def\putCtriangle{\@ifnextchar <{\putCtrianglep}{\putCtrianglep
    <\arrowtypea`\arrowtypeb`\arrowtypec;\height>}}
\def\Ctriangle{\@ifnextchar <{\Ctrianglep}{\Ctrianglep
    <\arrowtypea`\arrowtypeb`\arrowtypec;\height>}}
\def\Ctrianglep<#1>[#2`#3`#4;#5`#6`#7]{{
\settriparms[#1]
\width=\height                               
\diagram
\putCtrianglep<\arrowtypea`\arrowtypeb`
\arrowtypec;\height>
(0,0)[#2`#3`#4;#5`#6`{#7}]
\enddiagram
}}                                           
\def\putDtrianglep<#1>(#2,#3)[#4`#5`#6;#7`#8`#9]{{%
\settriparms[#1]%
\xpos=#2 \ypos=#3 \advance\xpos by\height \advance\ypos by\height
\putmorphism(\xpos,\ypos)(-1,-1)[``{#9}]{\height}{\arrowtypec}r%
\advance\xpos by-\height \advance\ypos by\height
\putmorphism(\xpos,\ypos)(1,-1)[`#5`{#8}]{\height}{\arrowtypeb}r%
{\multiply\height by 2
\putvmorphism(\xpos,\ypos)[#4`#6`{#7}]{\height}{\arrowtypea}l}%
}}
\def\putDtriangle{\@ifnextchar <{\putDtrianglep}{\putDtrianglep
    <\arrowtypea`\arrowtypeb`\arrowtypec;\height>}}
\def\Dtriangle{\@ifnextchar <{\Dtrianglep}{\Dtrianglep
   <\arrowtypea`\arrowtypeb`\arrowtypec;\height>}}
\def\Dtrianglep<#1>[#2`#3`#4;#5`#6`#7]{{
\settriparms[#1]
\width=\height                              
\diagram
\putDtrianglep<\arrowtypea`\arrowtypeb`
\arrowtypec;\height>
(0,0)[#2`#3`#4;#5`#6`{#7}]
\enddiagram
}}                                          
\def\setrecparms[#1`#2]{\width=#1 \height=#2}%
\def\recursep<#1`#2>[#3;#4`#5`#6`#7`#8]{{%
\width=#1 \height=#2 \settokens[#3]
\settowidth{\tempdimen}{$\tokena$} \ifdim\tempdimen=0pt
  \savebox{\tempboxa}{\hbox{$\tokenb$}}%
  \savebox{\tempboxb}{\hbox{$\tokend$}}%
  \savebox{\tempboxc}{\hbox{$#6$}}%
\else
  \savebox{\tempboxa}{\hbox{$\hbox{$\tokena$}\times\hbox{$\tokenb$}$}}%
  \savebox{\tempboxb}{\hbox{$\hbox{$\tokena$}\times\hbox{$\tokend$}$}}%
  \savebox{\tempboxc}{\hbox{$\hbox{$\tokena$}\times\hbox{$#6$}$}}%
\fi \ypos=\height \divide\ypos by 2 \xpos=\ypos \advance\xpos by
\width \bfig
\putCtrianglep<-1`1`1;\ypos>(0,0)[`\tokenc`;#5`#6`{#7}]%
\puthmorphism(\ypos,0)[\tokend`\usebox{\tempboxb}`{#8}]{\width}{-1}b%
\puthmorphism(\ypos,\height)[\tokenb`\usebox{\tempboxa}`{#4}]{\width}{-1}a%
\advance\ypos by \width
\putvmorphism(\ypos,\height)[``\usebox{\tempboxc}]{\height}1r%
\efig }}
\def\recurse{\@ifnextchar <{\recursep}{\recursep<\width`\height>}}
\def\puttwohmorphisms(#1,#2)[#3`#4;#5`#6]#7#8#9{{%
\puthmorphism(#1,#2)[#3`#4`]{#7}0a \ypos=#2 \advance\ypos by 20
\puthmorphism(#1,\ypos)[\phantom{#3}`\phantom{#4}`#5]{#7}{#8}a
\advance\ypos by -40
\puthmorphism(#1,\ypos)[\phantom{#3}`\phantom{#4}`#6]{#7}{#9}b }}
\def\puttwovmorphisms(#1,#2)[#3`#4;#5`#6]#7#8#9{{%
\putvmorphism(#1,#2)[#3`#4`]{#7}0a \xpos=#1 \advance\xpos by -20
\putvmorphism(\xpos,#2)[\phantom{#3}`\phantom{#4}`#5]{#7}{#8}l
\advance\xpos by 40
\putvmorphism(\xpos,#2)[\phantom{#3}`\phantom{#4}`#6]{#7}{#9}r }}
\def\puthcoequalizer(#1)[#2`#3`#4;#5`#6`#7]#8#9{{%
\setpos(#1)%
\puttwohmorphisms(\xpos,\ypos)[#2`#3;#5`#6]{#8}11%
\advance\xpos by #8
\puthmorphism(\xpos,\ypos)[\phantom{#3}`#4`#7]{#8}1{#9} }}
\def\putvcoequalizer(#1)[#2`#3`#4;#5`#6`#7]#8#9{{%
\setpos(#1)%
\puttwovmorphisms(\xpos,\ypos)[#2`#3;#5`#6]{#8}11%
\advance\ypos by -#8
\putvmorphism(\xpos,\ypos)[\phantom{#3}`#4`#7]{#8}1{#9} }}
\def\putthreehmorphisms(#1)[#2`#3;#4`#5`#6]#7(#8)#9{{%
\setpos(#1) \settypes(#8)
\if a#9 %
     \vertsize{\tempcounta}{#5}%
     \vertsize{\tempcountb}{#6}%
     \ifnum \tempcounta<\tempcountb \tempcounta=\tempcountb \fi
\else
     \vertsize{\tempcounta}{#4}%
     \vertsize{\tempcountb}{#5}%
     \ifnum \tempcounta<\tempcountb \tempcounta=\tempcountb \fi
\fi \advance \tempcounta by 60
\puthmorphism(\xpos,\ypos)[#2`#3`#5]{#7}{\arrowtypeb}{#9}
\advance\ypos by \tempcounta
\puthmorphism(\xpos,\ypos)[\phantom{#2}`\phantom{#3}`#4]{#7}{\arrowtypea}{#9}
\advance\ypos by -\tempcounta \advance\ypos by -\tempcounta
\puthmorphism(\xpos,\ypos)[\phantom{#2}`\phantom{#3}`#6]{#7}{\arrowtypec}{#9}
}}
\def\setarrowtoks[#1`#2`#3`#4`#5`#6]{%
\def\toka{#1}
\def\tokb{#2}
\def\tokc{#3}
\def\tokd{#4}
\def\toke{#5}
\def\tokf{#6}
}
\def\hex{\@ifnextchar <{\hexp}{\hexp<1000`400>}}
\def\hexp<#1`#2>[#3`#4`#5`#6`#7`#8;#9]{%
\setarrowtoks[#9] \yext=#2 \advance \yext by #2 \xext=#1
\advance\xext by \yext \bfig
\putCtriangle<-1`0`1;#2>(0,0)[`#5`;\tokb``\tokd] \xext=#1
\yext=#2 \advance \yext by #2
\putsquare<1`0`0`1;\xext`\yext>(#2,0)[#3`#4`#7`#8;\toka```\tokf]
\advance \xext by #2
\putDtriangle<0`1`-1;#2>(\xext,0)[`#6`;`\tokc`\toke] \efig }
\begin{document}
\newtheorem{theorem}{Theorem}[section]
\newtheorem{lemma}[theorem]{Lemma}
\newtheorem{corollary}[theorem]{Corollary}
\newtheorem{conjecture}[theorem]{Conjecture}
\newtheorem{remark}[theorem]{Remark}
\newtheorem{definition}[theorem]{Definition}
\newtheorem{problem}[theorem]{Problem}
\newtheorem{example}[theorem]{Example}
\newtheorem{proposition}[theorem]{Proposition}
\title{{\bf RICCI ITERATIONS AND CANONICAL K\"{A}HLER-EINSTEIN CURRENTS ON 
 LOG CANONICAL PAIRS}}
\date{November 4, 2012}
\author{Hajime TSUJI\footnote{Partially supported by Grant-in-Aid for Scientific Reserch (S) 17104001}}
\maketitle
\begin{abstract}
\noindent In this article we construct a canonical K\"{a}hler-Einstein current 
on an arbitrary  LC (log canonical) pair of log general type as the limit of a sequence of canonical K\"{a}hler-Einstein currents on KLT(Kawamata log terminal) pairs 
of log general type. 
We call the volume form of the canonical K\"{a}hler-Einstein current 
the canonical measure of the LC pair.
We prove that the relative canonical measure on a projective family of LC pairs  of log general type defines a singular hermitian metric on the relative log canonical bundle 
and the metric has semipositive curvature in the sense of current. 
This is the first semipositivity result for relative log canonical 
bundles of a family of LC pairs in general dimension. 

Our proof depends on  certain Ricci iterations and dynamical systems of 
Bergman kernels.  
\end{abstract}
\tableofcontents
\section{Introduction}
In \cite{KE}, I have constructed the canonical K\"{a}hler-Einstein metric on 
a smooth canonically polarized variety in terms of the dynamical systems of Bergman kernels. 
This dynamical construction is considered to be a polynomial approximation of the canonical K\"{a}hler-Einstein metric and naturally leads us to the logarithmic pluri-subharmonicity of relative canonical measure defined by \cite{s-t} on an algebraic fiber space with nonnegative relative Kodaira dimension (\cite{canonical})\footnote{Throughout  this article, an algebraic fiber space $f : X \to Y$ means that 
$X$,$Y$ are smooth projective and $f$ is a dominant morphism with connected fibers}. 

In algebraic geometry, it is natural and useful to consider 
 KLT pairs instead of projective varieties themselves.  
For example the finite generation of canonical rings has been proven 
not only for smooth projective varieties but also KLT pairs (\cite{b-c-h-m}). 
Moreover the proof of the finite generation  in \cite{b-c-h-m} also depends on the generalization to KLT pairs.   
Another example is the canonical bundle formula  which reduces 
the canonical rings of projective varieties of arbitrary Kodaira dimension 
to the log canonical rings of KLT pairs of log general type (\cite{f-m}).
  
In this way generalizations to KLT pairs are quite essential in algebraic 
geometry.  
Roughly speaking in differential geometric context, to consider KLT pairs corresponds to consider orbifold 
structures on the projective varieties.   

But so far  most of the results for KLT pairs have not yet been generalized  
to the case of LC pairs. For example, the finite generation of log canonical 
rings of LC pairs has not yet been proven.  The main reason is that 
the usual tools to study KLT pairs, such as the branched covering trick, the 
variation of Hodge structures and  the Kawamata-Viehweg vanishing theorem are no longer effective in the case of LC pairs.  At present there are big obstacles to 
go beyond KLT pairs. 
      
The purpose of this article is to construct a canonical K\"{a}hler-Einstein current (resp. a canonical measure) on a LC pair of log general type (resp. a KLT pair with nonnegative log Kodaira dimension) and prove the logarithmic plurisubharmonicity of the relative canonical K\"{a}hler-Einstein volume form (resp. the canonical measure) on a projective family.   This is a natural extension of the results in \cite{KE,canonical} to the case of LC pairs of log general type.
More precisely, first we extend the reuslt in \cite{KE,canonical} to 
the case of KLT pairs and then passing to the limits, we extend the results
to the case of LC pairs. 
  
The new feature  here is the use of Ricci iterations.  
In fact, in this case we cannot construct the K\"{a}hler-Einstein current 
in terms of a single dynamical system of Bergman kernels and the 
construction requires another dynamical system, namely Ricci iterations.

The main motivation to construct canonical K\"{a}hler-Einstein currents 
in terms of dynamical systems of Bergman kernels is to prove the 
semipositivity of relative (log) canonical bundles or 
the direct images of the pluri(log)canonical systems in terms of the 
recent break through due to Bo Berndtsson 
on  the logarithmic pluri-subharmonic variation properties  of Bergman kernels (\cite{b1,b2,b3}). 

We note that the corresponding result has already  been obtained in the context of 
algebraic geometry much earlier by Y. Kawamata in terms of  variation of Hodge structures as follows.
Let $f : X \to Y$ be an algebraic fiber space over a projective curve $Y$.  In \cite{ka1}, Y. Kawamata
proved that for every positive integer $m$, $f_{*}\mathcal{O}_{X}(K_{X/Y}^{\otimes{m}})$ is semipositive in the sense that any quotient of it has semipositive degree on $Y$.   Since then, this result has been used extensively 
in algebraic geometry.  Later Kawamata generalized his result
to the case of KLT pairs (\cite{ka2},\cite[p.175,Theorem 1.2]{ka3}). In particular he proved the subadjunction theorem (\cite{ka2}). We note that the KLT condition is quite essential in  his argument, bacause he applied the theory of variation of Hodge structures 
(\cite{griff, sch}) after taking a suitable branched covering of 
the original variety.  Hence in this sense LC pairs are out of reach of the Hodge theory (Maybe we may use mixed Hodge theory. But it is not clear now.).   

The advantage of the differential geometric method in this article is two folds. First, it gives a canonical metric with semipositive curvature in the sense of current on the relative log canonical bundle of a 
projective family of LC pairs. The canonical property of the metric is quite essential to consider the moduli problems.  See \cite{canAZD} for other canonical metrics. 
Second, it gives a general method to consider  LC pairs as a limit of 
KLT pairs. 
Although LC pairs are out of reach of the Hodge theory, nevertheless intutively, a LC pair is a limit of a sequence of KLT pairs. 
Hence it is natural to expect that if one obtains a theorem 
for KLT pairs, then  the corresponding result for LC pairs 
would be obtained by taking a limit in a suitable sense.   
In fact we implement this philosophy in terms of  canonical K\"{a}hler-Einstein  currents on KLT pairs and obtain the semipositivity of 
the relative log canonical bundle for a family of  LC pairs (Theorem \ref{variation3}).
The prototypes of this kind of argument were already observed in 
\cite{chern,ball} (cf. Example (\ref{orblimit})), where I considered a quasiprojective manifold as a limit of
orbifolds. I hope that we may justify the philosophy in much broader context in future and construct the theory of LC pairs such as finite generation 
of the log canonical rings of LC pairs. The new feature here is that the semipositivity 
is not on the direct image, but on the relative log canonical bundle itself. 
Unlike the KLT case, this is an essential difference, i.e., it seems to be 
difficult to obtain the semipositivity of the direct images of relative pluri log 
canonical systems for a projective family of LC pairs in general.  
    
The organization of this article is as follows.  In Section 1, I explain  
the results in this article.   In Section 2, I prove the existence of canonical K\"{a}hler-Einstein currents on LC pairs of log general type. 
In Section 3,  I construct a Ricci iterations which converges to 
the canonical K\"{a}hler-Einstein current on LC pairs of log general type. 
In Section 4, I decompose the Ricci iteration into a series of dynamical systems of Bergman kernels.
In Section 5, I prove the logarithmic pluri-subhamonicity of canonical K\"{a}hler-Einstein volume forms on a projective family of LC pais. 
In Section 6, I apply the pluri-subhamonic variation property of 
relative canonical measures to deduce the 
weak semistability of the direct images of relative pluri-log-canonial systems  for a projective family of KLT pairs. This is a direct application of 
the results in Section 5 and Viehweg's ingeneous ideas in \cite{v}. \vspace{3mm} \\
\noindent{\bf Acknowledgements} \vspace{2mm}\\  
I would like to express my sincere thanks to Professor Y. Rubinstein who sent me 
his preprint \cite{ru} and drew my attention to Ricci iterations. 
In fact the Ricci iterations used in this article are  variants of the Ricci iterations considered in \cite{ru}. 
I also would like to express my sincere thanks to Professor V. Tosatti who pointed out a serious incompleteness in the previous version.  
  \vspace{5mm}\\
\noindent{\bf Notations}
\begin{itemize}
\item  For a real number $a$, $\lceil a\rceil$ denotes the minimal integer greater than or equal to $a$ and $\lfloor a\rfloor$ denotes the maximal integer smaller than 
or equal to $a$.  
\item Let $X$ be a projective variety and let $D$ be a Weil divsor on $X$.
Let $D = \sum d_{i}D_{i}$ be the irreducible decomposition.  
We set 
\begin{equation}
\lceil D\rceil := \sum \lceil d_{i}\rceil D_{i}
,\lfloor D \rfloor := \sum \lfloor d_{i}\rfloor D_{i}.
\end{equation} 
\item Let $f : X \to Y$ be an algebraic fiber space and let $D$ be 
a $\mathbb{Q}$-divisor on $X$.  Let 
\begin{equation}\label{h-v}
D = D^{h} + D^{v}
\end{equation}
be  the decomposition such that an irreducible component of 
$\mbox{Supp}\, D$ is contained in $\mbox{Supp}\, D^{h}$ if and only if 
it is mapped onto $Y$.  $D^{h}$ is the horizontal part of $D$ 
and $D^{v}$ is the vertical part of $D$.  
\item Let $(X,D)$ be a pair of a normal variety and a $\mathbb{Q}$-divisor 
on $X$.  Suppose that $K_{X} + D$ is $\mathbb{Q}$-Cartier. 
Let $f : Y \to X$ be a log resolution.  Then we have the formula :
\[
K_{Y} = f^{*}(K_{X}+D) + \sum a_{i}E_{i}, 
\]
where $E_{i}$ is a prime divisor and $a_{i}\in \mathbb{Q}$. 
The pair $(X,D)$ is said to be {\bf subKLT}(resp. {\bf subLC}, if $a_{i} > -1$
(resp. $a_{i} \geqq -1$ holds for every $i$. 
$(X,D)$ is said to be {\bf KLT} (resp. {\bf LC}), if $(X,D)$ is subKLT(resp. subLC) and $D$ is effective.
\item Let $X$ be a projective variety and let $\mathcal{L}$ be an invertible
sheaf on $X$. $\mathcal{L}$ is said to be semiample, if there exists a positive
integer $m$ such that $|\mathcal{L}^{\otimes m}|$ is  base point free.  
\item $f : X \to Y$ be a morphism between projective varieties. 
Let $\mathcal{L}$ be an invertible sheaf on $X$.  $\mathcal{L}$ is said to 
be $f$-semiample, if for every $y \in Y$, $\mathcal{L}|f^{-1}(y)$ is 
semiample.  
\item Let  $D,D^{\prime}\in \mbox{Div}(X)\otimes\mathbb{Q}$ be $\mathbb{Q}$-Cartier divisors on a normal projective variety $X$.  
We denote 
\[
D \preceq D^{\prime},
\]
if $D^{\prime} - D$ is effective.
\item Let $L$ be a $\mathbb{Q}$-line bundle on a compact complex manifold $X$, i.e., $L$ is a formal fractional power of a genuine line bundle on $X$.  
A  singular hermitian metric $h$ on $L$ is given by
\[
h = e^{-\varphi}\cdot h_{0},
\]
where $h_{0}$ is a $C^{\infty}$-hermitian metric on $L$ and 
$\varphi\in L^{1}_{loc}(X)$ is an arbitrary function on $X$.
We call $\varphi$ the  weight function of $h$ with respect to $h_{0}$.  We note that $h$ makes sense,
 since a hermitian metric is a real object. 

The curvature current $\Theta_{h}$ of the singular hermitian $\mathbb{Q}$-line
bundle $(L,h)$ is defined by
\[
\Theta_{h} := \Theta_{h_{0}} + \partial\bar{\partial}\varphi,
\]
where $\partial\bar{\partial}\varphi$ is taken in the sense of current.
We define the multiplier ideal sheaf  ${\cal I}(h)$ of $(L,h)$  
by 
\[
\mathcal{I}(h)(U) := \{ f \in \mathcal{O}_{X}(U); \,\, |f|^{2}\,e^{-\varphi}
\in L^{1}_{loc}(U)\},
\]
where $U$ runs open subsets of $X$. 
\item For a closed positive $(1,1)$ current $T$, $T_{abc}$ denotes 
the abosolutely continuous part of $T$. 
\item  For  a Cartier divisor $D$, we denote the 
corresponding line bundle by the same notation. 
Let $D$ be an effective $\mathbb{Q}$-divisor on a smooth projective variety $X$.  Let $a$ be a positive integer such that $aD$ is Cartier. 
We identify $D$ with the formal $a$-th root of the line bundle $aD$.  
We say that $\sigma$ is a multivalued global holomorphic section of $D$ 
with divisor $D$, if  $\sigma$ is the formal $a$-th root of a global holomorphic section of $aD$ with divisor $aD$. And $1/|\sigma|^{2}$ denotes the singular hermitian metric on $D$ defined by 
\[
\frac{1}{|\sigma|^{2}} := \frac{h_{D}}{h_{D}(\sigma,\sigma)}, 
\]
where $h_{D}$ is an arbitrary $C^{\infty}$-hermitian metric on $D$.
\item For a singular hermitian line bundle $(F,h_{F})$ on a compact complex 
manifold $X$ of dimension $n$.   $K(X,K_{X}+F,h_{F})$ denotes (the diagonal part of) the Bergman kernel of $H^{0}(X,{\cal O}_{X}(K_{X} + F)\otimes {\cal I}(h_{F}))$ 
with respect to the $L^{2}$-inner product: 
\begin{equation}\label{inner}
(\sigma ,\sigma^{\prime}) := (\sqrt{-1})^{n^{2}}\int_{X}h_{F}\cdot\sigma\wedge \bar{\sigma}^{\prime}, 
\end{equation}  
i.e., 
\begin{equation}\label{BergmanK}
K(X,K_{X}+F,h_{F}) = \sum_{i=0}^{N}|\sigma_{i}|^{2}, 
\end{equation}
where $\{\sigma_{0},\cdots ,\sigma_{N}\}$ is a complete orthonormal basis 
of $H^{0}(X,{\cal O}_{X}(K_{X} + F)\otimes {\cal I}(h_{F}))$. 
It is clear that $K(X,K_{X}+F,h_{F})$ is independent of the choice of 
the complete orthonormal basis. 
\end{itemize}
     
\subsection{Twisted K\"{a}hler-Einstein equations for ample hermitian adjoint bundles}

First let us consider an easy case. Let $X$ be a smooth projective $n$-fold and let $(L,h_{L})$ be a 
$\mathbb{Q}$-line bundle on $X$ with a $C^{\infty}$-hermitian metric $h_{L}$. Later, we shall consider the case that $h_{L}$ is a singular hermitian metric. 
Here a $\mathbb{Q}$-line bundle means a formal  fractional power 
of a genuine line bundle on $X$.  
We shall assume that $K_{X} + L$ is ample.  
We consider the following equation:
\begin{equation}\label{ricci}
-\mbox{Ric}_{\omega} + \sqrt{-1}\,\Theta_{h_{L}} = \omega, 
\end{equation}
where $\omega$ is a $C^{\infty}$-K\"{a}hler form such that 
the K\"{a}hler class $[\omega]$ is equal to $2\pi c_{1}(K_{X}+L)$ 
in $H^{2}(X,\mathbb{R})$ and $\Theta_{h_{L}}$ denotes the curvature 
form of $h_{L}$.  We call (\ref{ricci}) {\bf the twisted K\"{a}hler-Einstein 
equation} associated with $(L,h_{L})$. 

The equation (\ref{ricci}) is reduced to a complex Monge-Amp\`{e}re 
equation as follows. Let $\omega_{0}$ be a K\"{a}hler form 
on $X$ with $[\omega_{0}] = 2\pi c_{1}(K_{X}+L)$.  
Let $\Omega$ be a $C^{\infty}$-volume form on $X$ such that 
\begin{equation}
\omega_{0} = -\mbox{Ric}\,\Omega + \sqrt{-1}\Theta_{h_{L}} 
\end{equation}
holds. 
Then there exists a $C^{\infty}$-function $u$ on $X$ such that
\begin{equation}
\omega = \omega_{0} + \sqrt{-1}\partial\bar{\partial}u
\end{equation}
and the equation (\ref{ricci}) is equivalent to the Monge-Amp\`{e}re equation:
\begin{equation}\label{MA}
\log\frac{(\omega_{0} +\sqrt{-1}\partial\bar{\partial}u)^{n}}{\Omega}
= u. 
\end{equation}
By the solution of Calabi's conjecture (\cite{a,y}) the  equation (\ref{MA}) has a unique $C^{\infty}$-solution $u$.  Hence (\ref{ricci}) has a $C^{\infty}$-solution $\omega$. 
\begin{definition}\label{TKE}
Let $X$ be a smooth projective variety and let $(L,h_{L})$ be a $C^{\infty}$-hermitian (holomorphic) line bundle on $X$ such that $K_{X}+L$ is ample.
The $C^{\infty}$-K\"{a}hler form  $\omega$ satisfying the equation (\ref{ricci})
is said to be {\bf the twisted K\"{a}hler-Einstein form} 
associated with $(X,(L,h_{L}))$. \fbox{}
\end{definition}
We call $(L,h_{L})$ the {\em boundary} of the pair $(X,(L,h_{L}))$. 
Of course  twisted K\"{a}hler-Einstein forms are not essentially new. 
But this new terminology gives us freedom to think about K\"{a}hler-Einstein-like forms on a smooth projective variety whose first Chern class is not 
necessarily definite.  
 
For example, in this terminology the canonical K\"{a}hler current on the base of an Iitaka fibration is considered to be the twisted K\"{a}hler-Einstein form   associated with the metrized Hodge $\mathbb{Q}$-line bundle (cf. \cite{canonical} or Section \ref{CM} below) as the boundary. 

\subsection{Ricci iterations}
Let $X$ be a compact K\"{a}hler manifold and let $L$ be a $\mathbb{Q}$-line 
bundle on $X$.  Let $h_{L}$ be a (possibly singular) hermitian metric 
on $L$.  
Let $\omega_{0}$ be a $C^{\infty}$-K\"{a}hler form on $X$ and let 
$t$ be a fixed positive number in $(0,1)$. 
In this article,  we shall consider the Ricci iteration of the form :
\begin{equation}\label{iterations1}
-\mbox{Ric}_{\omega_{m}} + t\omega_{m-1} + \sqrt{-1}\Theta_{h_{L}}  = \omega_{m}\end{equation}
for $m \geqq 1$.  We call the term $\sqrt{-1}\Theta_{h_{L}}$ {\em the 
drift term} of the Ricci iteration. 

Suppose that $h_{L}$ is $C^{\infty}$. 
In this case, if the iteration (\ref{iterations1}) has $C^{\infty}$-solution 
$\{\omega_{m}\}$ for $m= 0,\cdots, \ell$, then 
\begin{equation}\label{cohomology}
[\omega_{m}] = \left(\frac{1 - t^{m}}{1 - t}\right)2\pi c_{1}(K_{X}+L) + t^{m}[\omega_{0}]
\end{equation}
holds for every $m = 0,\cdots ,\ell$, where $[\omega_{m}]$ denotes the 
de Rham cohomology class of $\omega_{m}$.   In particular 
 $\{[\omega_{m}]\}$ moves on the segment connecting 
$[\omega_{0}]$ and $(1-t)^{-1}2\pi c_{1}(K_{X}+ D)$.   

Conversely by the solution of Calabi's conjecture(\cite{a,y}), if 
\[
\left(\frac{1 - t^{\ell}}{1 - t}\right)2\pi c_{1}(K_{X}+L) + t^{\ell}[\omega_{0}] 
\]  
is in the K\"{a}hler cone of $X$ for some positive integer  $\ell$, then we see that the sequence 
of $C^{\infty}$-K\"{a}hler forms $\{ \omega_{m}\}_{m=0}^{\ell}$ 
 satisfying (\ref{iterations1}) 
exists.   
Moreover if $c_{1}(K_{X}+L)$ is in the K\"{a}hler cone and $\sqrt{-1}\Theta_{h_{L}} \geqq 0$, then  $\{\omega_{m}\}_{m=0}^{\infty}$ 
exists and we expect that $\{\omega_{m}\}_{m=0}^{\infty}$  converge to a $C^{\infty}$-solution $\omega_{\infty}$  of 
\begin{equation}
- \mbox{Ric}_{\omega_{\infty}} + \sqrt{-1}\Theta_{h_{L}} = \frac{1}{1- t}\,\omega_{\infty}
\end{equation}
as $m$ tends to infinity (cf. Theorem \ref{main1} for a special case). 

In this article, we shall consider the case that 
$K_{X}+L$ is big and $h_{L}$ is a singular hermitian 
metric of the form :
\[
h_{L} = \frac{1}{|\sigma|^{2}},
\]
where $\sigma$ is a multivalued global holomorphic section of $L$ such that 
$(X,(\sigma))$ is KLT, where $(\sigma)$ denotes the effective $\mathbb{Q}$-divisor associated with $\sigma$  on $X$.    We note that in this case the solution is K\"{a}hler-Einstein 
on the locus (contained in $X \backslash \mbox{Supp}\,(\sigma)$), where $\omega_{\infty}$ is a $C^{\infty}$-K\"{a}hler form.

\subsection{Dynamical systems of Bergman kernels}\label{DC}
Here we shall review the construction in \cite{canonical}. 
Let $X$ be a smooth projective variety and let $(L,h_{L})$ be a 
$C^{\infty}$-hermitian $\mathbb{Q}$-line bundle on $X$ with semipositive curvature.  Let $a$ be a positive integer such that $aL$ is a genuine line bundle. 
We assume that $K_{X} + L$ is ample. 
Then as in \cite{canonical}, we may construct 
the twisted K\"{a}hler-Einstein form $\omega$ satisfying (\ref{ricci}) in terms of the dynamical system of Bergman kernels.

Let $A$ be an ample line bundle such that for every pseudo-effective 
singular hermitian line bundle $(F,h_{F})$,  
$\mathcal{O}_{X}(jK_{X} + A + F)\otimes \mathcal{I}(h_{F})$ is globally generated for every $0\leqq j\leqq a-1$. 
Such an ample line bundle $A$ exists by Nadel's vanishing theorem \cite[p.561]{n}.

We shall construct a sequence of Bergman kernels $\{ K_{m}\}_{m=1}^{\infty}$ and a sequence of singular hermitian metrics
$\{ h_{m}\}_{m=1}^{\infty}$  as follows. \vspace{3mm} \\

\noindent We set  
\begin{equation}
K_{1} := \left\{\begin{array}{ll} K(X,K_{X}+A,h_{A}), & \mbox{if}\,\, a > 1 \\ 
& \\ 
& \\
K(X,K_{X}+L+A,h_{L}\cdot h_{A}), & \mbox{if}\,\, a = 1 
\end{array}\right. 
\end{equation}
where $K(X,K_{X}+A,h_{A})$ and $K(X,K_{X}+L+A,h_{L}\cdot h_{A})$ are the Bergman kernels defined as (\ref{BergmanK}). 
Then we set 
\begin{equation}
h_{1} := (K_{1})^{-1}. 
\end{equation}
We continue this process. 
Suppose that we have constructed 
$K_{m}$ and the singular hermitian metric $h_{m}$ on 
$\lfloor\frac{m}{a}\rfloor a(K_{X}+ L) + (m -\lfloor\frac{m}{a}\rfloor a)K_{X}$. We define $K_{m+1}$ by 
\begin{equation}
K_{m+1}:= \left\{\begin{array}{ll} K(X,(m+1)K_{X}+\lfloor \frac{m+1}{a}\rfloor aL+A,h_{m}) & \mbox{if}\,\, m +1 \not{\equiv} 0 \,\,\mbox{mod}\,a \\
 & \\ 
 & \\  
K(X,(m+1)(K_{X}+L)+A,h_{L}^{a}\otimes h_{m}) & \mbox{if}\,\,m +1 \equiv 0 \,\,\mbox{mod}\,a \end{array}\right.  
\end{equation}
and 
\begin{equation}
h_{m+1}:= (K_{m+1})^{-1}.  
\end{equation}
Since $\sqrt{-1}\,\Theta_{h_{m}}$ is a closed positive current by 
definition, by the choice of $A$ and the fact that $K_{X}+L$ is 
pseudoeffective, inductively we can construct the sequences $\{ h_{m}\}_{m\geqq 1}$
and $\{ K_{m}\}_{m \geqq 1}$.  
In fact if we take an AZD $h$ of $K_{X}+L$, then 
\[
\mathcal{O}_{X}\left(A + mK_{X}+\lfloor \frac{m}{a}\rfloor aL\right)\otimes \mathcal{I}(h^{\lfloor\frac{m}{a}\rfloor a})
\hookrightarrow \mathcal{O}_{X}\left(A + mK_{X}+\lfloor \frac{m}{a}\rfloor aL\right)
\otimes\mathcal{I}(h_{m-1}) 
\]
holds for every $m \geqq 1$ and the left hand side is globally 
generated by the choice of $A$. 

This inductive construction is essentially the same one originated by the author in \cite{tu3}. 
The following theorem asserts that the above dynamical system yields the 
twisted K\"{a}hler-Einstein metric on $(X,(L,h_{L}))$. 

\begin{theorem}\label{DS}
 Let $X$ be a smooth projective $n$-fold and let $(L,h_{L})$ be a 
 $C^{\infty}$-hermitian $\mathbb{Q}$-line bundle on $X$ with semipositive 
 curvature $\sqrt{-1}\Theta_{h_{L}}$ such that $K_{X} + L$ is ample.  
 Let  $\{ h_{m}\}_{m \geq 1}$ be the sequence 
of hermitian metrics  as above and let $n$ denote $\dim X$. 
Then 
\begin{equation}
h_{\infty} := \liminf_{m\rightarrow\infty} \sqrt[m]{(m!)^{n}\cdot h_{m}}
\end{equation}
is a $C^{\infty}$-hermitian metric on $K_{X}+L$ such that 
\begin{equation}
\omega = \sqrt{-1}\,\Theta_{h_{\infty}}
\end{equation}
is the twisted K\"{a}hler-Einstein form (Definition \ref{TKE}).  

In particular $\omega = \sqrt{-1}\Theta_{h_{\infty}}$ (in fact $h_{\infty}$) is unique and  is independent of the choice of $A$ and $h_{A}$. 
\fbox{}
\end{theorem} 
The proof of Theorem \ref{DS} is completely the same as the one of 
\cite[Theorem 1.8]{canonical}.  Hence we omit it (The only difference is that here $h_{L}$ is $C^{\infty}$ while $h_{L}$ is singular in \cite[Theorem 1.8]{canonical}.). 
Theorem \ref{DS} naturally leads us to the following semipositivity 
theorem. 

\begin{theorem}\label{easysemipositive}
Let $f : X \to Y$ be an algebraic fiber space and let 
$L$ be a  $\mathbb{Q}$-line bundle on $X$ with 
the $C^{\infty}$-hermitian metric $h_{L}$.  
Suppose that $\sqrt{-1}\Theta_{h_{L}}$ is semipositive on $X$ and 
$K_{X/Y} + L$ is relatively ample over $Y_{0}$.  
Let $Y_{0}$ be the maximal Zariski open subset of $Y$ such that $f$ is smooth  
over $Y_{0}$.   
For $y\in Y_{0}$, let $\omega_{y}$ be the K\"{a}hler form 
on $X_{y} := f^{-1}(y)$ satisfying the equation: 
\begin{equation}
-\mbox{\em Ric}_{\omega_{y}} + \sqrt{-1}\Theta_{h_{L}}|X_{y} = \omega_{y}. 
\end{equation} 
Let $n$ be the relative dimension of $f : X \to Y$. 
Then the singular hermitian metric $h^{\circ}$ on $K_{X/Y} + L|f^{-1}(Y_{0})$ defined by 
\begin{equation}
h^{\circ}|X_{y} := \left(\frac{1}{n!}\,\omega_{y}^{n}\right)^{-1}\cdot h_{L} \hspace{10mm} (y\in Y_{0})
\end{equation}
extends to a singular hermitian metric $h$ on $K_{X/Y}+L$ and has semipositive curvature in the sense of current everywhere on $X$. 
\fbox{}
\end{theorem}
\begin{remark}
$\omega$ is $C^{\infty}$ over the complement of the discriminant locus
of $f$.  This is an easy consequence of the implicit function theorem. 
\fbox{}
\end{remark}
\noindent The proof of Theorem \ref{easysemipositive} is completely the same 
as the one of \cite[Theorem 4.1]{canonical}.  Hence we omit it. 
In fact Theorem \ref{DS} and the logarithmic pluri-subharmonicity of 
Bergman kernels (\cite{b3,b-p,KE}) implies Theorem \ref{easysemipositive}.   

\subsection{Canonical K\"{a}hler-Einstein currents on KLT pairs}

Now we shall consider the twisted K\"{a}hler-Einstein form for 
some singular hermitian $\mathbb{Q}$-line bundle $(L,h_{L})$. 

Let $X$ be a smooth projective $n$-fold and let let $D$ be an effective 
$\mathbb{Q}$-divisor on $X$ such that $(X,D)$ is KLT.  
We assume that $(X,D)$ is of log general type, i.e., $K_{X}+D$ is big.  
Let us consider the equation :
\begin{equation}\label{KE}
-\mbox{Ric}_{\omega_{K}} + 2\pi [D] = \omega_{K}, 
\end{equation}  
where $[D]$ denotes the closed positive current associated with 
the effective $\mathbb{Q}$-divisor $D$. 
Let $\sigma_{D}$ be a multivalued global holomorphic section of $D$ with divisor $D$ and let $h_{\sigma_{D}}$ be the singular hermitian metric on $D$
defined by 
\begin{equation}\label{sigmad}
h_{\sigma_{D}} := \frac{1}{|\sigma_{D}|^{2}}. 
\end{equation}
We construct a solution $\omega_{K}$ of (\ref{KE}) which satisfies the followings:
\begin{enumerate}
\item[(1)] $\omega_{K}$ is a closed positive current on $X$, 
\item[(2)] There exists a nonempty Zariski open subset $U$ of $X$ such that 
$\omega_{K}|U$ is $C^{\infty}$-K\"{a}hler form,
\item[(3)] We define the singular hermitian metric $h_{K}$ on $K_{X}+D$ 
by  
\[
h_{K}:= \left(\frac{1}{n!}\,\omega_{K,abc}^{n}\right)^{-1}\cdot h_{\sigma_{D}}.
\] 
Then $h_{K}$ is an AZD (cf. \cite{tu,tu2}) of $K_{X}+ D$, i.e., 
\begin{enumerate}
\item $\sqrt{-1}\,\Theta_{h_{K}}$ is a closed semipositive current on $X$, 
\item 
For every positive integer $m$ such that $m(K_{X}+D)$ is Cartier 
\[
H^{0}(X,\mathcal{O}_{X}(m(K_{X}+D))\otimes\mathcal{I}(h_{K}^{m}))
\simeq H^{0}(X,\mathcal{O}_{X}(m(K_{X}+D)))
\]
holds.
\end{enumerate} 
\end{enumerate} 
We note that if such $\omega_{K}$ exists, then $\omega_{K}$ is 
K\"{a}hler-Einstein on $U$, i.e., 
$-\mbox{Ric}_{\omega_{K}} = \omega_{K}$ holds on $U$. 

\begin{definition}\label{GKE}
Let $(X,D)$ be a KLT pair of log general type. A closed positive 
current $\omega_{K}$ satisfying the above properties: (1),(2),(3) is said to be 
{\bf the canonical K\"{a}hler-Einstein current} on $(X,D)$. \fbox{}  
\end{definition}   
The following theorem asserts that the canonical K\"{a}hler-Einstein current 
always exists on a KLT pair $(X,D)$ of log general type. 
\begin{theorem}\label{main}
Let $X$ be a smooth projective variety and let $D$ be an effective divisor 
such that $(X,D)$ is a KLT pair of log general type.  
Then there exists a unique canonical K\"{a}hler-Einstein current 
$\omega_{K}$ on $(X,D)$. \fbox{}
\end{theorem}
Later we shall prove a stronger uniqueness of the canonical K\"{a}hler-Einstein  currents on KLT pairs (Theorem \ref{unique}).
The proof of Theorem \ref{main} here is very similar to the case that $D = 0$ in \cite{canonical}. 
The proof is done by solving a Monge-Amp\`{e}re equation as the solution of 
Calabi's conjecture in \cite{a,y}.
The  Monge-Amp\`{e}re equation is 
quite similar to the one in \cite[p.409,Theorem 8]{y}. 
But to prove Theorem \ref{main}, we cannot apply \cite[p.409, Theorem 8]{y}
directly, since  the coefficients of $D$ may exceed $1/n (n = \dim X)$. 
We overcome this difficulty by introducing an orbifold structure on 
$(X,D)$, i.e., we resolve the singularities in terms of local 
cyclic branched coverings. 

\subsection{Canonical measures on  KLT pairs of nonnegative Kodaira dimension}
\label{CM}
Similar twisted K\"{a}hler-Einstein currents naturally appear in 
algebraic geometry of KLT pairs.

Let  $(X,D)$ be a KLT pair  of  nonnegative Kodaira dimension, i.e., 
$|m!(K_{X}+D)|\neq \emptyset$ for every sufficiently large $m$.

Let $f : X -\cdots\rightarrow Y$ be the Iitaka fibration associated with 
the log canonical divisor $K_{X} +D$.
By replacing $X$ and $Y$ by suitable modifications, we may assume the followings:
\begin{enumerate}
\item[(1)] $X$,$Y$ are smooth and $f$ is a morphism with connected fibers. 
\item[(2)] $\mbox{Supp}\,D$ is a divisor with normal crossings. 
\item[(3)] There exists an effective divisor $\Sigma$ on $Y$ such that 
$f$ is smooth over $Y - \Sigma$, $\mbox{Supp}\, D^{h}$ is relatively normal crossings
over $Y - \Sigma$ and $f(D^{v})\subset \Sigma$, where $D^{h}, D^{v}$ denote the horizontal and the vertical component of $D$ respectively (cf. (\ref{h-v})). 
\item[(4)] There exists a positive integer $m_{0}$ such that for every $m \geqq m_{0}$, $m!(K_{X}+D)$ is 
Cartier and $\left(f_{*}{\cal O}_{X}(m!(K_{X}+D))\right)^{**}$ is a line bundle on $Y$, where $**$ denotes the double dual.  
\end{enumerate}
We note that adding effective exceptional $\mathbb{Q}$-divisors does not 
change the log canonical ring.   
Such a modification exists by \cite[p.175, Proposition 4.2]{f-m}. 
We define the $\mathbb{Q}$-line bundle $L_{X/Y,D}$ on $Y$ by 
\begin{equation}
L_{X/Y,D} = \frac{1}{m_{0}!}\left(f_{*}{\cal O}_{X}(m_{0}!(K_{X}+D))\right)^{**}. 
\end{equation}
$L_{X/Y,D}$ is independent of the choice of $m_{0}$. 
Similarly as before we may define the singular hermitian metric 
$h_{L_{X/Y,D}}$ on $L_{X/Y,D}$ by 
\begin{equation}
h_{L_{X/Y,D}}^{m_{0}!}(\sigma,\sigma)(y):= \left(\int_{X_{y}}|\sigma|^{\frac{2}{m_{0}!}}\right)^{m_{0}!},
\end{equation}
where $y \in Y - \Sigma, X_{y}:= f^{-1}(y)$ and $\sigma
\in m_{0}!L_{X/Y,D,y}$. 
We call the singular hermitian $\mathbb{Q}$-line bundle 
$(L_{X/Y,D},h_{L_{X/Y,D}})$ {\bf the metrized Hodge $\mathbb{Q}$-line bundle} of the Iitaka fibration 
$f : X \to Y$ associated with the KLT pair $(X,D)$. 
We note that since $(X,D)$ is KLT, $h_{L_{X/Y,D}}$ is well defined. 
By the same strategy as in the proof of Theorem \ref{main}
and \cite[Theorem 1.6]{canonical}, we have  
the following theorem : 
\begin{theorem}\label{main2}
In the above notations, there exists a unique singular hermitian metric 
on $h_{K}$ on $K_{Y}+L_{X/Y,D}$ and a nonempty Zariski open subset 
$U$ of $Y$ such that 
\begin{enumerate}
\item[(1)] $h_{K}$ is an AZD of $K_{Y} + L_{X/Y,D}$, 
\item[(2)] $f^{*}h_{K}$ is an AZD of $K_{X}+D$,  
\item[(3)] $h_{K}$ is $C^{\infty}$ on $U$, 
\item[(4)] $\omega_{Y} = \sqrt{-1}\,\Theta_{h_{K}}$ is a K\"{a}hler form 
on $U$,
\item[(5)] $-\mbox{\em Ric}_{\omega_{Y}} + \sqrt{-1}\,\Theta_{L_{X/Y,D}} = \omega_{Y}$
holds on $U$. \fbox{} 
\end{enumerate}
\end{theorem}

\noindent We may construct $\omega_{Y}$ in Theorem \ref{main2} 
in terms of a family of  dynamical systems of Bergman kernels as
Theorem \ref{variation} below   
and we obtain the same semipositivity. 
We define {\bf the canonical measure $d\mu_{can}$ of $(X,D)$} by 
\begin{equation}\label{canmeas}
d\mu_{can} := \frac{1}{n!}\,f^{*}\left(\omega_{Y,abc}^{n}\cdot h_{L_{X/Y,D}}^{-1}\right), 
\end{equation} 
where $n = \dim Y$. 
Then $d\mu_{can}$ is considered to be a singular volume form on $X$. 

\subsection{Variation of canonical K\"{a}hler-Einstein currents}

One of the  important properties of canonical K\"{a}hler-Einstein currents 
constructed in Theorem \ref{main} is the following pluri-subharmonic variation property.   

\begin{theorem}\label{variation}
Let $f : X \to Y$ be an algebraic fiber space and let $D$ be an effective 
$\mathbb{Q}$-divisor on $X$.
Suppose that there exists a nonempty Zariski 
open subset $Y_{0}$ of $Y$ such that
\begin{enumerate}
\item[(1)] $f$ is smooth over $Y_{0}$, 
\item[(2)]  
$(X_{y},D_{y}) (X_{y} := f^{-1}(y),D_{y} := D\cap X_{y})$ is a KLT pair
of log general type. 
\end{enumerate}
Let $\omega_{K,y}$ be the canonical K\"{a}hler-Einstein current  on 
$(X_{y},D_{y}) (y\in Y_{0})$ constructed as in Theorem \ref{main}.
Let $n$ be the relative dimension of $f : X \to Y$. 

Then the singular hermitian metric $h_{K}$ on $K_{X/Y}+D|f^{-1}(Y_{0})$ defined by 
\[
h^{\circ}_{K}|X_{y} := \left(\frac{1}{n!}\,(\omega_{K,y})_{abc}^{n}\right)^{-1} \!\!\!\!\cdot h_{\sigma_{D}}|X_{y}\,\,\,\,\,\,(y\in Y_{0}) 
\]
extends to a singular hermitian metric $h_{K}$ on $K_{X/Y}+D$ and has  semipositive curvature  in the sense of current everywhere on $X$. \fbox{}
\end{theorem}
This theorem is a natural generalization of the result in 
\cite{KE}, where we deal with the case that $D = 0$.
We note that Theorem \ref{variation} is closely related to the 
semipositivity of the direct image of a relative pluri-log-canonical systems of a family of KLT pairs \cite[p.175, Theorem 1.2]{ka3}.  
In fact Theorem \ref{variation} is considered to be a generalization of it.   
More generally we have the following theorem.  

\begin{theorem}\label{variation2}
Let $f : X \to Y$ be an algebraic fiber space and let $D$ be an effective 
$\mathbb{Q}$-divisor on $X$.
Suppose that there exists a nonempty Zariski 
open subset $Y_{0}$ of $Y$ such that 
\begin{enumerate}
\item[(1)] $f$ is smooth over $Y_{0}$, 
\item[(2)] $(X_{y},D_{y}) (X_{s} := f^{-1}(y),D_{y} := D\cap X_{y})$ is a KLT pair of nonnegative Kodaira dimension.
\end{enumerate} 
Let $d\mu_{can,y}$ be the canonical measure  on 
$(X_{y},D_{y}) (y\in Y_{0})$ as in Section \ref{CM}. 
Then the singular hermitian metric $h^{\circ}_{K}$ on $K_{X/Y}+D|f^{-1}(Y_{0})$ defined by 
\[
h^{\circ}_{K}|X_{y} := d\mu_{can,y}^{-1}\cdot h_{\sigma_{D}}|X_{y} (y\in Y_{0}) 
\]
extends to a singular hermitian metric on $K_{X/Y} +D$  and
has  semipositive curvature in the sense of current everywhere on $X$. \fbox{}
\end{theorem}
This is a natural generalization of \cite[Theorem 4.1]{canonical}. 

To prove Theorem \ref{variation},  we use 
the dynamical construction of canonical K\"{a}hler-Einstein currents
as in \cite{canonical}.  But there is a major difference between 
the dyanamical construction here and the one in \cite{canonical}.

The most naive way to construct canonical K\"{a}hler current 
by a dynamical system of Bergman kernels is to replace the hermitian line bundle $(L,h_{L})$ in Section \ref{DC}  by 
$(D,|\sigma_{D}|^{-2})$, where $\sigma_{D}$ is a multivalued 
holomorphic section of the $\mathbb{Q}$-line bundle $D$ with divisor 
$D$. Let $a$ denote the minimal positive integer such that $aD$ is 
Cartier.  Then in the above naive dynamical construction (cf. Section \ref{DC}) we tensorize $(\mathcal{O}_{X}(aD),|\sigma_{D}|^{-2a})$ in every 
$a$-steps. But unfortunately this construction does not yield the 
desired canonical K\"{a}hler-Einstein current, since $|\sigma_{D}|^{-2a}$
has too large singularities in general.  Moreover the dynamical system 
itself is not well defined in general.  
Actually this construction  works only when $D$ is $\mathbb{Q}$-linear equivalent 
to a Cartier divisor essentially.   

To overcome this difficulty we consider a Ricci iteration and 
a family of dynamical systems of Bergman kernels associated with 
the Ricci iteration instead of a single 
dynamical system. 

\subsection{Semipositivity of the relative log canonical bundle on  
a family of LC pairs}

Let $X$ be a smooth projective variety and let $D$ be an effective 
$\mathbb{Q}$-divisor on $X$ such that $(X,D)$ is a LC pair. 
Then for every rational number $t \in [0,1)$, we see that 
$(X,tD)$ is a KLT pair.   
Suppose that $(X,D)$ is of log general type. 
Then for a sufficiently small rational number $0< \epsilon \ll  1$,
$(X,(1-\epsilon)D)$ is a KLT pair of log general type.  
Let $\omega_{\epsilon}$ be the canonical K\"{a}hler-Einstein current 
 on $(X,(1-\epsilon)D)$ for $0 < \epsilon \ll  1$ (cf. Theorem \ref{main}). 
 Then the K\"{a}hler-Einstein volume form $dV_{\epsilon}$ of 
 $\omega_{\epsilon}$ is monotone increasing 
 as $\epsilon$ tends to $0$ and 
 converges to a singular volume form  on $X$ (see Lemma \ref{monotonicity2} below). 
Using this simple fact, Theorem \ref{variation} can be generalized to the case of LC pairs.

\begin{theorem}\label{variation3}
Let $f : X \to Y$ be an algebraic fiber space and let $D$ be an effective 
$\mathbb{Q}$-divisor on $X$.
Suppose  that there exists a nonempty Zariski 
open subset $Y_{0}$ of $Y$ such that
\begin{enumerate}
\item[(1)] $f$ is smooth over $Y_{0}$,
\item[(2)]
For every $y\in Y_{0}$, $(X_{y},D_{y}) (X_{s} := f^{-1}(y),D_{y} := D\cap X_{y})$ is a LC pair of log general type.
\end{enumerate}
Let $\omega_{K,y}$ be the canonical K\"{a}hler-Einstein current  on 
$(X_{y},D_{y}) (y\in Y_{0})$ constructed as in Theorem \ref{limitvol}
below. 
Let $n$ be the relative dimension of $f : X \to Y$.

Then the singular hermitian metric $h_{K}$ on $K_{X/Y}+D|f^{-1}(Y_{0})$ defined by 
\[
h^{\circ}_{K}|X_{y} := \left(\frac{1}{n!}\,(\omega_{K,y})_{abc}^{n}\right)^{-1}\!\!\!\!\!\cdot h_{\sigma_{D}}|X_{y}\,\,\,\,\,\,\, (y\in Y_{0}) 
\]
extends to a singular hermitian metric $h_{K}$ on $K_{X/Y}+D$ and has  semipositive curvature in the sense of current everywhere on $X$. \fbox{}
\end{theorem}
For  families of LC pairs not necessarily of log general type, 
we have the following corollary (see \cite{b-p} for a special case). 

\begin{corollary}\label{pe}
Let $f : X \to Y$ be an algebraic fiber space and let $D$ be an effective 
$\mathbb{Q}$-divisor on $X$.
Suppose  that there exists a nonempty Zariski 
open subset $Y_{0}$ of $Y$ such that
\begin{enumerate}
\item[(1)] $f$ is smooth over $Y_{0}$, 
\item[(2)]
For every $y\in Y_{0}$, $(X_{y},D_{y}) (X_{s} := f^{-1}(y),D_{y} := D\cap X_{y})$ is a LC pair such that $K_{X_{y}} + D_{y}$ is pseudo-effective.
\end{enumerate}
Then $K_{X/Y} + D$ is pseudo-effective.  \fbox{}
\end{corollary}
{\em Proof of Corollary \ref{pe}}.
Let $A$ be an ample effective divisor on $X$. 
Then for every positive rational number $\varepsilon$, by Theorem \ref{variation3},
we see that $K_{X/Y}+ D + \varepsilon A$ is pseudo-effective. 
Letting $\varepsilon\downarrow 0$, we complete the proof of Corollary \ref{pe}.
q.e.d. \vspace{3mm} \\ 

\noindent The following corollary is a typical special case of Theorem \ref{variation3}. 
\begin{corollary}\label{complete}
Let $f : X \to Y$ be an algebraic fiber space and let $D$ be an 
effective divisor on $X$ with normal crossings.
Let $n$ denote the relative dimension of $f : X \to Y$.
Suppose that there exists a nonempty Zariski open subset $Y_{0}$ 
of $Y$ such that 
\begin{enumerate} 
\item[(1)] $K_{X/Y} + D$ is relatively ample over $Y_{0}$,
\item[(2)] $f$ is smooth over $Y_{0}$,
\item[(3)] For every $y\in Y_{0}$, $D_{y} = X_{y}\cap D$ (the scheme 
theoretic intersection) is a divisor with normal crossings in 
$X_{y}:= f^{-1}(y)$.  
\end{enumerate}
For every $y\in Y_{0}$ let $\omega_{E,y}$ be the complete 
K\"{a}hler-Einstein form on $X_{y}\backslash D_{y}$ 
such that $-\mbox{\em Ric}_{\omega_{E,y}} = \omega_{E,y}$ 
constructed as in \cite{ko}.   
We define the metric $h_{K}^{\circ}$ on $K_{X/Y} + D|f^{-1}(Y_{0})$ 
by 
\[
h_{K}^{\circ}|X_{y} := \left(\frac{1}{n!}\,\omega_{E,y}^{n}\right)^{-1}\!\!\cdot h_{\sigma_{D}}|X_{y}. 
\]
Then $h_{K}^{\circ}$ extends to a singular hermitian metric 
$h_{K}$ on $K_{X/Y} + D$ with semipositive curvature in the sense of 
current.  \fbox{}
\end{corollary}
\begin{remark}
In Corollary \ref{complete}, by the implicit function theorem, we see 
that $h_{K}$ is $C^{\infty}$ on a nonempty Zariski open subset of $X$. \fbox{}
\end{remark}
For a family of KLT pairs of not necessarily of log general type, we 
have the following useful theorem. 
\begin{theorem}\label{variation4}
Let $f : X \to Y$ be an algebraic fiber space and let $D$ be an effective 
$\mathbb{Q}$-divisor on $X$.
Suppose  that there exists a nonempty Zariski 
open subset $Y_{0}$ of $Y$ such that
\begin{enumerate}
\item[(1)] $f$ is smooth over $Y_{0}$, 
\item[(2)]
For every $y\in Y_{0}$, $(X_{y},D_{y}) (X_{s} := f^{-1}(y),D_{y} := D\cap X_{y})$ is a KLT pair of nonnegative Kodaira dimension. 
\end{enumerate}
Let $d\mu_{can,X/Y}$ be the relative canonical measure defined by
\[
d\mu_{can,X/Y}|X_{y} := d\mu_{can,y}\,\,\,\,\,\,\, (y\in Y_{0})
\]
where $d\mu_{can,y}$ denotes the canonical measure on  
$(X_{y},D_{y}) (y\in Y_{0})$ constructed as in Theorem \ref{main2}.
Then the singular hermitian metric 
\[
h^{\circ}_{K}|X_{y} := d\mu_{can,y}^{-1}\cdot h_{\sigma_{D}}|X_{y} \,\,\,\,\,\,\, (y\in Y_{0}) 
\]
on $K_{X/Y} + D|f^{-1}(Y_{0})$ extends to a singular hermitian metric $h_{K}$ on $K_{X/Y}+D$ and has  semipositive curvature in the sense of current everywhere on $X$. \fbox{}
\end{theorem}
There are numerous applications of the results in this article. 
They will be discussed in the subsequent paper(\cite{global}). 

\section{Existence of canonical K\"{a}hler-Einstein currents}\label{existence}

In this section we shall prove Theorem \ref{main}.  
The proof is not essentially new except the use of orbifold structures 
and the uniqueness.   This simple technique overcomes the difficulty arising from 
the poles of the righthand side of the Monge-Amp\`{e}re equation: (\ref{MP}) below, i.e., 
we can eliminate the poles of the Monge-Amp\`{e}re equation in terms of 
 local cyclic coverings.  
\subsection{Existence of canonical K\"{a}hler-Einstein currents on KLT pairs
of log general type}\label{ex1}

In this subsection, we shall prove Theorem \ref{main}. 
Let $X$ be a smooth projective $n$-fold and 
let $D$ be an effective $\mathbb{Q}$-divisor on $X$ such that 
$(X,D)$ is KLT and  $K_{X} + D$ is big. 
Then by \cite{b-c-h-m} the canonical ring:
\begin{equation}
R(X,K_{X}+D) = \oplus_{m=0}^{\infty}\Gamma (X,\mathcal{O}_{X}(\lfloor m(K_{X}+D)\rfloor))
\end{equation}
is finitely generated.  Hence there exists a log resolution
\begin{equation}
\mu : Y \to X
\end{equation}
of $(X,D)$ which satisfies  the followings:  
\begin{enumerate}
\item[(1)] If we write $K_{Y} = \mu^{*}(K_{X}+D) + \sum a_{i}E_{i}$, 
then $a_{i} > -1$ holds for every $i$, where $\{ E_{i}\}$ are 
prime divisors.  
\item[(2)] There exists a decomposition: 
\begin{equation}\label{ZD}
\mu^{*}(K_{X}+D) = P + N (P,N \in \mbox{Div}(Y)\otimes \mathbb{Q})
\end{equation}
into  $\mathbb{Q}$-divisors $P$ and $N$  
such that $P$ is semiample, $N$ is effective and  
\begin{equation}
H^{0}(Y,\mathcal{O}_{Y}(\mu^{*}(amP))
\simeq H^{0}(Y,\mathcal{O}(am\mu^{*}(K_{X} + D)))
\end{equation}
holds for every $m\geqq 0$, where $a$ is the minimal positive integer 
such that $aD, aP \in \mbox{Div}(X)$ hold. 
\end{enumerate} 
We call the decomposition (\ref{ZD}) a Zariski decomposition of 
$\mu^{*}(K_{X} + D)$.  

We note that adding effective exceptional $\mathbb{Q}$-divisors does not 
change the log canonical ring.  
We set  $I := \{ i|\, a_{i} < 0\}$. 
Then  replacing $X$ by $Y$ and $D$ by
\begin{equation}\label{DY}
D_{Y} := \sum_{i\in I} (-a_{i})E_{i}, 
\end{equation}
we obtain the new KLT pair $(Y,D_{Y})$ such that 
\begin{enumerate}
\item[(a)] $R(Y,K_{Y} + D_{Y}) \simeq R(X,K_{X} +D)$, 
\item[(b)] There exists a Zariski decomposition:  
$K_{Y} + D_{Y} = P + (N + \sum_{i\not{\in} I} a_{i}E_{i})$.
\end{enumerate}
Hence we may assume that $\mbox{Supp}\, D$ and $\mbox{Supp}\, N$ are divisors with normal crossings and the decomposition 
$K_{X} + D = P + N$ holds on $X$ from the beginning. 
In fact if we construct a canonical K\"{a}hler-Einstein current $\omega_{K}$ on $Y$, 
then the push forward $f_{*}\omega_{K}$ is a canonical K\"{a}hler-Einstein 
current on $(X,D)$.   

Let $\sigma_{D}$ be a multivalued holomorphic section of 
the $\mathbb{Q}$-line bundle $D$ with divisor $D$. 
Then $|\sigma_{D}|^{-2}$ is a singular hermitian metric on $L$. 
Similarly we consider the singular hermitian metric 
$|\sigma_{N}|^{-2}$ on the $\mathbb{Q}$-line bundle $N$.  
Let $h_{D}$,$h_{N}$ be $C^{\infty}$-hermitian metrics on $D$ and $N$ 
respectively.  
Then since $P$ is semiample, there exists a $C^{\infty}$-volume form 
$\Omega$ on $X$ such that 
\begin{equation}
h_{P} : = \Omega^{-1}\cdot h_{D}\cdot h_{N}^{-1}
\end{equation}  
is a $C^{\infty}$-hermitian metric on $P$ such that the curvature 
$\sqrt{-1}\Theta_{h_{P}}$ is the positive multiple of the pull back of 
the Fubini-Study form by the base point free linear system 
$|m_{0}!P|$ for $m_{0}\gg 1$. 
Let $D = \sum d_{i}D_{i}$ be the irreducible decomposition of $D$ and 
let $\sigma_{i}$ be a nontrivial global section of $\mathcal{O}_{X}(D_{i})$
with divisor $D_{i}$ and let $\parallel\sigma_{i}\parallel$ be the 
hermitian norm of $\sigma_{i}$ with respect to a $C^{\infty}$-hermitian metric 
on $\mathcal{O}_{X}(D_{i})$.  By multiplying a small positive constant to $\sigma_{i}$, we may and do assume that 
$\parallel\sigma_{i}\parallel < 1$  holds for every $i$ on $X$. 
Let $b$ be a positive integer such that 
\begin{equation}\label{choice}
d_{i} < \frac{b -1}{b}
\end{equation}
holds for every $i$ (Since $(X,D)$ is KLT, $d_{i} < 1$ holds for every $i$.). 
We set 
\begin{equation}
\omega_{P} := \sqrt{-1}\,\Theta_{h_{P}}. 
\end{equation}
We consider the Monge-Amp\`{e}re equation 
\begin{equation}\label{MP}
(\omega_{P} + \sqrt{-1}\partial\bar{\partial}u)^{n}= \frac{\parallel\sigma_{N}\parallel_{h_{N}}^{2}}{\parallel\sigma_{D}\parallel_{h_{D}}^{2}}\cdot \Omega\cdot e^{u}
\end{equation}
on $X$.   
To solve (\ref{MP}), we consider the perturbation of the equation (\ref{MP}) 
and construct the solution  as the limit of the 
solutions of the perturbed equations. 
We may assume that there exists an effective exceptional $\mathbb{Q}$-divisor $E$ with respect to the morphism $\Phi_{|m!P|}: X \to \mathbb{P}^{\nu}$
for a sufficiently large $m$
 such that 
\begin{enumerate}
\item[(1)] $P - \delta E$ is ample for every  
$\delta \in (0,1]$, 
\item[(2)] If we define
\begin{equation}\label{ud}
U := \{ x\in X|\mbox{$|m!P|$ defines an embedding for every $m \gg 1$ on a neighbourhood of $x$}\},
\end{equation}
then 
\begin{equation}\label{ed}
X \backslash \mbox{Supp}\, E  = U
\end{equation}
holds. 
\end{enumerate}
The existence of such $E$ follows from the definition (\ref{ud}) of $U$ and the trivial fact that 
for the composition of any  successive blowing ups 
\[
\varpi : \tilde{\mathbb{P}^{\nu}}\to \mathbb{P}^{\nu}
\]
of  $\mathbb{P}^{\nu}$ with smooth centres, 
there exists an effective $\mathbb{Q}$-divisor $B$  supported on the exceptional  divisors of $\varpi$ such that $\varpi^{*}\mathcal{O}(1) - B$ 
is ample.  Hence by taking a suitable successive blowing ups 
with smooth centres over $X \backslash U$, 
we may assume the existence of such an effective $\mathbb{Q}$-divisor 
$E$. Then there exists a $C^{\infty}$-hermitian metric $h_{E}$ on $E$ 
such that $h_{P}\cdot h_{E}^{-1}$ is a metric with strictly positive 
curvature on $X$.
For every $0 < \delta \leqq 1$, we define the orbifold K\"{a}hler form $\omega_{P,\delta}$ on $X$ by 
\begin{equation}\label{Pdelta}
\omega_{P,\delta}:= (1-\delta)\omega_{P} + \delta\left(\omega_{P}- \sqrt{-1}\,\Theta_{h_{E}}
 +\varepsilon \sum_{i}\sqrt{-1}\partial\bar{\partial}
\log (1 - \parallel\sigma_{i}\parallel^{\frac{2}{b}})\right), 
\end{equation}  
where $\varepsilon$ is a fixed sufficiently small positive number 
so that $\omega_{P,\delta}$ is an orbifold K\"{a}hler form on 
$X$ branching along $D$ with order $b$ for every $\delta\in (0,1]$.

More precisely, $\omega_{P,\delta} (\delta\in (0,1])$ is an orbifold
K\"{a}hler form in the following sense. 
Let $(V,(z_{1},\cdots,z_{n}))$ be a local coordinate on $X$ 
such that $V$ is biholomorphic to the unit open polydisk $\Delta^{n}$ 
in $\mathbb{C}^{n}$ with centre $O$ and  
\[
V \cap D = \{p\in V| z_{1}(p)\cdots z_{k}(p) = 0\}
\]
holds for some $k$ (We have assumed that $\mbox{Supp}\, D$ is a divisor with normal crossings on $X$).  
Let 
\begin{equation}\label{piu}
\pi_{V} : \Delta^{n} \to V
\end{equation}
be the morphism defined by 
\begin{equation}\label{piu2}
\pi_{V}(t_{1},\cdots ,t_{n}) = (t_{1}^{b},\cdots ,t_{k}^{b},t_{k+1},\cdots,t_{n}). 
\end{equation}
If we take a positive number $\varepsilon$ sufficiently small, then  $\pi_{V}^{*}(\omega_{P,\delta}|V) (\delta\in (0,1])$  
is a $C^{\infty}$-K\"{a}hler form on $\Delta^{n}$. 
In this sense $\omega_{P,\delta}(\delta \in (0,1])$ is an orbifold K\"{a}hler 
form on $(X,D)$. 
  
Now we consider the perturbed equation: 
\begin{equation}\label{MD}
(\omega_{P,\delta} + \sqrt{-1}\partial\bar{\partial}u_{\delta})^{n}= \frac{\parallel\sigma_{N}\parallel_{h_{N}}^{2}\left((\prod_{i}(1 - \parallel\sigma_{i}\parallel^{\frac{2}{b}}))^{-\varepsilon}\parallel\sigma_{E}\parallel_{h_{E}}^{2}\right)^{\delta}}{\parallel\sigma_{D}\parallel^{2}_{h_{D}}}\cdot \Omega\cdot e^{u_{\delta}} 
\end{equation}
for $\delta\in (0,1]$.
Then for $\pi_{V}:\Delta^{n}\to V$ as above, pulling back (\ref{MD}) by $\pi_{V}$, we obtain 
\begin{equation}\label{MDV}
(\pi_{V}^{*}\omega_{P,\delta} + \sqrt{-1}\partial\bar{\partial}\pi_{V}^{*}u_{\delta})^{n}= \pi_{V}^{*}\left(\frac{\parallel\sigma_{N}\parallel_{h_{N}}^{2}\left((\prod_{i}(1 - \parallel\sigma_{i}\parallel^{\frac{2}{b}}))^{-\varepsilon}\parallel\sigma_{E}\parallel_{h_{E}}^{2}\right)^{\delta}}{\parallel\sigma_{D}\parallel_{h_{D}}^{2}}\cdot \Omega\right)\cdot e^{\pi_{V}^{*}u_{\delta}}. 
\end{equation}
In (\ref{MDV}), 
\[
\pi_{V}^{*}\left(\frac{\parallel\sigma_{N}\parallel_{h_{N}}^{2}\left((\prod_{i}(1 - \parallel\sigma_{i}\parallel^{\frac{2}{b}}))^{-\varepsilon}\parallel\sigma_{E}\parallel_{h_{E}}^{2}\right)^{\delta}}{\parallel\sigma_{D}\parallel_{h_{D}}^{2}}\cdot \Omega\right)
\]
degenerates along $\{(t_{1},\cdots ,t_{n})\in \Delta^{n}| t_{1}\cdots t_{k} = 0\}$ by the choice of $b$ (cf. (\ref{choice})). 
Hence the equation (\ref{MD}) is considered to be a complex Monge-Amp\`{e}re
equation with degeneracy along $\mbox{Supp}\, D$ as an equation on  the orbifold branching along $\mbox{Supp}\,D$ with order $b$.  Hence by \cite[p.387,Theorem 6]{y}, there exists 
a solution $u_{\delta}$ of (\ref{MD}) on $X$ such that 
\begin{enumerate}
\item $u_{\delta}$ is $C^{\infty}$ on the complement of $\mbox{Supp}\,D \cup \mbox{Supp}\, N$,
\item $\sup |u_{\delta}| < + \infty$, 
\item $|\Delta_{\omega_{P,\delta}}u_{\delta}|$  
 is bounded, where $\Delta_{\omega_{P,\delta}}$ denotes the Laplacian with respect to the orbifold K\"{a}hler form 
$\omega_{P,\delta}$, 
\item The solution $u_{\delta}$ satisfying the above properties is unique. 
\end{enumerate} 
We set 
\begin{equation}\label{tilde}
\tilde{\omega}_{P,\delta}:= \omega_{P,\delta} +\sqrt{-1}\partial\bar{\partial}u_{\delta}. 
\end{equation}
The following lemma asserts that $\{\omega_{P,\delta}^{n}\}$ 
is (weakly) monotone decreasing with respect to $\delta$.  
\begin{lemma}\label{monotonicity}
For $0< \delta < \delta^{\prime} < 1$, we have that 
\[
\tilde{\omega}_{P,\delta}^{n} \geqq \tilde{\omega}_{P,\delta^{\prime}}^{n}
\]
holds on $U = X \backslash \mbox{\em Supp}\, E$. \fbox{}
\end{lemma}
{\em Proof.}  Let $\delta < \delta^{\prime}$ be positive numbers 
as above. By (\ref{MD}), 
\begin{equation}\label{diff}
\frac{(\omega_{P,\delta} + \sqrt{-1}\partial\bar{\partial}u_{\delta})^{n}}
{(\omega_{P,\delta^{\prime}} + \sqrt{-1}\partial\bar{\partial}u_{\delta^{\prime}})^{n}}
= \left((\prod_{i}(1 - \parallel\sigma_{i}\parallel^{\frac{2}{b}}))^{-\varepsilon}\parallel\sigma_{E}\parallel_{h_{E}}^{2}\right)^{\delta-\delta^{\prime}}\!\!\!\!\!\!\!\!\cdot e^{u_{\delta}-u_{\delta^{\prime}}}. 
\end{equation}
holds on $X$. 
We set 
\begin{equation}\label{diff2}
w_{\delta,\delta^{\prime}}:= u_{\delta} - u_{\delta^{\prime}}
+ (\delta - \delta^{\prime})\log \left((\prod_{i}(1 - \parallel\sigma_{i}\parallel^{\frac{2}{b}}))^{-\varepsilon}\parallel\sigma_{E}\parallel_{h_{E}}^{2}\right). 
\end{equation}
Then by (\ref{diff}) we see that
\begin{equation}\label{lap}
\int_{0}^{1}\tilde{\Delta}_{t}w_{\delta,\delta^{\prime}}\,dt 
= w_{\delta,\delta^{\prime}} 
\end{equation}
holds, where $\tilde{\Delta}_{t}$ denotes 
$\mbox{trace}_{(1-t)\tilde{\omega}_{\delta}+ t\omega_{\delta^{\prime}}}\sqrt{-1}\partial\bar{\partial}$. 
Since 
\[
(\delta - \delta^{\prime})\log \left((\prod_{i}(1 - \parallel\sigma_{i}\parallel^{\frac{2}{b}}))^{-\varepsilon}\parallel\sigma_{E}\parallel_{h_{E}}^{2}\right) \to  + \infty \,\,\,\,\,\,\,\,\mbox{as $x\to\mbox{Supp}\,E$},  
\]  
by the boundedness of $u_{\delta},u_{\delta^{\prime}}$ and the definition of 
$w_{\delta,\delta^{\prime}}$ (cf. (\ref{diff2})), there exists a point 
$p_{0}$ where $w_{\delta,\delta^{\prime}}$ takes its minimum.
Then by  (\ref{lap}), we see that $w_{\delta,\delta^{\prime}}(p_{0}) \geqq 0$ 
holds.   Hence $w_{\delta,\delta^{\prime}}\geqq 0$ holds on $U$. 
Hence by (\ref{diff}) and (\ref{diff2}), we see that 
\[
\tilde{\omega}_{\delta}^{n} \geqq \tilde{\omega}_{\delta^{\prime}}^{n}
\]
holds on $U$. This completes the proof of Lemma \ref{monotonicity}. 
q.e.d. \vspace{5mm} \\

\noindent To obtain a uniform estimate of $u_{\delta}$  with respect to $\delta$, we set 
\begin{equation}\label{vde}
v_{\delta} := u_{\delta} - (1-\delta)\cdot\log \left((\prod_{i}(1 - \parallel\sigma_{i}\parallel^{\frac{2}{b}}))^{-\varepsilon}\parallel\sigma_{E}\parallel_{h_{E}}^{2}\right) 
\end{equation}
and estimate $v_{\delta}$. 
By (\ref{MD}) $v_{\delta}$ satisfies the equation:  
\begin{equation}\label{veq}
(\omega_{P,1} + \sqrt{-1}\partial\bar{\partial}v_{\delta})^{n}= \frac{\parallel\sigma_{N}\parallel_{h_{N}}^{2}\left((\prod_{i}(1 - \parallel\sigma_{i}\parallel^{\frac{2}{b}}))^{-\varepsilon}\parallel\sigma_{E}\parallel_{h_{E}}^{2}\right)}{\parallel\sigma_{D}\parallel^{2}_{h_{D}}}\cdot \Omega\cdot e^{v_{\delta}}. 
\end{equation}
Since $u_{\delta}$ is bounded on $X$, we see that
\[
\log \frac{(\omega_{P,1} + \sqrt{-1}\partial\bar{\partial}v_{\delta})^{n}}
{\frac{\parallel\sigma_{N}\parallel_{h_{N}}^{2}\left((\prod_{i}(1 - \parallel\sigma_{i}\parallel^{\frac{2}{b}}))^{-\varepsilon}\parallel\sigma_{E}\parallel_{h_{E}}^{2}\right)}{\parallel\sigma_{D}\parallel^{2}_{h_{D}}}\cdot \Omega} = v_{\delta}
\]
is bounded from below and blows up along $E$. 
Hence we see that there exists a point $q_{0}$ on $X\backslash \mbox{Supp}\,E$
where $v_{\delta}$ takes minimum. 
Then by the maximum (minimum) principle, we see that 
\begin{equation}
\tilde{\omega}_{\delta}(q_{0})\geqq \omega_{P,1}(q_{0})
\end{equation}
holds.  Hence if we set 
\begin{equation}
C_{-} := \inf \frac{\omega_{P,1}^{n}}{\frac{\parallel\sigma_{N}\parallel_{h_{N}}^{2}\left((\prod_{i}(1 - \parallel\sigma_{i}\parallel^{\frac{2}{b}}))^{-\varepsilon}\parallel\sigma_{E}\parallel_{h_{E}}^{2}\right)}{\parallel\sigma_{D}\parallel^{2}_{h_{D}}}\cdot \Omega} > 0, 
\end{equation}
then 
\begin{equation}\label{lowerv}
v_{\delta} > C_{-}
\end{equation}
holds.  We note that $C_{-}$ is independent of $\delta \in (0,1]$.
Let us fix $0 < \delta_{0} < 1$. 
We set 
\begin{equation}\label{vd}
v_{\delta,\delta_{0}} := u_{\delta} - (\delta_{0}-\delta)\cdot\log \left((\prod_{i}(1 - \parallel\sigma_{i}\parallel^{\frac{2}{b}}))^{-\varepsilon}\parallel\sigma_{E}\parallel_{h_{E}}^{2}\right).
\end{equation}
Then by the same argument, we see that there exists a  constant 
$C_{-}(\delta_{0})$ such that for every $\delta\in (0,\delta_{0})$
\[
v_{\delta,\delta_{0}} \geqq C_{-}(\delta_{0})
\]
holds on $X$.  Hence by (\ref{vde}) for every $\delta\in (0,\delta_{0})$ 
\begin{equation}\label{delta0}
v_{\delta} \geqq C_{-}(\delta_{0}) - (1-\delta_{0})
\cdot\log \left((\prod_{i}(1 - \parallel\sigma_{i}\parallel^{\frac{2}{b}}))^{-\varepsilon}\parallel\sigma_{E}\parallel_{h_{E}}^{2}\right).
\end{equation}
holds on $X$. \vspace{3mm} \\  

\noindent Next we shall estimate $v_{\delta}$ from above.
We note that by (\ref{MD}) 
\begin{eqnarray*}
\int_{X}\frac{\parallel\sigma_{N}\parallel_{h_{N}}^{2}\left((\prod_{i}(1 - \parallel\sigma_{i}\parallel^{\frac{2}{b}}))^{-\varepsilon}\parallel\sigma_{E}\parallel_{h_{E}}^{2}\right)^{\delta}}{\parallel\sigma_{D}\parallel^{2}_{h_{D}}}\cdot \Omega\cdot e^{u_{\delta}}&  =  &\int_{X}(\omega_{P,\delta} + \sqrt{-1}\partial\bar{\partial}u_{\delta})^{n}\\ & = & (2\pi)^{n}(P -\delta E)^{n}
\end{eqnarray*}
hold. 
Then by the concavity of logarithm, 
\begin{equation}\label{concave}
\int_{X} \left(u_{\delta}+\log\frac{\parallel\sigma_{N}\parallel_{h_{N}}^{2}\left((\prod_{i}(1 - \parallel\sigma_{i}\parallel^{\frac{2}{b}}))^{-\varepsilon}\parallel\sigma_{E}\parallel_{h_{E}}^{2}\right)^{\delta}}{\parallel\sigma_{D}\parallel^{2}_{h_{D}}}\right)\Omega  
\leqq \left(\log (2\pi)^{n}(P -\delta E)^{n}\right)\int_{X}\Omega
\end{equation}
holds.  Hence there exists a positive constant $C$ independent of $\delta\in (0,1]$
such that 
\begin{equation}\label{intineq}
\int_{X}u_{\delta}\,\,\Omega < C
\end{equation}
holds for every $\delta \in (0,1]$. Then since $\omega_{P,\delta} +\sqrt{-1}\partial\bar{\partial}u_{\delta}$ is a closed positive current on $X$ and 
$\omega_{P,\delta}$ is an orbifold K\"{a}hler form on $X$, 
we see that $u_{\delta}$ is almost pluri-subharmonic function on $X$ in the sense of orbifold.  Hence by the submeanvalue inequality for subharmonic functions, by (\ref{intineq}) there exists a positive constant $C_{+}$ independent 
of $\delta\in (0,\delta]$ such that 
\begin{equation}\label{upper}
u_{\delta} \leqq C_{+} 
\end{equation}
holds on $X$.  
Hence by (\ref{vde}), for every $\delta\in (0,1]$  
\begin{equation}\label{vupper}
v_{\delta} \leqq C_{+} - (1-\delta)\cdot\log \left((\prod_{i}(1 - \parallel\sigma_{i}\parallel^{\frac{2}{b}}))^{-\varepsilon}\parallel\sigma_{E}\parallel_{h_{E}}^{2}\right). 
\end{equation}
holds on $X$. 

Now we shall estimate the $C^{2}$-norm of $v_{\delta}$.
First we note that $\omega_{P,1}$ is a $C^{\infty}$-orbifold K\"{a}hler 
form on $X$.  In particular, the bisectional curvature $R_{\alpha\bar{\alpha}\beta\bar{\beta}}$ of $\omega_{P,1}$ 
is bounded on $X$.  
To get the $C^{2}$-estimate, we shall estimate
$e^{-Cv_{\delta}}(n + \Delta_{P,1}v_{\delta})$, 
where $C$ is a positive constant which will be specified later and  
and $\Delta_{P,1}$ denotes the Laplacian with respect to the orbifold
K\"{a}hler form $\omega_{P,1}$. 

\begin{lemma}\label{c2}(\cite[p. 127, Lemma 2.2]{tu0}))
We set 
\begin{equation}
f := \log \frac{\omega_{P,1}^{n}}{\Omega}.
\end{equation}
Let $C$ be a positive number such that 
\begin{equation}
C + \inf_{\alpha\neq \beta}R_{\alpha\bar{\alpha}\beta\bar{\beta}} > 1
\end{equation}
holds on $Y$, where $R_{\alpha\bar{\alpha}\beta\bar{\beta}}$ denotes the 
bisectional curvature of $\omega_{P,1}$. 

Then 
\begin{equation}
e^{Cv_{\delta}}\tilde{\Delta}_{\delta}(e^{-Cv_{\delta}}(n + \Delta_{P,1}\,v_{\delta}))
\geqq (n + \Delta_{P,1}\,v_{\delta}) 
\end{equation}
\begin{equation*}
 +  \Delta_{P,1}\left(f+ \log \frac{\parallel\sigma_{N}\parallel_{h_{N}}^{2}\left((\prod_{i}(1 - \parallel\sigma_{i}\parallel^{\frac{2}{b}}))^{-\varepsilon}\parallel\sigma_{E}\parallel_{h_{E}}^{2}\right)}{\parallel\sigma_{D}\parallel^{2}_{h_{D}}}\right) 
\end{equation*}
\begin{equation*}
-(n+n^{2}\inf_{\alpha\neq \beta}R_{\alpha\bar{\alpha}\beta\bar{\beta}} ) - C\cdot n(n+\Delta_{P,1}\,v_{\delta}) + 
\end{equation*}
\vspace{-5mm}
\begin{equation*}
(n+\Delta_{P,1}\,v_{\delta})^{\frac{n}{n-1}}
\left(\frac{\parallel\sigma_{N}\parallel_{h_{N}}^{2}\left((\prod_{i}(1 - \parallel\sigma_{i}\parallel^{\frac{2}{b}}))^{-\varepsilon}\parallel\sigma_{E}\parallel_{h_{E}}^{2}\right)}{\parallel\sigma_{D}\parallel^{2}_{h_{D}}}\right)^{-\frac{1}{n-1}}
\!\!\!\!\!\!\exp\left(-\frac{1}{n-1}(v_{\delta} + f)\right)
\end{equation*}
holds, where $\Delta_{P,1}$ denotes the Laplacian with respect to $\omega_{P,1}$
\mbox{\em (}i.e., $\Delta_{P,1}= \mbox{\em trace}_{\omega_{P,1}}\sqrt{-1}\partial\bar{\partial}$\mbox{\em )} and 
$\tilde{\Delta}_{\delta}$ denotes the Laplacian with respect to $\tilde{\omega}_{\delta}$. 

\fbox{} \vspace{5mm}
\end{lemma}

\noindent Since $u_{\delta}$ is bounded, by the definition of $v_{\delta}$ (cf. (\ref{vde})) and the boundedness of $\Delta_{1,\delta}u_{\delta}$, we see 
that there exists a point $x_{0}\in X \backslash \mbox{Supp}\, E$ such that 
$e^{-Cv_{\delta}}(n +\Delta_{P,1}v_{\delta})$ takes its maximum 
at $x_{0}$. 
Then by Lemma \ref{c2} and the lower estimate (\ref{delta0})(taking $\delta_{0} = 2/3$ for example), 
if we take $C$ sufficiently large, we see that there exists a positive constant
 $C_{2}$ such that for every $\delta \in (0,1/2)$ 
\begin{equation}
e^{-Cv_{\delta}(x_{0})}(n + \Delta_{P,1}v_{\delta})(x_{0}) \leqq C_{2} 
\hspace{10mm}(\delta\in (0,1/2))
\end{equation}
holds. 
This implies that for every $\delta \in (0,1/2)$  
\begin{equation}\label{osc}
n + \Delta_{P,1}v_{\delta}
\leqq C_{2}e^{Cv_{\delta}} \leqq C_{2}\cdot e^{CC_{+}} \left((\prod_{i}(1 - \parallel\sigma_{i}\parallel^{\frac{2}{b}}))^{-\varepsilon}\parallel\sigma_{E}\parallel_{h_{E}}^{2}\right)^{-C(1-\delta)} 
\end{equation}
holds on $X$.   
Hence we see that $\{v_{\delta}|\,\,\delta\in (0,1/2)\}$ is uniformly $C^{2}$-bounded on every compact subset of $X \backslash \mbox{Supp}\,E$. 
Then by the general theory of fully nonlinear equations of 2nd order (\cite{tr}), we see that $\{v_{\delta}|\,\,\delta\in (0,1/2)\}$ is uniformly $C^{2,\alpha}$ bounded for some $\alpha \in (0,1)$ on every compact subset of $X \backslash \mbox{Supp}\,E$.   
Then by the standard theory for linear elliptic partial differential equations of  2nd order, we see that for every $k \geqq 0$, $\{v_{\delta}|\delta\in (0,1/2)\}$ is uniformly $C^{k}$-bounded on every compact subset of $X \backslash \mbox{Supp}\,E$. Hence there exists a sequence $\{\delta_{k}\},\delta_{k}\downarrow 0$ such that \begin{equation}
v: = \lim_{k\rightarrow\infty}v_{\delta_{k}}
\end{equation}
exists in $C^{\infty}$-topology on every compact subset of $X \backslash \mbox{Supp}\, E$.   
Then by (\ref{vde}) and (\ref{MD}), we see that 
\begin{equation}
u := v +  \log \left((\prod_{i}(1 - \parallel\sigma_{i}\parallel^{\frac{2}{b}}))^{-\varepsilon}\parallel\sigma_{E}\parallel_{h_{E}}^{2}\right)
\end{equation}
is $C^{\infty}$ on $X \backslash \mbox{Supp}\,E$ and satisfies the Monge-Amp\`{e}re
equation (\ref{MP}) and $\omega_{K}:= \omega_{P}+\sqrt{-1}\partial\bar{\partial}u$ is a well defined closed semipositive current by the estimates (\ref{upper}) and (\ref{lowerv}).   Moreover by Lemma \ref{monotonicity}, $u$ is independent of the choice of 
the subsequence $\{\delta_{k}\}$. 
We set 
\begin{equation}
\omega_{K} := \omega_{P} + \sqrt{-1}\partial\bar{\partial}u.
\end{equation}
To prove that $\omega_{K}$ is the canonical K\"{a}hler-Einstein current 
on $(X,D)$, we need to show that 
\[
h_{K} := \left(\frac{1}{n!}\,\omega_{K}^{n}\right)^{-1}\!\!\cdot h_{\sigma_{D}}
\]
is an AZD of $K_{X} + D$.   
By (\ref{vde}),(\ref{delta0}) and (\ref{upper}), we obtain the 
following {\bf almost boundedness} of $u_{\delta}$ and $u$. 
\begin{lemma}\label{almost}
For every $\delta_{0}\in (0,1]$, 
there exists a constant $C_{-}(\delta_{0})$ such that 
for every  $\delta \in (0,\delta_{0})$ 
\[
C_{-}(\delta_{0}) + (\delta_{0}-\delta)\log 
\left((\prod_{i}(1 - \parallel\sigma_{i}\parallel^{\frac{2}{b}}))^{-\varepsilon}\parallel\sigma_{E}\parallel_{h_{E}}^{2}\right)  
\leqq u_{\delta} \leqq C_{+}
\]
holds.  And in particular 
\[
C_{-}(\delta_{0}) + \delta_{0}\log 
\left((\prod_{i}(1 - \parallel\sigma_{i}\parallel^{\frac{2}{b}}))^{-\varepsilon}\parallel\sigma_{E}\parallel_{h_{E}}^{2}\right)  
\leqq u \leqq C_{+}
\] 
holds.  
\fbox{} 
\end{lemma}
Then by (\ref{MP}) and Lemma \ref{almost}, letting $\delta_{0}$ tend 
to $0$, we see that $h_{K}$ is an AZD of $K_{X} + D$. 

The uniqueness is the direct consequence of the dynamical construction 
(cf. Theorem \ref{unique} below).
This completes the proof of Theorem \ref{main}.  q.e.d.

\subsection{Monotonicity Lemma and K\"{a}hler-Einstein currents on LC pairs}
Using Theorem \ref{main}, we shall construct a canonical K\"{a}hler-Einstein 
current on a LC pair of log general type.

Let $X$ be a smooth projective variety and let $D$ be an effective 
$\mathbb{Q}$-divisor on $X$ such that $(X,D)$ is a LC pair. 
Then for every rational number $t \in [0,1)$, we see that 
$(X,tD)$ is a KLT pair.   
Suppose that $(X,D)$ is of log general type. 
Then for a sufficiently small rational number $0< \epsilon \ll  1$,
$(X,(1-\epsilon)D)$ is a KLT pair of log general type.
Let $\epsilon_{0}$ be a positive number such that 
for every $t\in (1-\epsilon_{0},1)$, $(X,tD)$ is of log general type.   
For $t\in (1-\epsilon_{0},1)$, let $\omega_{t}$ be the canonical K\"{a}hler-Einstein current on $(X,tD)$  as in Theorem \ref{main}. 
We set 
\begin{equation}
d\mu_{can,t}: = \frac{1}{n!}\,\omega_{t,abc}^{n}
\end{equation} 
and call it the canonical K\"{a}hler-Einstein volume form on 
$(X,tD)$. 
$d\mu_{can,t}$ is nothing but the canonical measure on $(X,tD)$.
The following monotonicity lemma is essential for our purpose.   
\begin{lemma}(Monotonicity Lemma)\label{monotonicity2}
$d\mu_{can,t}$ is (weakly) monotone increasing with respect to $t \in (1-\epsilon_{0},1)\cap \mathbb{Q}$. \fbox{}
\end{lemma}
{\em Proof of Lemma \ref{monotonicity2}}. 
Let $t < t^{\prime}$ be positive numbers in $(1-\epsilon_{0},1)$.
Taking a suitable modification of $X$,   
we may assume that $\mbox{Supp}\, D$ is a divisor with normal crossings and  there exist Zariski decompositions:
\[
K_{X} + tD = P_{t} + N_{t} 
\]
and 
\[
K_{X} + t^{\prime}D = P_{t^{\prime}} + N_{t^{\prime}}. 
\]
We note that adding exceptional divisor does not affect the log 
canonical ring. 
It is clear that
\begin{equation}\label{ineq}
N_{t^{\prime}} \preceq N_{t}
\end{equation}
holds. 
Let $h_{P_{t}}$ be a $C^{\infty}$-hermitian metric on $P_{t}$ induced by 
the Fubini-Study metric by the morphism $\Phi_{|\nu!P_{t}|}$ 
associated with the base point free linear system
$|\nu!P_{t}|$ from $X$ into a projective space for some sufficiently large
$\nu$.
We set 
\begin{equation}
\omega_{P_{t}} := \sqrt{-1}\,\Theta_{h_{P_{t}}}. 
\end{equation}
Let $h_{tD}$ be a $C^{\infty}$-hermitian metric on $tD$ and let 
$h_{N_{t}}$ be a $C^{\infty}$-hermitian metric on $N_{t}$.
Let $\sigma_{tD}$ be a multivalued global holomorphic section on 
$tD$ with divisor $tD$ and let $\sigma_{N_{t}}$ be a multivalued global holomorphic section of $N_{t}$ with divisor $N_{t}$. 
Let $\Omega$ be a $C^{\infty}$-volume form on $X$ such that 
\[
h_{P_{t}} = \Omega^{-1}\cdot h_{tD}\cdot h_{N_{t}}^{-1}
\]
holds. 
We consider 
\begin{equation}\label{MPt}
(\omega_{P_{t}} + \sqrt{-1}\partial\bar{\partial}u_{t})^{n}= \frac{\parallel\sigma_{N_{t}}\parallel_{h_{N_{t}}}^{2}}{\parallel\sigma_{tD}\parallel_{h_{tD}}^{2}}\cdot \Omega\cdot e^{u_{t}}
\end{equation}
on $X$ as (\ref{MP}) in Section 2.1 such that 
\begin{equation}\label{omegat}
\omega_{t} = \omega_{P_{t}} + \sqrt{-1}\partial\bar{\partial}u_{t} 
+ 2\pi t[D]
\end{equation}
is the canonical K\"{a}hler-Einstein current on $X$.  
As in Section 2.1, let $E$ be an effective $\mathbb{Q}$-divisor on $X$ such that  $P_{t} - E$ is ample and  $X \backslash \mbox{Supp}\, E$ 
is contained in 
\begin{equation}
W:= \{x\in X| d\mu_{can,t}(x), d\mu_{can,t^{\prime}}(x) > 0\} \backslash \mbox{Supp}\, D
\end{equation}
By the proof of Theorem \ref{main}, $W$ is a nonempty Zariski open subset of $X$ and 
$\omega_{t},\omega_{t^{\prime}}$ are $C^{\infty}$ on $W$. 
Let $D = \sum d_{i}D_{i}$ be the irreducible decomposition of $D$ and 
let $\sigma_{i}$ be a nontrivial global section of $\mathcal{O}_{X}(D_{i})$
with divisor $D_{i}$ and let $\parallel\sigma_{i}\parallel$ be the 
hermitian norm of $\sigma_{i}$ with respect to a $C^{\infty}$-hermitian metric 
on $\mathcal{O}_{X}(D_{i})$,  We may and do assume that 
$\parallel\sigma_{i}\parallel < 1$  holds for every $i$ on $X$. 
Let $b$ be a positive integer such that 
\[
d_{i} < \frac{b - 1}{b}
\]
holds for every $i$ as (\ref{choice}). 
Let $\sigma_{E}$ be a global multivalued holomorphic section of $E$ with divisor $E$ and 
let $h_{E}$ be a $C^{\infty}$-hermitian metric on $E$ such that 
$h_{P_{t}}\cdot h_{E}^{-1}$ has strictly positive curvature on $X$.
For every $0 < \delta \leqq 1$, we define the orbifold K\"{a}hler form $\omega_{P,\delta}$ on $X$ by 
\begin{equation}\label{Pdelta}
\omega_{P_{t},\delta}:= (1-\delta)\omega_{P_{t}} + \delta\left(\omega_{P_{t}}- \sqrt{-1}\,\Theta_{h_{E}}
 +\varepsilon \sum_{i}\sqrt{-1}\partial\bar{\partial}
\log (1 - \parallel\sigma_{i}\parallel^{\frac{2}{b}})\right), 
\end{equation}  
where $\varepsilon$ is a fixed sufficiently small positive number 
so that $\omega_{P,\delta}$ is an orbifold K\"{a}hler form on 
$X$ branching along $D$ with order $m$ for every $\delta\in (0,1]$.
 
For $\delta > 0$ we consider the perturbed equation: 
\begin{equation}\label{MDt}
(\omega_{P_{t},\delta} + \sqrt{-1}\partial\bar{\partial}u_{t,\delta})^{n}= \frac{\parallel\sigma_{N_{t}}\parallel_{h_{N_{t}}}^{2}\left((\prod_{i}(1 - \parallel\sigma_{i}\parallel^{\frac{2}{b}}))^{-\varepsilon}\parallel\sigma_{E}\parallel_{h_{E}}^{2}\right)^{\delta}}{\parallel\sigma_{tD}\parallel^{2}_{h_{tD}}}\cdot \Omega\cdot e^{u_{t,\delta}} 
\end{equation}
We see that (\ref{MDt}) has a unique bounded solution $u_{t,\delta}$ whose 
$C^{2}$-norm with respect to $\omega_{P_{t},\delta}$ 
is bounded by \cite[p.387,Theorem 6]{y}.  
By the uniqueness of the solution of (\ref{MPt}) (cf. Theorem \ref{unique}) 
and the uniform weighted $C^{2}$-estimate of $\{ u_{t,\delta}| t\in (0,1)\}$ 
parallel to the one of the solution of (\ref{MD}) in Section 2.1, 
we see that 
\begin{equation}\label{conv}
u_{t} := \lim_{\delta\downarrow 0} u_{t,\delta}
\end{equation}
exists in $C^{\infty}$-topology on $X \backslash \mbox{Supp}\, E$ and is the unique solution of (\ref{MPt}) under the condition that 
$u_{t}$ is almost bounded on $X$ (cf. Lemma \ref{almost}) or equivalently 
$\omega_{P_{t}}+\sqrt{-1}\partial\bar{\partial}u_{t}$ is the canonical K\"{a}hler-Einstein current on $(X,tD)$. 
We set 
\begin{equation}
\omega_{t,\delta} := \omega_{P_{t},\delta} + \sqrt{-1}\partial\bar{\partial}u_{t,\delta}. 
\end{equation}
Then by the equation (\ref{MDt}), we see that 
\begin{equation}
-\mbox{Ric}_{\omega_{t,\delta}} = \omega_{t,\delta}
\end{equation}
holds on $X \backslash \mbox{Supp}\, E$.   

Now we shall compare $\omega_{t,\delta}^{n}(t\in (0,1))$ and 
$\omega^{n}_{t^{\prime}} (t < t^{\prime})$.    We note that  
\[
F_{\delta}(x) := \frac{\omega^{n}_{t^{\prime}}}{\omega^{n}_{t,\delta}}(x)
\]
tends to $+\infty$ as $x$ tends to $\mbox{Supp}\,E$  by (\ref{ineq}),(\ref{MDt}) and the equation (\ref{MPt}) replacing $t$ by $t^{\prime}$. 
Hence there exists a point $x_{0}\in X \backslash \mbox{Supp}\,E$  where  
$F(x)$ takes its minimum.  
Then by the maximum (minimum) principle and the K\"{a}hler-Einstein condition, 
we see that 
\begin{equation}
\omega_{t^{\prime}}(x_{0})\geqq \omega_{t,\delta}(x_{0})
\end{equation}  
holds.  Hence we see that $F_{\delta}(x) \leqq 1$ holds on $X$ and 
\begin{equation}
\omega_{t,\delta}^{n} \leqq \omega_{t^{\prime}}^{n}
\end{equation}
holds on $X$.  
Hence letting $\delta$ tend to $0$, by the convergence (\ref{conv}), we see that 
\begin{equation}
\omega_{t}^{n} \leqq \omega_{t^{\prime}}^{n}
\end{equation}
holds on $X$.  This completes the proof of Lemma \ref{monotonicity2}. 
q.e.d. \vspace{3mm}\\
Now we shall construct an AZD of $K_{X} + D$ by using Lemma \ref{monotonicity2}. 
\begin{theorem}\label{limitvol}
In the above notations 
\begin{equation}
d\mu_{can} := \lim_{t\uparrow 1}d\mu_{can,t}
\end{equation}
exists.  And $d\mu_{can}^{-1}\cdot h_{\sigma_{D}}$ is an AZD of $K_{X} + D$. \fbox{}
\end{theorem}
{\em Proof}. First we shall prove the convergence. 
 By Lemma \ref{monotonicity2}, it is enough to prove that 
$\{d\mu_{can,t}| t\in (1-\epsilon_{0},1)\}$ is locally uniformly bounded
on some nonempty Zariski open subset of $X$. 
Taking modification of $X$, we may assume that $\mbox{Supp}\, E$ 
is a divisor with normal crossings.  
Let $H$ be an ample divisor such that 
$\mbox{Supp}\, E + H$ is  a divisor with normal crossings on $X$.
We set 
\begin{equation}
U_{0}:= X \backslash (\mbox{Supp}\, E\cup H)
\end{equation}
Then by \cite{ko}, there exists a complete K\"{a}hler-Einstein form 
$\omega_{H}$ on $U_{0}$ such that 
\[
-\mbox{Ric}_{\omega_{H}} = \omega_{H}
\]
and it  extends to a closed positive current on $X$. 
Then by Yau's Schwarz lemma (\cite{y2}), we have that 
\begin{equation}\label{dom}
d\mu_{can,t} \leqq \frac{1}{n!}\,\omega_{H}^{n}
\end{equation}
holds on $U_{0}$. 
Hence by the monotonicity of $\{ d\mu_{can,t}\}$ (Lemma \ref{monotonicity2}), this completes the proof of the convergence of $d\mu_{can,t}$ 
as $t\uparrow 1$.
Hereafter we shall identify $d\mu_{can,t}^{-1}$ with the 
singular hermitian metric $d\mu_{can,t}^{-1}\cdot h_{\sigma_{D}}^{t}$ 
on $K_{X}+tD$. 

Next we shall prove that $d\mu_{can}^{-1}$ is an AZD of $K_{X}+D$.
Let us fix an arbitrary $t_{0} \in (1-\epsilon_{0},1)\cap \mathbb{Q}$.   
Then $K_{X} + t_{0}D$ is big.     
Let $\ell_{0}$   be a sufficiently large positive integer such that 
 $B : = \ell_{0}(K_{X}+ t_{0}D)$ is  an integral divisor and  $|B| \neq \emptyset$. 
Let $m$ be an arbitrary positive integer such that $m(K_{X}+D)$ is Cartier.
We note that $d\mu_{can,t}^{-1}$ is an AZD of $K_{X} + tD$ for every $t\in 
(1-\varepsilon_{0},1)\cap \mathbb{Q}$. 
Then for any positive integer $\ell$, if we set 
\begin{equation}
s := \frac{\ell_{0}t_{0} + \ell m}{\ell_{0} + \ell m}.
\end{equation}
\begin{equation}
H^{0}(X,\mathcal{O}_{X}(m\ell(K_{X}+D)+B)) 
\subseteq H^{0}(X,\mathcal{O}_{X}(m\ell+\ell_{0})
(K_{X}+D))\otimes 
\mathcal{I}(d\mu_{can,s}^{-(\ell m + \ell_{0})}))
\end{equation}
holds by the definition of $s$. 
We note that  $d\mu_{can,s} \leqq d\mu_{can}$ holds on $X$ by Lemma \ref{monotonicity2}. Hence we see that for every nonzero global section $\tau \in H^{0}(X,\mathcal{O}_{X}(B))$ and  arbitrary section  $\sigma \in H^{0}(X,\mathcal{O}_{X}(m(K_{X}+D)))$,  
\begin{equation}
\sigma^{\ell}\otimes \tau \in H^{0}(X,\mathcal{O}_{X}(m\ell+\ell_{0})
(K_{X}+D))\otimes 
\mathcal{I}(d\mu_{can}^{-(\ell m + \ell_{0})}))
\end{equation}
holds. 
Let us fix a $C^{\infty}$-volume form $dV$ on $X$. 
Let us take a positive integer $\ell$ sufficiently large so that 
\begin{equation}
\int_{X} (|\tau|^{2}\cdot d\mu_{can}^{-\ell_{0}})^{-\frac{1}{\ell}}dV < \infty
\end{equation} 
holds. 
Then by H\"{o}lder's inequality, we see that 
\begin{equation}
\int_{X}|\sigma|^{2}\cdot (d\mu_{can})^{-m} dV
\leqq \left(\int_{X}|\sigma^{\ell}\otimes \tau|^{2}(d\mu_{can})^{-(\ell m+\ell_{0})}dV\right)^{\frac{1}{\ell}}\cdot 
\left(\int_{X} (|\tau|^{2}\cdot d\mu_{can}^{-\ell_{0}})^{-\frac{1}{\ell}}dV\right)^{\frac{\ell -1}{\ell}} < + \infty
\end{equation}
hold.   Since $\sigma$ is arbitrary, this means that $d\mu_{can}^{-1}\cdot h_{\sigma_{D}}$ is an AZD 
of $K_{X} + D$. 
This completes the proof of Lemma \ref{limitvol}. 
q.e.d. 
\begin{remark}
There are many other ways to approximate the LC pair $(X,D)$ 
by a sequence KLT pairs $\{ (X,D_{k})| D_{k}\preceq D\}$ 
such that $D_{k}\preceq D_{k+1}$ and \\ $D_{k}\uparrow D$ 
as $k$ tends to infinity.  We may easily generalize Lemma \ref{monotonicity2}
including all such approximations. Hence we may generalize Theorem \ref{limitvol} including such approximations. 
\fbox{} 
\end{remark}
\begin{definition}\label{LCcanmeasure}
Let $(X,D)$ be a LC pair of log general type. 
Then $d\mu_{can}$ constructed as above is said to be 
{\bf the canonical measure} on $(X,D)$ \fbox{}
\end{definition} 
I believe that $d\mu_{can}$ constructed is $C^{\infty}$ on a 
nonempty Zariski open subset of $X$.  But at present it is not clear. 
In the following examples, the canonical measures are generically $C^{\infty}$. 
\begin{example}
Let $X$ be a smooth projective variety and let $D$ be a divisor with normal 
crossings on $X$ such that $K_{X} + D$ is ample. 
Then $(X,D)$ is a LC pair. 
Then by \cite{ko}, there exists a complete K\"{a}hler-Einstein form 
$\omega_{E}$ on $X \backslash D$ such that 
$-\mbox{\em Ric}_{\omega_{E}} = \omega_{E}$ holds on $X \backslash D$.   
Then $\omega_{E}$ extends to be a closed positive current on $X$
cohomologous to $2\pi c_{1}(K_{X}+D)$ and is a canonical K\"{a}hler-Einstein current on $X$. \fbox{} 
\end{example}

\begin{example}\label{orblimit}
Let $X$ be a smooth projective $n$-fold and let $D$ be a divisor 
with normal crossings on $X$ such that for every sufficiently small positive 
rational number $\varepsilon$, $K_{X} + (1-\varepsilon)D$ is ample.
Then as in \cite{chern},  for every sufficiently large positive 
integer $m$, there exists an orbifold K\"{a}hler-Einstein form 
$\omega_{m}$ on the KLT pair $(X,\frac{m-1}{m}D)$. 
Then as in \cite{ball}, we see that $\omega_{\infty}:= \lim_{m\to\infty}
\omega_{m}$ exists on $X$ as a closed positive current and 
in $C^{\infty}$-topology on $X \backslash D$. 
Then  $(n!)^{-1}\omega_{\infty}^{n}$ is the canonical measure on 
the LC pair $(X,D)$. \fbox{}
\end{example}

\section{Ricci iterations}\label{DsKE}
In this section, we shall construct  solutions of (\ref{MA}) and  (\ref{MP}) as the limit 
of dynamical systems of K\"{a}hler forms or closed semipositive currents.  
The dynamical systems are defined by  Ricci iterations.
This construction will be used to define a family of  dynamical systems
of Bergman kernels which converge to the canonical K\"{a}hler-Einstein current on a KLT pair of log general type. 
More precisely we use the Ricci iteration to eliminate the effect of 
singularities of singular hermitian metrics on the $\mathbb{Q}$-line bundle. 
In this section, we shall use the notations in the previous section, if 
without fear of confusion. 
 
\subsection{The case of twisted K\"{a}hler-Einstein forms}
Let $X$ be a smooth projective $n$-fold  and let $(L,h_{L})$ be a $C^{\infty}$-hermitian $\mathbb{Q}$-line bundle on $X$ such that $\sqrt{-1}\Theta_{h_{L}}\geqq 0$ on $X$.  
We shall assume that $K_{X} + L$ is ample. 
Let $a$ be a positive integer such that $aL$ is a genuine line bundle 
on $X$.   
Let $\omega_{0}$ be a $C^{\infty}$-K\"{a}hler form representing 
$2\pi a c_{1}(K_{X}+L)$.
For $m\geqq 0$, we shall define a sequence of K\"{a}hler forms 
$\{\omega_{m}\}_{m= 0}^{\infty}$ inductively by :
\begin{equation}\label{ricciind}
-\mbox{Ric}_{\omega_{m}} + \frac{a - 1}{a}\,\omega_{m-1} + \sqrt{-1}\,\Theta_{h_{L}} = \omega_{m}  (m = 1,2,\cdots  ). 
\end{equation}
Let us reduce (\ref{ricciind}) to the sequence of Monge-Amp\`{e}re equations 
 as follows. 
Let $\Omega$ be a $C^{\infty}$-volume form on $X$ such that 
\begin{equation}\label{Omega}
\omega_{0} = a\left(-\mbox{Ric}\,\Omega + \sqrt{-1}\,\Theta_{h_{L}}\right)
\end{equation}
holds.  
We define the sequence of $C^{\infty}$-functions $\{ u_{m}\}^{\infty}_{m=0}$ by  $u_{0} = 0$ and for $m\geqq 1$, by the successive equations:  
\begin{equation}\label{MAind}
\log \frac{(\omega_{0} +\sqrt{-1}\partial\bar{\partial}u_{m})^{n}}
{\Omega}  = u_{m} - \frac{a - 1}{a}\,u_{m-1}. 
\end{equation}
The existence of $\{ u_{m}\}_{m=0}^{\infty}$ follows from the solution 
of Calabi's conjecture (\cite{a,y}).
Then we  see that the sequence of K\"{a}hler forms: 
\begin{equation}
\omega_{m}:= \omega_{0} + \sqrt{-1}\partial\bar{\partial}u_{m}\,\,\,\,\,\,\, (m\geqq 0) 
\end{equation}
satisfy the successive equations (\ref{ricciind}).
Now we shall state the result. 
\begin{theorem}\label{main1}
Let $\{ \omega_{m}\}_{m=0}^{\infty}$ be the sequence of K\"{a}hler forms 
defined inductively by (\ref{ricciind}) as above. 
Then 
\begin{equation}
\omega := \frac{1}{a}\lim_{m\rightarrow\infty}\omega_{m}
\end{equation} 
exists on $X$ in $C^{\infty}$-topology and $\omega$ is 
a $C^{\infty}$-K\"{a}hler form on $X$. 
And $\omega$ satisfies 
\begin{equation}
-\mbox{\em Ric}_{\omega} + \sqrt{-1}\Theta_{h_{L}} = \omega 
\end{equation} 
on $X$, i.e., $\omega$ is the twisted K\"{a}hler-Einstein form 
associated with  $(X,(L,h_{L}))$ (cf. Definition \ref{TKE}). 
\fbox{}
\end{theorem} 
{\em Proof of Theorem \ref{main1}}. 
By (\ref{MAind}), we have that 
\begin{equation}\label{succ}
\log \frac{(\omega_{0} +\sqrt{-1}\partial\bar{\partial}u_{m})^{n}}
{(\omega_{0} +\sqrt{-1}\partial\bar{\partial}u_{m-1})^{n}}  = u_{m} - u_{m-1}
-\frac{a - 1}{a}\,(u_{m-1}- u_{m-2})
\end{equation}
holds for every $m\geqq 2$. 
We note that 
\begin{equation}\label{intlap}
\log \frac{(\omega_{0} +\sqrt{-1}\partial\bar{\partial}u_{m})^{n}}
{(\omega_{0} +\sqrt{-1}\partial\bar{\partial}u_{m-1})^{n}}
= \int_{0}^{1}\Delta_{t}(u_{m} - u_{m-1})dt
\end{equation}
holds, where $\Delta_{t} (t\in [0,1])$ denotes the Laplacian 
with respect to the K\"{a}hler form $(1-t)\omega_{m-1} + t\omega_{m}$.
Then by the maximum principle and (\ref{succ}), we see that 
\begin{equation}\label{sup}
\sup_{X}(u_{m} -u_{m-1}) \leqq \frac{a - 1}{a}\sup_{X}(u_{m-1} - u_{m-2})
\end{equation}
and 
\begin{equation}\label{inf}
\inf_{X}(u_{m} -u_{m-1}) \geqq \frac{a - 1}{a}\inf_{X}(u_{m-1} - u_{m-2})
\end{equation}
hold for every $m \geqq 2$. 
Hence by (\ref{sup}) and (\ref{inf}), 
\begin{equation}\label{zero}
\lim_{m\to\infty} (u_{m}- u_{m-1}) = 0
\end{equation}
holds and  there exists a positive constant $C_{0}$ independent of $m$  such that 
\begin{equation}\label{c0}
|u_{m}| \leqq C_{0}
\end{equation}
holds on $X$. 

Now we shall consider the uniform $C^{2}$-estimate of $\{ u_{m}\}$.
We follow the argument in  \cite[Theorem 3.3]{ru} modeled after \cite{b-k}.  
Let $f : (M,g) \to (N,h)$ be a holomorphic map between K\"{a}hler manifolds.
The Chern-Lu inequality in the context of Yau's Schwarz lemma (\cite{lu,y2}) gives: 
\begin{equation}\label{schwarz}
\Delta_{g}\log |\partial f|^{2}\geqq \frac{\mbox{Ric}_{g}(\partial f,\bar{\partial} f)}{|\partial f|^{2}} - \frac{\mbox{Bisec}_{h}(\partial f,\bar{\partial}f,\partial f,\bar{\partial}f)}{|\partial f|^{2}},
\end{equation}
where $\Delta_{g}$ denotes the Laplacian with respect to $g$. 
Let $f$ be the identity map $1_{X}:(X,\omega_{m})\to (X,\omega_{0})$.
By rewriting the equation (\ref{ricciind}) as: 
\[
\mbox{Ric}_{\omega_{m}}= \frac{a - 1}{a}\omega_{m-1} + \sqrt{-1}\Theta_{h_{L}}
- \omega_{m}, 
\]
we see that Ricci curvature of $(X,\omega_{m})$ is uniformly bounded from 
below along the Ricci iteration by our assumption that $\sqrt{-1}\Theta_{h_{L}}$ is semipositive. 
Since 
\[
|\partial f|^{2} = \mbox{tr}_{\omega_{m}}\omega_{0} = n - \Delta_{\omega_{m}}u_{m},
\]
we see that 
\begin{equation}\label{important}
\Delta_{\omega_{m}}\log (n -\Delta_{\omega_{m}}u_{m})
\geqq - C(1 + n -\Delta_{\omega_{m}}u_{m})
\end{equation}
holds for some positive constant $C$ independent of $m$. 
Hence 
\begin{equation}
\Delta_{\omega_{m}}\left(\log (n-\Delta_{\omega_{m}}u_{m})
- (C+1)u_{m}\right) \geqq - n - C(1+n) + (n-\Delta_{\omega_{m}}u_{m})
\end{equation}
holds. 
Hence by the maximum principle and the uniform estimate (\ref{c0}) we see that there exists a positive constant $C_{2}$ such that 
\begin{equation}\label{c2m}
\parallel u_{m}\parallel_{C^{2}} \leqq C_{2}
\end{equation}
holds, where $\parallel\,\,\,\,\parallel_{C^{2}}$ denotes 
the $C^{2}$-norm with respect to $\omega_{0}$. 
Then by the general theory of nonlinear elliptic equations of 2nd order (\cite{tr}), 
we see that there exist a positive number $\alpha < 1$ and a positive constant 
$C_{2,\alpha}$ such that 
\begin{equation}
\parallel u_{m}\parallel_{C^{2,\alpha}} \leqq C_{2,\alpha}
\end{equation} 
holds,  where $\parallel\,\,\,\,\parallel_{C^{2,\alpha}}$ denotes 
the $C^{2,\alpha}$-norm with respect to $\omega_{0}$.
Hence there exists a subsequence $\{m_{k}\}$ of $\mathbb{N}$ such that 
\begin{equation}
u := \lim_{k\to\infty}u_{m_{k}}
\end{equation}
exists in $C^{2,\alpha}$-topology.  
Then by (\ref{MAind}) and (\ref{zero}),
we see that  $u$ satisfies the equation:
\begin{equation}\label{eq2}
\log \frac{(\omega_{0} +\sqrt{-1}\partial\bar{\partial}u)^{n}}
{\Omega}  = \frac{1}{a}\,u
\end{equation}
on $X$.  We note that (\ref{eq2}) has a unique solution as in 
\cite{y}.  
Hence by (\ref{sup}) and (\ref{inf}) the limit $u$ exists in $C^{2,\alpha}$-topology without taking the subsequence $\{ m_{k}\}$.   
By taking $\sqrt{-1}\partial\bar{\partial}$ of the both sides of 
(\ref{eq2}), by  the definition of $\Omega$ (cf. (\ref{Omega})), 
we see that if we define $\omega$ by 
\begin{equation}
\omega := \frac{1}{a}(\omega_{0} + \sqrt{-1}\partial\bar{\partial}u),
\end{equation}
then 
\begin{equation}
-\mbox{Ric}_{\omega} + \sqrt{-1}\Theta_{h_{L}} = \omega
\end{equation}   
holds.   This completes the proof of Theorem \ref{main1}. q.e.d.  

\subsection{The case of canonical K\"{a}hler-Einstein currents}\label{CKE}
Next we shall construct the  canonical K\"{a}hler-Einstein current 
in Theorem \ref{main} in terms of Ricci iterations.   
The major difference is that the $C^{\infty}$-hermitian $\mathbb{Q}$-line bundle $(L,h_{L})$  is replaced by a singular hermitian $\mathbb{Q}$-line bundle $(D,h_{\sigma_{D}})$ (cf. (\ref{sigmad}) for the definition of $h_{\sigma_{D}}$). 

Let $X$ be a smooth projective  $n$-fold  and let $D$ be an effective $\mathbb{Q}$-divisor on $X$ such that $(X,D)$ is a KLT pair. 
We assume that $\mbox{Supp}\, D$ is a divisor with normal crossings on $X$ and 
$K_{X}+D$ is big. 
We assume that there exists a Zariski decomposition :
\begin{equation}\label{ZD2}
K_{X} + D = P + N  \,\,\,\,\,(P,N \in \mbox{Div}(X)\otimes \mathbb{Q})
\end{equation}
i.e.,  $P$ is semiample, $N$ is effective  and 
\begin{equation}
H^{0}(X,\mathcal{O}_{X}(maP))
\simeq H^{0}(X,\mathcal{O}(ma(K_{X} + D)))
\end{equation}
holds for every $m\geqq 0$, where $a$ is the minimal positive integer 
such that $aD, aP, aN \in \mbox{Div}(X)$. 
We may assume that $(X,D)$ satisfies the above conditions without loss 
of generality as was observed in Section 2.1.  

Moreover since adding exceptional effective $\mathbb{Q}$-divisors 
does not change the log canonical ring and the canonical K\"{a}hler-Einstein 
current depends only on the log canonical ring essentially as is seen 
in Theorem \ref{unique} below.  
\vspace{2mm} \\

\noindent Let $h_{P}$ be a $C^{\infty}$-hermitian metric on $P$ with semipositive 
curvature.  
We set 
\begin{equation}\label{u}
U : = \{ x\in X|\,\,\,\mbox{$|\nu!P|$ is  very ample on a neighbourhood of 
$x$ for every $\nu \gg 1$}\}.
\end{equation} 
Let $a$ be a positive integer such that $aD\in \mbox{Div}(X)$. 
For $m\geqq 0$, we shall define inductively a sequence of closed positive current 
$\{\omega_{m}\}_{m=0}^{\infty}$ satisfying the following conditions : 
\begin{itemize}
\item[(P1)] $\omega_{0} = a\left(\sqrt{-1}\Theta_{h_{P}} + 2\pi[N]\right)$, 
where $[N]$ denotes the current of integration over $N$.  
\item[(P2)] The cohomology class of $[\omega_{m}]$ of $\omega_{m}$ 
is equal to $2\pi a\cdot c_{1}(K_{X}+D)$.
\item[(P3)] $\omega_{m}$ is $C^{\infty}$ on $U$.
\item[(P4)] $\{\omega_{m}\}_{m=0}^{\infty}$ satisfies the equation :
\begin{equation}\label{ricciind2}
-\mbox{Ric}_{\omega_{m}} + \frac{a - 1}{a}\,\omega_{m-1} + 2\pi[D] = \omega_{m}
\end{equation}
on $U$ for every $m\geqq 1$.
\item[(P5)] We define the singular hermitian metric $h_{m}$ by 
\begin{equation}
h_{m}:= \left(\left(\frac{1}{n!}\omega_{m,abc}^{n}\right)^{-1}\cdot h_{m-1}^{a-1}\cdot h_{\sigma_{D}}\right)^{\frac{1}{a}}
\end{equation}
on $K_{X}+D$. 
Then $h_{m}$ is an AZD of $K_{X}+D$ for every $m\geqq 1$. 
\end{itemize}
Now we shall state the main result in this subsection.  
\begin{theorem}\label{iteration2}
The dynamical system $\{ \omega_{m}\}_{m=0}^{\infty}$ 
satisfying (P1)-(P5) above exists 
and is unique and the limit 
\begin{equation}
\omega_{K} := \frac{1}{a}\lim_{m\rightarrow\infty}\omega_{m}
\end{equation}
exists in $C^{\infty}$-topology on every compact subset of $U$ (cf.{\rm (\ref{u})}) and also in the sense of current on $X$. 
The closed positive current $\omega_{K}$ satisfies the equation :
\[
-\mbox{\em Ric}_{\omega_{K}} + 2\pi[D] = \omega_{K}
\]
on $X$ and  $\omega_{K}$ is the canonical K\"{a}hler-Einstein current on 
$(X,D)$ (cf. Definition \ref{GKE}).  \fbox{}
\end{theorem}
Let us reduce (\ref{ricciind2}) to  successive Monge-Amp\`{e}re
equations. 
Let $\sigma_{N}$ be a multivalued holomorphic section of $N$ with divisor 
$N$. 
Let $h_{D},h_{N}$ are $C^{\infty}$-hermitian metric on $D$ and $N$ respectively.
Let $\Omega$ be a $C^{\infty}$-volume form on $X$ such that 
\begin{equation}
h_{P} := \Omega^{-1}\cdot h_{D}\cdot h_{N}^{-1}. 
\end{equation}
holds. Then 
\begin{equation}\label{Omega}
\omega_{0} = a\left(-\mbox{Ric}\,\Omega + 2\pi[N]\right)
\end{equation}
holds.  
We set 
\begin{equation}
\omega_{P} := \sqrt{-1}\Theta_{h_{P}}. 
\end{equation}
As in Section 2.1 let $E$ be an effective $\mathbb{Q}$-divisor such that 
$P - \delta E$ is ample for every  
$\delta \in (0,1]$. As before we may and do assume that $X \backslash \mbox{Supp}\, E = U$ 
holds (cf. (\ref{ed}) in Section 2.1). 
Then there exists a $C^{\infty}$-hermitian metric $h_{E}$ on $E$ 
such that $h_{P}\cdot h_{E}^{-1}$ is a metric with strictly positive 
curvature on $X$.
We define the sequence of functions $\{ u_{m}\}^{\infty}_{m=0}$ by 
$u_{0} = 0$ and for $m\geqq 1$ as follows.  
\begin{itemize}
\item[(Q1)] $\{u_{m}\}$ are almost bounded in the sense of 
Lemma \ref{almost} in Section 2.1, i.e.,
for every $\delta_{0}\in (0,1]$ 
there exists a constant $C_{-}(\delta_{0})$ depending only on $\delta_{0}\in (0,1]$ and a constant $C_{+}$ independent of $m$  such that 
\[
C_{-}(\delta_{0}) + \delta_{0}\log 
\left((\prod_{i}(1 - \parallel\sigma_{i}\parallel^{\frac{2}{b}}))^{-\varepsilon}\parallel\sigma_{E}\parallel_{h_{E}}^{2}\right)  
\leqq u_{m} \leqq C_{+}
\] 
holds, where we have used the notations in Section \ref{ex1} (cf. (\ref{choice}),(\ref{Pdelta}) etc.). 
 
\item[(Q2)] $u_{m} \in C^{\infty}(X\backslash \mbox{Supp} E)$.
\item[(Q3)] $\{ u_{m}\}$ satisfy the successive equations:  
\begin{equation}\label{MAind2}
\log \frac{(a\omega_{P} +\sqrt{-1}\partial\bar{\partial}u_{m})^{n}}
{\Omega}  = \log \frac{\parallel\sigma_{N}\parallel^{2}_{h_{N}}}{\parallel\sigma_{D}\parallel^{2}_{h_{D}}}+ \left(u_{m} - \frac{a - 1}{a}\,u_{m-1}\right)  
\end{equation}
\end{itemize}
for $m \geqq 1$. 
Then as in the Section 2, (\ref{MP}), (\ref{MAind2}) is 
equivalent to (\ref{ricciind2}), i.e.,  the sequence of closed positive 
currents $\{\omega_{m}\}$ defined by 
\begin{equation}\label{um}
\omega_{m} := a(\omega_{P} + 2\pi[N]) + \sqrt{-1}\partial\bar{\partial}u_{m}
\end{equation}
satisfies the conditions (P1) to (P5) above. 
For example the almost boundedness (Q1) implies the condition (P5) and 
(Q2) implies (P3) etc. 

The construction of  $\{ u_{m}\}_{m=0}^{\infty}$ is 
very much similar to the one of Theorem \ref{main} in 
Section 2.  
The only difference is that we need to estimate inductively.  
To construct $\{ u_{m}\}_{m=0}^{\infty}$ we consider the perturbed equation as follows.
For a fixed $\delta\in (0,1]$, we set 
\begin{equation}
\omega_{m,\delta,0}:= \omega_{P} + \delta\left(\frac{a -1}{a}\right)^{m}
\left(\sqrt{-1}\,\Theta_{h_{E}} + \varepsilon \sqrt{-1}\partial\bar{\partial}\sum_{i} \log (1 -\parallel\sigma_{i}\parallel^{\frac{2}{b}})\right),
\end{equation} 
where $\varepsilon$ is a sufficiently small positive number 
such that $\omega_{m,\delta,0}$ is an orbifold K\"{a}hler form on 
$X$ for every $\delta\in (0,1]$ as in Secion 2.1. 
We consider the successive equations:
\begin{equation}\label{MAind3}
\log \frac{(a\omega_{m,\delta,0} +\sqrt{-1}\partial\bar{\partial}u_{m,\delta})^{n}}
{\Omega}   = \log \frac{\parallel\sigma_{N}\parallel^{2}_{h_{N}}}{\parallel\sigma_{D}\parallel^{2}_{h_{D}}}+ \left(u_{m,\delta} - \frac{a - 1}{a}\,u_{m-1,\delta}\right) 
\end{equation}
for $m \geqq 1$ and we set 
\begin{equation}
\omega_{m,\delta}:= a\omega_{m,\delta,0} + \sqrt{-1}\partial\bar{\partial}u_{m,\delta}.
\end{equation} 
We set $u_{0,\delta} = 0$.  Then inductively by using 
the same strategy as in the proof of Theorem \ref{main} in Section 2, 
by \cite[p.387,Theorem 6]{y}, we see that (\ref{MAind3}) has a sequence of bounded solutions 
$\{ u_{m,\delta}\}_{m=0}^{\infty}$ such that for every $m$,  $u_{m,\delta}$ has
 a bounded  $C^{2}$-norm on $X$ with respect to the orbifold K\"{a}hler form 
 $\omega_{1,1,0}$.  
\begin{lemma}\label{u1}
For every $\epsilon > 0$ there exists a positive constant
$C_{0,\epsilon}$ such that  
\[
u_{1,\delta} \geqq -C_{0,\epsilon} + \epsilon\log \parallel\sigma_{E}\parallel_{h_{E}}^{2}
\]
holds for every $\delta \in (0,1]$. 
\fbox{}
\end{lemma} 
{\em Proof}.
We set 
\begin{equation}\label{epsdef}
u_{1,\delta,\epsilon}:= u_{1,\delta} -   
\epsilon \left(\log \parallel\sigma_{E}\parallel_{h_{E}} -\varepsilon\sum_{i} \log (1 -\parallel\sigma_{i}\parallel^{\frac{2}{b}})\right)
\end{equation}
and 
\[
\omega_{1,\delta,0,\epsilon}:= \omega_{1,\delta,0}
+ \epsilon \sqrt{-1}\left(-\Theta_{h_{E}}
- \varepsilon\partial\bar{\partial}\sum_{i}\log (1 -\parallel\sigma_{i}\parallel^{\frac{2}{b}})\right). 
\]
Then $u_{1,\delta,\epsilon}$ satisfies the equation:
\begin{equation}\label{e1}
\log \frac{(a\omega_{1,\delta,0,\epsilon} +\sqrt{-1}\partial\bar{\partial}u_{1,\delta,\epsilon})^{n}}
{\omega_{1,\delta,0,\epsilon}^{n}}   = \log \frac{\parallel\sigma_{N}\parallel^{2}_{h_{N}}\cdot\parallel\sigma_{E}\parallel^{2\epsilon}_{h_{E}}\cdot\Omega}{\parallel\sigma_{D}\parallel^{2}_{h_{D}}\prod_{i}(1 -\parallel\sigma_{i}\parallel^{\frac{2}{b}})^{\epsilon\varepsilon}\cdot\omega_{1,\delta,0,\epsilon}^{n}}+ u_{1,\delta,\epsilon}. 
\end{equation}
We note that since $\omega_{1,\delta,0,\epsilon}$ 
is an orbifold K\"{a}hler form, by the choice of $b$ (cf. (\ref{choice})), there exists a positive constant  $C_{0,\epsilon}$ 
depending only on $\epsilon$ such that 
\[
\log \frac{\parallel\sigma_{N}\parallel^{2}_{h_{N}}\cdot\parallel\sigma_{E}\parallel^{2\epsilon}_{h_{E}}\cdot\Omega}{\parallel\sigma_{D}\parallel^{2}_{h_{D}}\cdot\prod_{i}(1 -\parallel\sigma_{i}\parallel^{\frac{2}{b}})^{\epsilon\varepsilon}\cdot\omega_{1,\delta,0,\epsilon}^{n}} \leqq C_{0,\epsilon}
\]
holds on $X$.  We note that since $u_{1,\delta}$ is bounded on $X$, 
$u_{1,\delta,\epsilon}(x)$ tends to $+\infty$ as $x$ tends to $\mbox{Supp}\, E$.  
Hence by applying the maximum priciple to (\ref{e1}), we obtain that 
\[
u_{1,\delta,\epsilon} \geqq - C_{0,\epsilon}.
\]
holds on $X$. 
By (\ref{epsdef}), we obtain the lemma. 
q.e.d. \vspace{3mm}\\
As (\ref{succ}) we have that 
\begin{equation}\label{succ2}
\log \frac{(a\omega_{m,\delta,0} +\sqrt{-1}\partial\bar{\partial}u_{m,\delta})^{n}}
{(a\omega_{m-1,\delta,0} +\sqrt{-1}\partial\bar{\partial}u_{m-1,\delta})^{n}}  = \ (u_{m,\delta} - u_{m-1,\delta}) - \frac{a - 1}{a}(u_{m-1,\delta}-u_{m-2,\delta}) 
\end{equation}
for $m \geqq 2$. 
We note that on $X \backslash \mbox{Supp}\,E$, 
\begin{equation}
\log \frac{(a\omega_{m,\delta,0} +\sqrt{-1}\partial\bar{\partial}u_{m,\delta})^{n}}
{(a\omega_{m-1,\delta,0} +\sqrt{-1}\partial\bar{\partial}u_{m-1,\delta})^{n}}
=
\log \frac{(a\omega_{m-1,\delta,0} +\sqrt{-1}\partial\bar{\partial}\left(u_{m,\delta}-\alpha_{m}\delta\log \parallel\sigma_{E}\parallel^{2}_{h_{E}}\right))^{n}}
{(a\omega_{m-1,\delta,0} +\sqrt{-1}\partial\bar{\partial}u_{m-1,\delta})^{n}}
\end{equation}
holds, where 
\begin{equation}\label{alpham}
\alpha_{m} := \frac{1}{a}\left(\frac{a-1}{a}\right)^{m-1}. 
\end{equation}
The definition of $\{\alpha_{m}\}$ follows from the successive equations 
in cohomology classes:   
\begin{equation}
2\pi[P] +\frac{a -1}{a}[\omega_{m-1,\delta}] = [\omega_{m,\delta}] \,\,\,\,\,\,\,\, (m\geqq 1)
\end{equation} 
arising from (\ref{ricciind2}). \\
We note that the equality:
\begin{equation}\label{trick}
\!\!\!\!\!\!\!\!(u_{m,\delta} - u_{m-1,\delta}) - \frac{a - 1}{a}(u_{m-1,\delta}-u_{m-2,\delta})   =  
\end{equation}
\[
  \left(u_{m,\delta} - u_{m-1,\delta}- \alpha_{m}\delta\log \parallel\sigma_{E}\parallel^{2}_{h_{E}}\right)
  - \frac{a - 1}{a}\left(u_{m-1,\delta}-u_{m-2,\delta}- \alpha_{m-1}\delta\log \parallel\sigma_{E}\parallel^{2}_{h_{E}}\right) 
\]
holds.  We also note that  by the boundedness of $|u_{m,\delta}|$ and $|u_{m-1,\delta}|$ 
on $X$,   
\begin{equation}
\left(u_{m,\delta} -u_{m-1,\delta}
- \alpha_{m}\delta\log \parallel\sigma_{E}\parallel^{2}_{h_{E}}\right)(x)
\to +\infty \,\,\,\,\,\,\,\,\,\,\mbox{as $x\to \mbox{Supp}\,E$}
\end{equation}
holds. 
Hence applying  the minimum principle to (\ref{succ2}), by using the similar formula as (\ref{intlap}), 
\begin{equation}\label{inf2}
\inf_{X}(u_{m,\delta} -u_{m-1,\delta}
- \alpha_{m}\delta\log \parallel\sigma_{E}\parallel^{2}_{h_{E}}) \geqq \frac{a - 1}{a}\inf_{X}(u_{m-1,\delta} - u_{m-2,\delta} -\alpha_{m-1}\delta\log \parallel\sigma_{E}\parallel^{2}_{h_{E}})
\end{equation}
holds for every $m \geqq 2$.
Hence  $\{u_{m,\delta}\}$ is an almost monotone increasing sequence in this sense  and by Lemma \ref{u1} (taking $\epsilon > 0$ sufficiently small), there exists a positive constant $C_{0}(\delta)$ depending on $\delta$ such that 
\begin{equation}\label{inf3}
u_{m,\delta} - u_{m-1,\delta} - \alpha_{m}\delta\log\parallel\sigma_{E}\parallel^{2}_{h_{E}}
\geqq  -C_{0}(\delta)\left(\frac{a-1}{a}\right)^{m-1}
\end{equation}
holds for every $m\geqq 1$. 
By (\ref{alpham}) and (\ref{inf3}), we see that
\begin{equation}\label{lowerum}
u_{m,\delta} \geqq - C_{0}(\delta)\cdot a + \delta \log \parallel\sigma_{E}\parallel^{2}_{h_{E}}
\end{equation}
holds.  But this estimate depends on $\delta$.

To get the uniform lower estimate with respect to $\delta$, we set 
\begin{equation}\label{umdeltaepsilon}
u_{m,\delta}(\epsilon):= u_{m,\delta} - \epsilon\left(\frac{a-1}{a}\right)^{m}
\log\parallel\sigma_{E}\parallel^{2}_{h_{E}}
\end{equation}
for  $m \geqq 0$. 
Then by Lemma \ref{u1}, the similar argument as above replacing 
$u_{m,\delta}$ by $u_{m,\delta}(\epsilon)$ and 
$\omega_{m,\delta,0}$ by 
\begin{equation}
\omega_{m,\delta,0}(\epsilon) : = \omega_{m,\delta,0} - \epsilon\sqrt{-1}\left(\frac{a -1}{a}\right)^{m}\Theta_{h_{E}}
\end{equation}
we see that 
for every $\epsilon > 0$, 
there exists a positive constant $C_{0,\epsilon}$ (same as in Lemma \ref{u1}) independent of $\delta$ 
such that 
\begin{equation}\label{indlower}
u_{m,\delta}(\epsilon) - u_{m-1,\delta}(\epsilon)- \alpha_{m}(\delta+\epsilon)\log\parallel\sigma_{E}\parallel^{2}_{h_{E}} \geqq -C_{0,\epsilon}\left(\frac{a-1}{a}\right)^{m-1}
\end{equation}
holds.
Hence we have the uniform lower estimate (with respect to $\delta$): 
\begin{equation}\label{uniflower}
u_{m,\delta}\geqq - C_{0,\epsilon}\cdot a + \left\{\left(\frac{a-1}{a}\right)^{m}\epsilon + (\delta +\epsilon)\right\}\log\parallel\sigma_{E}\parallel_{h_{E}}^{2}. 
\end{equation} 

On the other hand, since 
\[
\int_{X}\exp\left(u_{m.\delta}-\frac{a-1}{a}u_{m-1,\delta}\right)\Omega
 = (2\pi)^{n}\left(P - \left(\frac{a-1}{a}\right)^{m}\delta E\right)^{n}
\]
holds, by the concavity of logarithm as (\ref{concave}) and (\ref{intineq}), we see that there exists a 
positive constant $C_{0}$ independent of $m$ and $\delta$ such that 
\begin{equation}
\int_{X}\left(u_{m.\delta}-\frac{a-1}{a}u_{m-1,\delta}\right)\Omega
\leqq C_{0}
\end{equation}
holds.   Hence by the almost monotonicity and the almost pluri-subharmonicity:
$\omega_{m,\delta} = \omega_{m,\delta,0} + \sqrt{-1}\partial\bar{\partial}u_{m,\delta} \geqq 0$,   
we see that there exists a positive constant $C_{+}$ independent of 
$m$ and $\delta$ such that 
\begin{equation}\label{upperbdd}
u_{m,\delta} \leqq C_{+}
\end{equation} 
holds.   Hence by the almost monotonicity of $\{ u_{m,\delta}\}_{m=0}^{\infty}$, we see that 
\begin{equation}
u_{\delta} := \lim_{m\to\infty}u_{m,\delta}
\end{equation}
and 
\begin{equation}
u := \lim_{m\to\infty}u_{m}
\end{equation}
exist.  
Let $\omega_{P,\delta} (\delta\in (0,1]$) be the orbifold K\"{a}hler form 
defined as (\ref{Pdelta})  in Section 2, i.e.,
\begin{equation}
\omega_{P,\delta}:= \omega_{P} + \delta \left(-\sqrt{-1}\Theta_{h_{E}}
+ \varepsilon\sum_{i}\sqrt{-1}\partial\bar{\partial}\log (1 -\parallel\sigma_{i}\parallel^{\frac{2}{b}})\right).  
\end{equation}  
As in Section 2, we consider the perturbed equation: 
\begin{equation}\label{MD4}
(a\omega_{P,1} + \sqrt{-1}\partial\bar{\partial}v_{m,\delta})^{n}= \frac{\parallel\sigma_{N}\parallel_{h_{N}}^{2}\left((\prod_{i}(1 - \parallel\sigma_{i}\parallel^{\frac{2}{b}}))^{-\varepsilon}\parallel\sigma_{E}\parallel_{h_{E}}^{2}\right)^{\frac{1}{a}}}{\parallel\sigma_{D}\parallel^{2}_{h_{D}}}\cdot 
\exp\left(v_{m,\delta}-\frac{a-1}{a}v_{m-1,\delta}\right)\cdot \Omega, 
\end{equation}
where 
\begin{equation}\label{vmdelta}
v_{m,\delta}:= u_{m,\delta} - (1-\alpha_{m}\delta)\cdot\log \left((\prod_{i}(1 - \parallel\sigma_{i}\parallel^{\frac{2}{b}}))^{-\varepsilon}\parallel\sigma_{E}\parallel_{h_{E}}^{2}\right).  
\end{equation}
This equation is equivalent to 
\begin{equation}\label{ricci33}
-\mbox{Ric}_{\omega_{m,\delta}} + \frac{a - 1}{a}\omega_{m-1,\delta} +2\pi [D]
= \omega_{m,\delta}
\end{equation}
as before, where 
\begin{equation}
\omega_{m,\delta}: = a(\omega_{P,1} + 2\pi [N]) + \sqrt{-1}\partial\bar{\partial}v_{m,\delta}.
\end{equation}

\noindent Then  by the weighted inductive $C^{2}$-estimate as in the proof of Theorem \ref{main} (or the  uniform weighted $C^{2}$-estimate below as 
the proof of Theorem \ref{main1} using Schwarz type lemma) 
letting $\delta\downarrow 0$, we see that $\{ u_{m,\delta}\}$ exists 
for every $m\geqq 0$ and $\delta \in (0,1]$. 
Here we note that $\Delta_{P,1}v_{m,\delta}$ is bounded on $X$ for every 
$m$ and $\delta > 0$, hence, we may take the constant $C$ in Lemma \ref{c2}
independent of $m$ and $\delta$, since $e^{-Cv_{m,\delta}}(n +\Delta_{P,1}v_{m,\delta})$ is bounded on $X$. 
Then as in the proof of Theorem \ref{main}, we have the following lemma.  

\begin{lemma}\label{exseq}
\begin{equation}
u_{m}: = \lim_{\delta\downarrow 0} u_{m,\delta}
\end{equation}
exists on $X \backslash \mbox{\em Supp}\,E$ compact uniformly in $C^{\infty}$-topology for every $m$  and $\{u_{m}\}$ satisfies the family of  equations 
(\ref{ricciind2}).   And $\{ u_{m}\}$ satisfies the properties (Q1),(Q2) and (Q3) above. \fbox{}
\end{lemma}
Then by Theorem \ref{unique} below, we see that $\{ u_{m}\}_{m=1}^{\infty}$  is unique and gives a sequence of closed positive currents 
$\{\omega_{m}\}$ satisfying (P1)-(P5) above. 
But we should note that the above inductive estimate may not be uniform 
with respect to $m$.  Apriori the estimate may get worse as 
$m$ tends to infinity, because to estimate the $C^{2}$-norm of $u_{m,\delta}$
by using Lemma \ref{c2}, the (weighted) $C^{2}$-norm of $u_{m-1,\delta}$ comes into the estimate. The following estimate gives a better $C^{2}$-estimate
of $\{ u_{m}\}$ (or $\{ u_{m,\delta}\}$).    

To prove the convergence, we proceed just as the proof of Theorem \ref{main1}.
But here we need to handle the incompleteness of $X \backslash \mbox{Supp}\, E$
with respect to the orbifold K\"{a}hler form $\omega_{P,1}$ (cf. (\ref{Pdelta})). 
First we note that by (\ref{ricci33}), $\mbox{Ric}_{\omega_{m,\delta}}$ is 
uniformly bounded from below along the Ricci iteration just as before. 
We note that the bisectional curvature of the K\"{a}hler orbifold
 $(X,a\omega_{P,1})$  is uniformly bounded.   
Then applying the Schwarz type lemma: (\ref{schwarz}) to the identity map 
$1_{X}:(X,\omega_{m,\delta})\to (X,a\omega_{P,1})$, as (\ref{important}), we see that
there exists a positive constant $C_{1}$ independent of $m$ and $\delta$ such that  
\begin{equation}\label{keyineq}
\Delta_{\omega_{m,\delta}}\left(\log (n-\Delta_{\omega_{m,\delta}}v_{m,\delta})
- (C_{1}+1)v_{m,\delta}\right) \geqq  n - C_{1}(1+n) + (n-\Delta_{\omega_{m,\delta}}v_{m,\delta}),
\end{equation}
holds on $U$.
We note that $v_{m,\delta}$ blows up (to $+\infty$) toward $\mbox{Supp}\, E$ by the definition (\ref{vmdelta}) and the $C^{0}$-estimate (\ref{uniflower})
(taking $\epsilon$ sufficiently small).
Then by the inductive $C^{2}$-estimate of $v_{m,\delta}$ using Lemma \ref{c2} as in 
Section 2 (cf. (\ref{osc})), if we take $C_{1}$ (independent of $m$ and $\delta$) sufficiently large,  then 
\[
\log (n-\Delta_{\omega_{m,\delta}}v_{m,\delta})
- (C_{1}+1)v_{m,\delta} \to -\infty \hspace{10mm}\mbox{as $x\to \mbox{Supp}\, E$}
\]
holds and 
there exists a point $x_{m,\delta}$ on $X \backslash\mbox{Supp}\, E$, where 
$\log (n-\Delta_{\omega_{m,\delta}}v_{m,\delta})
- (C_{1}+1)v_{m,\delta}$ takes it maximum.  
Now we are ready to apply the maximum principle to (\ref{keyineq}) and we see  
 that 
\begin{equation}
\log (n- \Delta_{\omega_{m,\delta}}v_{m,\delta})(x_{m,\delta}) \leqq  C_{1}(1+n) -n   
\end{equation}
holds.  Since $v_{m,\delta}$ is uniformly bounded from below by 
the definition of $v_{m,\delta}$ (cf. (\ref{vmdelta})) and 
the uniform lower estimate (\ref{uniflower}), there exists a positive constant $C_{2}$ independent of $m$ and $\delta$ such that 
\begin{equation}
\log (n-\Delta_{\omega_{m,\delta}}v_{m,\delta}) - (C_{1}+1)v_{m,\delta} \leqq  C_{2}
\end{equation}
holds on $U$.
Letting $\delta$ tend to $0$, this gives a weighted uniform $C^{2}$-estimate of $\{ u_{m}\}$ on every compact subset of $X \backslash \mbox{Supp}\, E$. 
Hence using \cite{tr} and the almost monotonicity of $\{ u_{m}\}$, we conclude that 
\begin{equation}
u = \lim_{m\rightarrow\infty}u_{m}
\end{equation} 
exists on $X \backslash \mbox{Supp}\,E$ compact uniformly in $C^{\infty}$-topology without taking a subsequence and $u$ is almost bounded by (\ref{uniflower})
and (\ref{upperbdd}).
   
Hence by (\ref{MAind2}), $u$  satisfies the Monge-Amp\`{e}re equation: 
\begin{equation}\label{P-eq}
\log \frac{(a\omega_{P} +\sqrt{-1}\partial\bar{\partial}u)^{n}}
{\Omega}  = \log \frac{\parallel\sigma_{N}\parallel_{h_{N}}}{\parallel\sigma_{D}\parallel_{h_{D}}}+ \frac{1}{a}\,u. 
\end{equation}
Then it is clear that 
\[
\omega_{K}:=\frac{1}{a}\left(a\omega_{P}  + \sqrt{-1}\partial\bar{\partial}u\right) + 2\pi [N] 
\]
satisfies the equation:
\[
-\mbox{Ric}_{\omega_{K}} + 2\pi[D] = \omega_{K}.
\] 
Now we define the singular hermitian metric: 
\[
h_{K} : = \left(\frac{1}{n!}\omega_{K,abc}^{n}\right)^{-1}\!\!\cdot h_{\sigma_{D}}
\]
on $K_{X} + D$. 
Then by the equation (\ref{P-eq}) and the uniform $C^{0}$-estimates
 (\ref{uniflower}) from below, we see that $h_{K}$ is an AZD of $K_{X} + D$. 
Hence $\omega_{K}$ is nothing but the canonical K\"{a}hler-Einstein
current on $(X,D)$. 
This completes the proof of Theorem \ref{iteration2}. q.e.d.

\section{Dynamical systems of Bergman kernels}\label{Bergman}

In this section, we shall construct twisted K\"{a}hler-Einstein forms 
and canonical K\"{a}hler-Einstein currents (on KLT pairs of log general type)  
in terms of dynamical systems of Bergman kernels. 
The proof is essentially the same as in \cite{KE,canonical}.
But here the dynamical construction cannot be applied 
directly to the twisted K\"{a}hler-Einstein forms or 
canonical K\"{a}hler-Einstein currents. 
In fact we apply the dynamical construction as in \cite{KE,canonical}
to the Ricci iterations: (\ref{ricciind}) or (\ref{ricciind2}). 
Hence the dynamical construction here has two indices.

\subsection{Dynamical construction of twisted K\"{a}hler-Einstein forms}
Let $X$ be a smooth projective  $n$-fold and let $(L,h_{L})$ is a $C^{\infty}$-hermitian $\mathbb{Q}$-line bundle on $X$ such that 
\begin{enumerate}
\item[(1)] $\sqrt{-1}\,\Theta_{h_{L}}$ is semipositive, 
\item[(2)] $K_{X} + L$ is ample. 
\end{enumerate}  
Let $a$ be a positive number such that $aL$ is a genuine line bundle. 
Let $h_{0}$ be a $C^{\infty}$-hermitian metric on 
$a(K_{X} + L)$ such that 
\begin{equation}
\omega_{0}:= \sqrt{-1}\,\Theta_{h_{0}}
\end{equation}
is a $C^{\infty}$-K\"{a}hler form on $X$. 
Let us consider the dynamical system of Bergman kernel $\{K_{\ell,1}\}_{\ell=1}^{\infty}$ as follows. 
Let $A$ be a sufficiently ample line bundle on $X$ such that 
for every pseudo-effective singular hermitian line bundle $(F,h_{F})$ on 
$X$,  
${\cal O}_{X}(jK_{X}+ A + F)\otimes {\cal I}(h_{F})$ is globally generated 
for every $0 \leqq j\leqq a-1$.
Such an ample line bundle $A$ exists by Nadel's vanishing theorem (cf. \cite[p.561]{n}). 
Let $h_{A}$ be a $C^{\infty}$-hermitian metric on $A$ with strictly positive 
curvature.  
We set  
\begin{equation}\label{k11}
K_{1,1}:= K(X,A+a(K_{X}+L),h_{A}\cdot h_{0}^{\frac{a-1}{a}}\cdot h_{L}), 
\end{equation}
where $K(X,A+a(K_{X}+L),h_{A}\cdot h_{0}^{\frac{a-1}{a}}\cdot h_{L})$
denotes the Bergman kernel defined as (\ref{BergmanK}). 
By the choice of $A$, $\mathcal{O}_{X}(A + a(K_{X}+L))$ is globally generated.
 Hence we may  define the $C^{\infty}$-hermitian metric $h_{1,1}$ on $A+ a(K_{X}+L)$ by 
\begin{equation}
h_{1,1} := \frac{1}{K_{1,1}}. 
\end{equation} 
where $K(X,A+a(K_{X}+L),h_{A}\cdot h_{0}^{\frac{a-1}{a}}\cdot h_{L})$. 
We define the sequence of Bergman kerenels $\{ K_{\ell,1}\}_{\ell\geqq 1}$ 
and the sequence of hermitian metrics $\{ h_{\ell,1}\}_{\ell\geqq 1}$ 
inductively as follows.   
Suppose that we have already defined a $C^{\infty}$-hermitian metric 
on $A + (\ell-1)a(K_{X}+L)$. 
We define $K_{\ell,1}$ by 
\begin{equation}
K_{\ell,1} := K(X,A +\ell a(K_{X}+L),h_{\ell-1,1}\cdot h_{0}^{\frac{a-1}{a}}\cdot h_{L}).
\end{equation} 
By the choice of $A$, $K_{\ell,1}$ does not vanish everywhere on $X$.  
Then we define a $C^{\infty}$-hermitian metric $h_{\ell,1}$ on $A + \ell a(K_{X}+L)$ by 
\begin{equation}
h_{\ell,1}:= \frac{1}{K_{\ell,1}}. 
\end{equation}
In this way we define the dynamical system of Bergman kernels 
$\{ K_{\ell,1}\}_{\ell\geqq 1}$ and the dynamical system of $C^{\infty}$-hermitian metrics $\{ h_{\ell,1}\}_{\ell\geqq 1}$. 
Then as in \cite{canonical}, we have the following lemma.
\begin{lemma}
\begin{equation}
K_{1}:= \limsup_{\ell\rightarrow\infty}\sqrt[\ell ]{(\ell !)^{-n}K_{\ell,1}}
\end{equation}
is bounded and  
\begin{equation}
h_{1}:= \frac{n!}{(2\pi)^{n}}\,K_{1}^{-1}
\end{equation}
is a $C^{\infty}$-hermitian metric on 
$a(K_{X}+L)$ (Here the constant $(2\pi)^{-n}n!$ is not so important
anyway.  
The reason why we put such a constant will become clear in the next subsection.).
If we set 
\begin{equation}
\omega_{1}:= \sqrt{-1}\Theta_{h_{1}},
\end{equation}
then $\omega_{1}$ is a K\"{a}hler form on $X$ and satisfies the equation:
\begin{equation}
-\mbox{\em Ric}_{\omega_{1}} + \frac{a-1}{a}\omega_{0} + \sqrt{-1}\Theta_{h_{L}}
= \omega_{1}. 
\end{equation}
\fbox{}
\end{lemma}
Next replacing $h_{0},\omega_{0}$ by $h_{1},\omega_{1}$ respectively. 
We obtain the dynamical system of Bergman kernels 
$\{ K_{\ell,2}\}_{\ell=1}^{\infty}$.  
Again  
\begin{equation}
K_{2}:= \limsup_{\ell\rightarrow\infty}\sqrt[\ell ]{(\ell !)^{-n}K_{\ell,2}}
\end{equation}
exists and  $h_{2}:= \frac{n!}{(2\pi)^{n}}K_{2}^{-1}$ is a $C^{\infty}$-hermitian metric on 
$a(K_{X}+L)$.  If we set 
\begin{equation}
\omega_{2}:= \sqrt{-1}\,\Theta_{h_{2}},
\end{equation}
then $\omega_{2}$ is a K\"{a}hler form on $X$ and satisfies the equation:
\begin{equation}
-\mbox{Ric}_{\omega_{2}} + \frac{a-1}{a}\omega_{1} + \sqrt{-1}\Theta_{h_{L}}
= \omega_{2}. 
\end{equation}
Inductively, for every positive integers  $m$, we define the dynamical system of Bergman kernels
$\{ K_{\ell,m}\}_{\ell=1}^{\infty}$ and a $C^{\infty}$-hermitian metric $h_{m}$ on $a(K_{X}+L)$.   Then by Theorem \ref{main1} we have the following theorem. 

\begin{theorem}\label{doubleinduction}
Let $\{ K_{\ell,m}\}_{\ell=1}^{\infty}$ be the system of Bergman kernels 
defined as above.  If we define 
\begin{equation}\label{Km}
K_{m}:= \mbox{the upper-semi-continuous envelope of}\,\,\,\, 
\limsup_{\ell\rightarrow\infty}\sqrt[\ell]{(\ell !)^{-n}K_{\ell,m}},
\end{equation}
then 
\begin{equation}
h_{m} := \frac{n!}{(2\pi)^{n}}\,K_{m}^{-1}
\end{equation}
is a $C^{\infty}$-hermitian metric on $a(K_{X}+L)$ with strictly positive curvature (The upper-semi-continuous envelope of the supremum 
of a family of  pluri-subharmonic functions is again pluri-subharmonic by 
\cite[p.26, Theorem 5]{l}, if the family is locally uniformly bounded from above.  This operation is often used in this article.
But anyway the adjustment occurs only on the set of measure $0$.).
And if we define a K\"{a}hler form $\omega_{m}$ by 
\begin{equation}
\omega_{m}:= \sqrt{-1}\,\Theta_{h_{m}},
\end{equation}
then $\{\omega_{m}\}_{m=1}^{\infty}$ satisfy  the equations: 
\begin{equation}
-\mbox{\em Ric}_{\omega_{m}} + \frac{a-1}{a}\,\omega_{m-1} + \sqrt{-1}\Theta_{h_{L}}= \omega_{m} 
\end{equation}
on $X$ for every $m \geqq 1$ as {\em (\ref{ricciind})}. 
And 
\[
\omega:=\frac{1}{a}\lim_{m\rightarrow\infty}\omega_{m}
\]
exists on $X$ in $C^{\infty}$-topology and $\omega$ is the 
unique $C^{\infty}$-solution of the equation: 
\[
-\mbox{\em Ric}_{\omega} + \sqrt{-1}\Theta_{h_{L}} = \omega.
\]
on $X$, i.e., $\omega$ is the twisted K\"{a}hler-Einstein form 
associated with $(X,(L,h_{L}))$. 
\fbox{}
\end{theorem} 
\begin{remark}
By the proof below the convergence (\ref{Km}) exists  in $L^{1}$-topology. \fbox{}
\end{remark}

\subsection{Dynamical construction of canonical K\"{a}hler-Einstein currents}
In this subsection, we shall construct canonical K\"{a}hler-Einstein currents 
in terms of dynamical systems of Bergman kernels. 

Let $X$ be a smooth projective variety of dimension $n$ and let $D$ be an effective $\mathbb{Q}$-divisor on $X$ such that $(X,D)$ is a KLT pair. 
We assume that $\mbox{Supp}\, D$ is a divisor with normal crossings and 
$K_{X}+D$ is big. 
We also assume that there exists a Zariski decomposition :
\begin{equation}\label{ZD3}
K_{X} + D = P + N \,\,\,\,\,\,\,\,(P,N \in \mbox{Div}(X)\otimes \mathbb{Q}),  
\end{equation} 
i.e., $P$ is semiample, $N$ is effective and  
\begin{equation}
H^{0}(X,\mathcal{O}_{X}(maP))
\simeq H^{0}(X,\mathcal{O}(ma(K_{X} + D)))
\end{equation} 
holds for every $m\geqq 0$, where $a$ is the minimal positive integer 
such that $aD, aP, aN \in \mbox{Div}(X)$ hold. 
These assumptions are not actual restrictions as is explained in the beginning of Section \ref{ex1}.
 
We set 
\begin{equation}\label{u2}
U : = \{ x\in X|\,\,\mbox{$|\nu !P|$ is  very ample on a neighbourhood of 
$x$ for every $\nu \gg 1$}\}. 
\end{equation}
Let $\sigma_{D}$ be a multivalued holomorphic section of $D$ with divisor 
$D$.  And we define the singular hermitian metric $h_{\sigma_{D}}$ by 
\begin{equation}
h_{\sigma_{D}} := \frac{1}{|\sigma_{D}|^{2}}. 
\end{equation}  
Let $h_{P}$ be a $C^{\infty}$-hermitian metric on 
$P$ such that 
\begin{equation}
\omega_{P}:= a\sqrt{-1}\,\Theta_{h_{P}}
\end{equation}
is a $C^{\infty}$-K\"{a}hler form on $U$.  
Since $P$ is semiample, this is possible.  
Let $\sigma_{N}$ be a global multivalued holomorphic section of $N$ with divisor $N$. Let $h_{\sigma_{N}}$ be the singular hermitian metric on $N$ 
defined by 
\begin{equation}
h_{\sigma_{N}}:= \frac{1}{|\sigma_{N}|^{2}}. 
\end{equation}
Let us consider the dynamical system of Bergman kernel $\{K_{\ell,1}\}_{\ell=1}^{\infty}$ as follows. 
As before, let  $A$ be a sufficiently ample line bundle on $X$ such that 
for every pseudo-effective singular hermitian line bundle $(F,h_{F})$ on 
$X$, $\mathcal{O}_{X}(A + F)\otimes \mathcal{I}(h_{F})$ is  globally generated. Such a line bundle  $A$ exists by Nadel's vanishing theorem (\cite[p.561]{n}).

The following lemma is trivial.  But it explains why we can handle KLT pairs 
in the dynamical system of Bergman kernels. 
\begin{lemma}\label{AAA}
\begin{equation}\label{inje}
\otimes \sigma_{N}^{\ell a}: \mathcal{O}_{X}(\ell aP) \hookrightarrow 
\mathcal{O}_{X}(\ell a(K_{X}+D))\otimes \mathcal{I}((h_{P}\cdot h_{\sigma_{N}})^{\ell a-1}\cdot h_{\sigma_{D}}))
\end{equation}
is a well defined injection for every $\ell \geqq 0$. \fbox{}
\end{lemma} 
\noindent{\em Proof of Lemma \ref{AAA}}. 
Let $\sigma$ be an arbitrary holomorphic section of $\ell aP$ on an open subset  $V$ in $X$.   
Then $\sigma^{\frac{1}{\ell a}}\cdot \sigma_{N}$ is a multivalued holomorphic section
 of $K_{X} + D$ on $V$.  Hence we see that 
\[
\frac{\sigma^{\frac{1}{\ell a}}\cdot \sigma_{N}}{\sigma_{D}}
\] 
is locally $L^{2}$-integrable multivalued meromorphic form on $V$, because $(X,D)$ is KLT.  We note that the equality: 
\[
((h_{P}\cdot h_{\sigma_{N}})^{\ell a-1}\cdot h_{\sigma_{D}})|\sigma\cdot \sigma_{N}^{\ell a}|^{2} = \left|\frac{\sigma^{\frac{1}{\ell a}}\cdot\sigma_{N}}{\sigma_{D}}\right|^{2}\cdot \left(
(h_{P}\cdot h_{N})^{\ell a-1}\cdot\left|\sigma^{1 -\frac{1}{\ell a}}\cdot \sigma_{N}^{\ell a-1}\right|^{2}\right)
\] 
holds. 
Since the second factor of the right hand side: $(h_{P}\cdot h_{N})^{\ell a-1}\left|\sigma^{1 -\frac{1}{\ell a}}\cdot \sigma_{N}^{\ell a-1}\right|^{2}
= h_{P}^{\ell a- 1}|\sigma|^{2(1-1/\ell a)}$ is locally bounded, we see that (\ref{inje}) is a well defined injection.  q.e.d. \vspace{2mm}\\

\noindent Let $h_{A}$ be a $C^{\infty}$-hermitian metric on $A$ with strictly positive 
curvature.  
We set  
\begin{equation}
K_{1,1}:= K(X,A+a(K_{X}+D),h_{A}\cdot (h_{P}\cdot h_{\sigma_{N}})^{a-1}\cdot h_{\sigma_{D}})
\end{equation}
and define a singular hermitian metric $h_{1,1}$ on $A+ a(K_{X}+D)$ by 
\begin{equation}
h_{1,1} := \frac{1}{K_{1,1}}.
\end{equation}
By Lemma \ref{AAA} we have the natural injection: 
\begin{equation}
H^{0}(X,\mathcal{O}(A+ aP))\hookrightarrow H^{0}(X,\mathcal{O}_{X}(A+a(K_{X}+D))\otimes\mathcal{I}(h_{A}\cdot (h_{P}\cdot h_{\sigma_{N}})^{a-1}\cdot h_{\sigma_{D}})). 
\end{equation}
We note that by the choice of $A$, $\mathcal{O}(A+ aP)$ is globally generated.
 Hence we have that 
\begin{equation}
h_{1,1} = O (h_{A}\cdot h_{P}^{a}\cdot h_{\sigma_{N}}^{a})
\end{equation}
holds. 
Suppose that we have already defined $K_{1,\ell-1}$ and the 
singular hermitian metric $h_{1,\ell-1}$ on $A + (\ell-1)a(K_{X}+D)$ such that\begin{equation}\label{OM-1}
h_{\ell-1,1} = O(h_{A}\cdot h_{P}^{(\ell-1)a}\cdot h_{\sigma_{N}}^{(\ell-1) a}). 
\end{equation} 
Then  we define $K_{\ell,1}$ by 
\begin{equation}
K_{\ell,1} := K(X,A +\ell a(K_{X}+D),h_{\ell-1,1}\cdot (h_{P}\cdot h_{\sigma_{N}})^{a-1}\cdot h_{\sigma_{D}})
\end{equation}  
and define the singular hermitian metric $h_{\ell,1}$ on $A + \ell a(K_{X}+D)$ by 
\begin{equation}
h_{\ell,1}:= \frac{1}{K_{\ell,1}}.
\end{equation}
By the choice of $A$, we see that
$\mathcal{O}_{X}(A + \ell aP)$ is globally generated and by (\ref{OM-1}) and Lemma \ref{AAA}, there is a natural injection:
\[
H^{0}(X,\mathcal{O}_{X}(A + \ell aP))\hookrightarrow 
H^{0}(X,\mathcal{O}_{X}(A +\ell a(K_{X}+D))\otimes\mathcal{I}(h_{\ell-1,1}\cdot (h_{P}\cdot h_{\sigma_{N}})^{a-1}\cdot h_{\sigma_{D}})). 
\]
Hence 
$h_{\ell,1} = O(h_{A}\cdot h_{P}^{\ell a}\cdot h_{\sigma_{N}}^{\ell a})$  
holds. 
In this way by the induction on $\ell$,  we define 
the dynamical systems $\{ K_{\ell,1}\}_{\ell \geqq 1}$ and  $\{ h_{\ell,1}\}_{\ell\geqq 1}$.  By the above inductive construction, 
\begin{equation}\label{OM}
h_{\ell,1} = O(h_{A}\cdot h_{P}^{\ell a}\cdot h_{\sigma_{N}}^{\ell a})
\end{equation} 
holds for every $\ell \geqq 1$. \vspace{2mm} \\
Let us make (\ref{OM}) quantitative.
Let $\{ \omega_{m}\}_{m=1}^{\infty}$ be the Ricci iteration constructed 
as in Theorem \ref{iteration2} in  Section \ref{DsKE}.  
We set 
\begin{equation}\label{key}
h_{1}:= \left(\frac{1}{n!}\,\omega_{1,abc}^{n}\right)^{-1}\cdot (h_{P}\cdot h_{\sigma_{N}})^{a-1}\cdot h_{\sigma_{D}}.  
\end{equation}
Then by the equation (\ref{ricciind2}), we see that 
\begin{equation}\label{key2}
\omega_{1} = \sqrt{-1}\,\Theta_{h_{1}} = -\mbox{Ric}_{\omega_{1}} + \frac{a-1}{a}\omega_{0}
\end{equation}
hold on $X$, where $\omega_{0} = a(\sqrt{-1}\,\Theta_{h_{P}} + 2\pi[N])$ 
as in Theorem \ref{iteration2}.   

Let $E$ be an effective $\mathbb{Q}$-divisor on $X$ such that 
$M := P - E$ is ample.  Let $h_{M}$ be a $C^{\infty}$-hermitian metric 
on $M$ with strictly positive curvature. Then 
\[
\hat{h}_{M} := h_{M}\cdot h_{\sigma_{E}}\cdot h_{\sigma_{N}}
\]
is an singular hermitian metric with strictly positive curvature on $K_{X}+D$.  
For every positive number $\varepsilon < 1$, we set 
\[
h_{0,\varepsilon} := h_{0}^{1-\varepsilon}\cdot \hat{h}_{M}^{a\varepsilon}.
\]
and 
\[
\omega_{0,\varepsilon} = \sqrt{-1}\Theta_{h_{0,\varepsilon}}. 
\]
Let us fix a positive rational number $\varepsilon < 1$.  And we consider 
the auxiliary Ricci iteration: 
\begin{equation}
\omega_{m+1,\varepsilon}:= -\mbox{Ric}_{\omega_{m,\varepsilon}}+\frac{a-1}{a}\omega_{m,\varepsilon} + 2\pi[D]
\end{equation} 
for every $m$ and we require $[\omega_{m.abc}] = 2\pi [a(P -\varepsilon E)]$. 
We set 
\begin{equation}
h_{m,\varepsilon} := \left(\frac{1}{n!}\omega_{1,\varepsilon}^{n}\right)^{-1}\cdot h_{m-1,\varepsilon}^{\frac{a-1}{a}}. 
\end{equation}
Let $\ell \geqq 2$ and suppose that there exists a positive constant $C_{\ell-1}$ such that 
\begin{equation}\label{assumption}
K_{\ell-1,1} \geqq C_{\ell-1}\cdot h_{1,\varepsilon}^{-(\ell-1)}\cdot h_{A}^{-1}
\end{equation}
holds.  We note that by the extremal property of Bergman kernels,  
\begin{equation}\label{extremal}
K_{\ell ,1}(x) := \sup \left\{|\sigma|^{2}(x); 
\int_{X}|\sigma|^{2}\cdot \left(h_{\ell-1,1}\cdot\left(h_{P}\cdot h_{\sigma_{N}}\right)^{a-1}\cdot h_{\sigma_{D}}\right) = 1\right\}
\end{equation}
holds for every $x\in X$, where $\sigma$ runs all the elements in 
\[
H^{0}\left(X,\mathcal{O}_{X}(A+\ell a(K_{X}+D))\otimes \mathcal{I}(h_{\ell-1,1}\cdot\left(h_{P}\cdot h_{\sigma_{N}}\right)^{a-1}\cdot h_{\sigma_{D}})\right).
\]
And this vector space contains $H^{0}(X,\mathcal{O}_{X}(A+\ell aP))$ by Lemma \ref{AAA} and (\ref{OM}). 

Let us estimate $K_{\ell,1}$ from below by using the $L^{2}$-estimate for 
$\bar{\partial}$-operators and (\ref{extremal}) above. 
Let $p \in U (= X \backslash \mbox{Supp}\, E)$ be 
an arbitrary point.  
Then since $\omega_{1}$ is strictly positive on a neighbourhood of $p$, 
by (\ref{key}), there exists a holomorphic normal coordinate
$(V,(z_{1},\cdots ,z_{n}))$ of $(X,\omega_{1})$ around $p$  and a 
multivalued holomorphic functions $f_{D},f_{N}$ defining 
effective $\mathbb{Q}$-divisors $D,N$ respectively on $V$ such that 
\begin{equation}\label{expansion}
h_{1,\varepsilon} =  \{ \prod_{i=1}^{n}(1 -\mid z_{i}\mid^{2}) + O(\parallel z\parallel^{3})\}\cdot \frac{|f_{D}|^{2a}}{2^{-2an}|dz_{1}\wedge\cdots\wedge dz_{n}|^{2a}\cdot |f_{N}|^{2a}}
\end{equation} 
holds on $V$.
We may assume that $V$ is biholomorphic to the polydisk $\Delta^{n}(r)$ of radius $r$ with center $O$ in $\mathbb{C}^{n}$ for some $0 < r < 1$ via $(z_{1},\cdots ,z_{n})$ and $\bar{V} \subseteq U$. 
Taking $r < 1$ sufficiently small, we may assume that 
there exists a  $C^{\infty}$-function  $\rho$ on $X$ such that 
\begin{enumerate}
\item[(1)] $\rho$ is identically $1$ on $\Delta^{n}(r/3)$,
\item[(2)] $0\leqq \rho \leqq 1$, 
\item[(3)] $\mbox{Supp}\,\rho \subset\subset V$,
\item[(4)] $\mid d\rho\mid < 3/r$, where $\mid\,\,\,\,\,\mid$ denotes 
the pointwise norm with respect to $\omega_{1}$, 
\item[(5)] There exists a local holomorphic frame $\mbox{\bf e}_{A}$ of $A$ on $V$ such that\\ $h_{A}(\mbox{\bf e}_{A},\mbox{\bf e}_{A})(p) = 1$.
\end{enumerate} 
We set 
\begin{equation}
\tau_{\ell}:= \mbox{\bf e}_{A}\otimes\frac{(dz_{1}\wedge\cdots \wedge dz_{n})^{\otimes \ell a}\cdot f_{N}^{\ell a}}
{f_{D}^{\ell a}}. 
\end{equation}
Then the $L^{2}$-norm $\parallel\rho\cdot\tau_{\ell}\parallel$ of the cut off $\rho\cdot\tau_{\ell}$ with respect to 
$h_{1}^{\ell}\cdot h_{A}$ and $\omega_{1}$ concentrates around $p$ 
as $\ell$ tends to infinity. 
More precisely,  noting the equality :
\begin{equation}\label{ellintegral}
\left|\int_{\Delta}(1 - |t|^{2})^{\ell}dt\wedge d\bar{t}\,\right| = 
\frac{2\pi}{\ell +1}
\end{equation} 
(where $\Delta := \{ t\in\mathbb{C}; |t| < 1\}$), we see that 
\begin{equation}\label{(2)}
\parallel\rho\cdot\tau_{\ell}\parallel^{2} 
\sim 2^{2na\ell}\left(\frac{2\pi}{\ell}\right)^{n}
\end{equation}
as $\ell$ tends to infinity, where $\sim$ means that  the ratio of the both sides converges to $1$ as $\ell$ tends to infinity. 
We set 
\begin{equation}
\phi := n\rho\log \sum_{i=1}^{n}\mid z_{i}\mid^{2}. 
\end{equation}
We may and do assume that $\ell$ is sufficiently large so that
\begin{equation}
(\ell -1)\omega_{1,\varepsilon} +\sqrt{-1}\,\Theta_{h_{A}} +\sqrt{-1}\Theta_{h_{P}}
+ \sqrt{-1}\partial\bar{\partial}\phi > 0
\end{equation}
holds on $X$, where $\omega_{1,\varepsilon} = \sqrt{-1}\Theta_{h_{1,\varepsilon}}$. 
We note that $\bar{\partial}(\rho\cdot\tau_{\ell})$  vanishes 
on the polydisc of radius $r/3$ with center $p$ as above. 
Then by (\ref{(2)}), the $L^{2}$-norm 
\begin{equation}
\parallel\bar{\partial}(\rho\cdot\tau_{\ell})\parallel_{\phi,\varepsilon}
\end{equation}
of  $\bar{\partial}(\rho\cdot\tau_{\ell})$ 
with respect to $e^{-\phi}\cdot h_{A}\cdot h_{1,\varepsilon}^{\ell-1}$ and $\omega_{1,\varepsilon}$ 
satisfies the inequality: 
\begin{equation}\label{(3)}
\parallel\bar{\partial}(\rho\cdot\tau_{\ell})\parallel_{\phi,\varepsilon}^{2}\leqq C_{0}\cdot\left(\frac{3}{r}\right)^{2n+2}2^{na\ell}\left(\frac{2\pi}{\ell}\right)^{n}
\end{equation}
for every $\ell$, where $C_{0}$ is a positive constant independent of $\ell$, 
but it may depend on the point $p$.  
By  H\"{o}rmander's $L^{2}$-estimate applied to the adjoint line bundle of the hermitian line bundle:
\begin{equation}
(\{(\ell -1)a + (a-1)\}(K_{X}+D)+D,e^{-\phi}\cdot h_{A}\cdot (h_{P}\cdot h_{\sigma_{N}})^{a-1}\cdot h_{1}^{\ell-1}\cdot h_{\sigma_{D}}),
\end{equation}
 we see that for every 
sufficiently large $\ell$, there exists a $C^{\infty}$- solution $u$ of 
the equation: 
\begin{equation}
\bar{\partial}u = \bar{\partial}(\rho\cdot\tau_{\ell})
\end{equation}
such that 
\begin{equation}
u(p) = 0
\end{equation}
and  
\begin{equation}
\parallel u\parallel_{\phi,\varepsilon}^{2} \leqq \frac{2}{\ell}\parallel\bar{\partial}(\rho\cdot\tau_{\ell})\parallel_{\phi,\varepsilon}^{2}
\end{equation}
hold, where $\parallel\,\,\,\,\parallel_{\phi,\varepsilon}$'s  denote the $L^{2}$-norms
with respect to $e^{-\phi}\cdot h_{A}\cdot h_{1,\varepsilon}^{\ell-1}\cdot h_{\sigma_{D}}$ and $\omega_{1}$
respectively.  
Here we note that $(U,\omega_{1}|U)$ is not complete.   But since $U$ is 
quasiprojective, we may approximate $\omega_{Y}$ by the K\"{a}hler form 
$\omega_{1} + \epsilon \omega_{U}$, where $\omega_{U}$ is a complete
 K\"{a}hler form with Poincar\`{e} growth along $Y \backslash U$ and 
 $\epsilon > 0$ is a small positive number.   
Hence we may apply the H\"{o}rmander's $L^{2}$-estimate on $(U,\omega_{Y} +\epsilon \omega_{U})$.   And then letting $\epsilon$ tend to $0$, we may apply 
the $L^{2}$-estimate on $(U,\omega_{1}|U)$. 
Then $\rho\cdot\tau_{\ell}- u$
is a holomorphic section of $\ell a(K_{X}+D)+A$ such that 
\begin{equation}
(\rho\cdot\tau_{\ell}- u)(p)
= \tau_{\ell}(p)
\end{equation}
and
\begin{equation}
\parallel\rho\cdot\tau_{\ell}- u\parallel^{2}
 \leqq  \left(1+ C_{0}\cdot \left(\frac{3}{r}\right)^{2n+2}\cdot\sqrt{\frac{2}{\ell}}\right)\cdot 2^{2na\ell}\cdot\left(\frac{2\pi}{\ell}\right)^{n}.
\end{equation}
Hence by the inductive assumption (\ref{assumption}) and (\ref{extremal}), this implies that there exists a positive constant 
$C$ independent of $\ell$ such that  
\begin{equation}\label{lowp}
K_{\ell,1}(p) \geqq \left(1 - \frac{C}{\sqrt{\ell}}\right)\cdot \ell^{n}\cdot (2\pi)^{-n}\cdot C_{\ell-1}\cdot \left(h_{A}^{-1}\cdot h_{1,\varepsilon}^{-\ell}\right)(p)
\end{equation} 
holds, since the point norm of $\tau_{\ell}$ at $p$ (with respect to $h_{A}\cdot h_{1}^{\ell}$) is asymptotically equal to $2^{\ell n}$. 
We note that in the estimate (\ref{lowp}), $C$ may depend on $p$, because the 
radius $r > 0$ depdends on $p \in U$.  This may prevent to obtain the estimate 
independent of $p \in U$.   

But this difficulty can be overcome by the following lemma.  

\begin{lemma}\label{nearE}
There exists a positive constant $C_{\delta}$ depending on $\delta > 0$ such that 
\[
K_{\ell,1}\geqq C_{\delta}^{n}\cdot (\ell !)^{n}\cdot h_{0,\delta}^{-\ell}\cdot h_{A}^{-1}
\]
holds for every $\ell \geqq 1$. \fbox{}
\end{lemma}
\noindent{\em Proof.}  Let $x$ be a point on $X$, then there exists a local holomorphic frame $\mbox{\bf e}_{P}$ of $aP$,  a normal coordinate $(z_{1},\cdots ,z_{n})$ with respect to the K\"{a}hler form $\omega_{0,\delta}$ on a neighbourhood $U$ with center $x$ and a such that 
\[
h_{0,\delta}(\mbox{\bf e}_{P},\mbox{\bf e}_{P})(x) \leqq \left(\prod_{i=1}^{n}(1 - c|z_{i}|^{2})\right) + O(|z|^{3})
\]  
holds.  We may take $c > 0$ independent of $x$. 
Then by the successive $L^{2}$-estimates as above,  we see that there exists a positive 
number $C_{\delta}$ idependent of $x$ such that  
\[
K_{\ell,1}\geqq C_{\delta}^{n}\cdot (\ell !)^{n}\cdot h_{1,\delta}^{-\ell}\cdot h_{A}^{-1}
\]
holds for every $m$.     \fbox{} \vspace{3mm} \\ 

\noindent By Lemma \ref{nearE}, taking $\delta = \varepsilon/2$, there exists 
a neighbourhood $W$ of $E$ such that  
the esitimte
\begin{equation}\label{W}
K_{\ell,1} \geqq (\ell !)^{n}\cdot(2\pi )^{-\ell n}\cdot h_{1,\varepsilon}^{-\ell}\cdot h_{A}^{-1}
\end{equation}
holds on $W$ for every $\ell \geqq 1$. 

Then by induction on $\ell$, using (\ref{extremal}) and 
(\ref{assumption}), we see that there exists a positive constant $C^{\prime}$
and a positive intger $\ell_{1}$  such that $C/\sqrt{\ell_{1}}< 1$ and for every $\ell >  \ell_{1}$ 
\begin{equation}\label{induc}
K_{\ell,1}(p) \geqq C ^{\prime}\cdot\left(\left(\prod_{k=\ell_{1}}^{\ell}\left(1 -\frac{C}{\sqrt{k}}\right)\right)\cdot (\ell !)^{n}\cdot (2\pi )^{-\ell n}\cdot h_{A}^{-1}\cdot h_{1,\varepsilon}^{-\ell}\right)(p)
\end{equation}   
holds for every $p\in U$, thanks to the estimate (\ref{W}) on $W$. 

Since $p\in U$ is arbitrary, we have the following lower estimate. 

\begin{lemma}\label{lower}
\begin{equation}\label{loweres}
\limsup_{\ell\rightarrow\infty}\sqrt[\ell ]{(\ell !)^{-n}K_{\ell,1}} \geqq (2\pi)^{-n}\cdot h_{1}^{-1}
\end{equation}
holds on $U$.  \fbox{}  
\end{lemma}
\noindent{\em Proof}.  
\begin{equation}\label{loweres}
\limsup_{\ell\rightarrow\infty}\sqrt[\ell ]{(\ell !)^{-n}K_{\ell,1}} \geqq (2\pi)^{-n}\cdot h_{1,\varepsilon}^{-1}
\end{equation}
holds for every $0<\varepsilon <1$. 
Letting $\varepsilon$ tend to $0$, we have the lemma. \fbox{}\vspace{3mm} \\
 
Next we shall estimate $K_{\ell,1}$ from above as in \cite{canonical}. 
We set 
\begin{equation}\label{vol}
\mu (X,K_{X}+D) = n!\cdot a^{-n}\limsup_{\ell\rightarrow\infty}\frac{\dim H^{0}\left(X,\mathcal{O}_{X}(\ell a(K_{X}+D))\right)}{\ell^{n}}. 
\end{equation}
We call $\mu(X,K_{X}+D)$ the {\bf volume} of $X$ with respect to $K_{X}+D$.

\noindent We set for every $\ell \geqq 1$
\begin{equation}\label{volumem}
dV_{\ell}:= K_{\ell,1}^{\frac{1}{\ell}}\cdot (h_{P}\cdot h_{\sigma_{N}})^{a-1}\cdot h_{\sigma_{D}}\cdot h_{A}^{\frac{1}{\ell}}.
\end{equation} 
Then $dV_{\ell}$ is a singular volume form on $X$.  
But by the construction, $\pi_{V}^{*}dV_{\ell}$ is a locally bounded volume form on $V$ for 
any $\pi_{V}: \Delta^{n}\to V$ as (\ref{piu}) and (\ref{piu2}).
By using H\"{o}lder's inequality, we have the following lemma.
\begin{lemma}\label{holder}
\begin{equation}\label{upperbd}
\limsup_{\ell\rightarrow\infty}\,\,(\ell !)^{-\frac{n}{\ell}}\int_{X}dV_{\ell}  \leqq a^{n}\mu(X,K_{X}+D)
\end{equation}
holds.  
\end{lemma}
{\em  Proof}. 
The following proof is essentially the same as \cite[Lemma 3.2]{canonical}.
By H\"{o}lder's inequality we have that 
\begin{equation}\label{holder2}
\int_{X}dV_{\ell} 
\leqq \left(\int_{X}\left(\frac{dV_{\ell}}{dV_{\ell-1}}\right)^{\ell}\!\! dV_{\ell-1}\right)^{\frac{1}{\ell}}\cdot \left(\int_{X}dV_{\ell-1}\right)^{\frac{\ell-1}{\ell}}
\end{equation}
holds for every $m\geqq 2$.
By direct calculation, we see that  
\begin{equation}\label{predim}
\left(\frac{dV_{\ell}}{dV_{\ell-1}}\right)^{\ell}\!\! dV_{\ell-1}
=  K_{\ell,1}\cdot h_{\ell-1}\cdot(h_{P}\cdot h_{\sigma_{N}})^{a-1}\cdot h_{\sigma_{D}}\cdot h_{A} 
\end{equation}
holds.   Hence by the definition of the Bergman kernel $K_{\ell,1}$ and (\ref{predim}) we have that  
\begin{equation}\label{dim}
\int_{X}\left(\frac{dV_{\ell}}{dV_{\ell-1}}\right)^{\ell}dV_{\ell-1} 
\leqq \dim |\ell a(K_{X}+D) + A| + 1
\end{equation}
holds. Combining (\ref{holder2}) and (\ref{dim}), we have that 
\begin{equation}
\int_{X}dV_{\ell} 
\leqq \left(\dim |\ell a(K_{X}+D) + A| + 1\right)^{\frac{1}{\ell}}\cdot \left(\int_{X}dV_{\ell-1}\right)^{\frac{\ell-1}{\ell}}
\end{equation}
holds. 
Repeating the same estimate (if $m \geqq 3$), 
\begin{equation}
\int_{X}dV_{\ell} 
\leqq \left(\dim |\ell a(K_{X}+D) + A| + 1\right)^{\frac{1}{\ell}}\cdot
\left(\dim |(\ell-1)a(K_{X}+D) + A| + 1\right)^{\frac{1}{\ell}}\cdot \left(\int_{X}dV_{\ell-2}\right)^{\frac{\ell-2}{\ell}}
\end{equation}
holds. 
Continueing this process, we obtain the inequality: 
\begin{equation}\label{upperind}
\int_{X}dV_{\ell} 
\leqq \left(\prod_{k=1}^{\ell}(\dim |A+ka(K_{X}+D)| + 1)\right)^{\frac{1}{\ell}}\cdot \left(\int_{X}dV_{1}\right)^{\frac{1}{\ell}}. 
\end{equation}
 Then by (\ref{upperind}) and (\ref{vol}),
\begin{equation}\label{upperbd}
\limsup_{\ell\rightarrow\infty}(\ell !)^{-\frac{n}{\ell}}\int_{X} dV_{\ell} \leqq a^{n}\mu(X,K_{X}+D)
\end{equation}
holds.  This completes the proof of Lemma \ref{holder}. q.e.d. \vspace{5mm} \\  
By the assumption that $K_{X}+D$ is big, we see that $\mu (X,K_{X}+D)$ 
is positive.  Moreover by the existence of the Zariski decomposition (\ref{ZD3}),  we see that 
\begin{equation}\label{P}
\mu (X,K_{X}+D) = P^{n}
\end{equation}
holds.  

On the other hand,  since $h_{1}$ is an AZD of $a(K_{X} + D)$, we see that the absolutely continuous part $\omega_{1,abc}$ of $\omega_{1}$ is cohomologous
to $2a\pi c_{1}(P)$ and we have that 
\begin{equation}\label{intabc}
\int_{X}\omega^{n}_{1,abc} = (2\pi)^{n}a^{n}P^{n}
\end{equation}
holds.    Hence by (\ref{P}) and (\ref{intabc}), we have that 
\begin{equation}\label{growth}
\int_{X}\omega^{n}_{1,abc}  = (2\pi)^{n}a^{n}\mu (X,K_{X}+D)
\end{equation}
holds. 
Then by (\ref{upperbd}) and (\ref{growth}), we see that 
\begin{equation}\label{upper2}
\limsup_{\ell\to\infty}(\ell !)^{-\frac{n}{\ell}}\int_{X}dV_{\ell} 
\leqq \frac{1}{(2\pi)^{n}}\int_{X}\omega^{n}_{1,abc}
\end{equation}
holds.    
By Lemma \ref{lower} and (\ref{key}), we see that 
\begin{equation}\label{lower2}
\limsup_{\ell\to\infty}(\ell !)^{-\frac{n}{\ell}}dV_{\ell} 
\geqq (2\pi)^{-n}h_{1}^{-1}\cdot (h_{P}\cdot h_{\sigma_{N}})^{-(a-1)}\cdot h_{\sigma_{D}}^{-1} = (2\pi)^{-n}\omega_{1}^{n} 
\end{equation}
hold. 
We note that $(\ell !)^{-\frac{n}{\ell}}K_{\ell,1}^{\frac{1}{\ell}}$ is 
the $\ell$-th root of sum of absolute value squares of holomorphic sections 
and is considered to be pluri-subharmonic with respect to 
a (hence any) local holomorphic trivialization  of $A + \ell a(K_{X}+D)$.  
Then by the definition of $dV_{\ell}$ (cf.(\ref{volumem})) and (\ref{upper2}), we see that $(\ell !)^{-\frac{n}{\ell}}\pi_{V}^{*}dV_{\ell}$ is locally uniformly bounded and semipositive volume form on $\Delta^{n}$ 
for any local cyclic branched covering $\pi_{V} : \Delta^{n}\to V$ as (\ref{piu}) and (\ref{piu2})
in Section 2.   Hence by Lebesgue's bounded convergence theorem,  
\begin{equation}\label{Lebesugue}
\int_{X}\limsup_{\ell\to\infty}(\ell !)^{-\frac{n}{\ell}}dV_{\ell}
= \limsup_{\ell\to\infty}(\ell !)^{-\frac{n}{\ell}}\int_{X}dV_{\ell} 
\end{equation}
holds. 
Hence by (\ref{upper2}) and (\ref{lower2}), we have that
\begin{equation}\label{equality}
\limsup_{\ell\to\infty}(\ell !)^{-\frac{n}{\ell}}dV_{\ell} = \frac{1}{(2\pi)^{n}}\omega_{1,abc}^{n}
\end{equation} 
holds almost everywhere on  $X$.
Hence by (\ref{equality}), (\ref{volumem}) and (\ref{key}), we see that 
\begin{equation}
\limsup_{\ell\to\infty}\sqrt[\ell]{(\ell !)^{-n}K_{\ell,1}}
= \left(\frac{1}{(2\pi)^{n}}\omega_{1,abc}^{n}\right)\cdot(h_{P}\cdot h_{\sigma_{N}})^{-(a-1)}\cdot h_{\sigma_{D}}^{-1} =\frac{n!}{(2\pi)^{n}}\,h_{1}^{-1}
\end{equation}
hold almost everywhere on $X$. 
This implies the following lemma.   
\begin{lemma}\label{K1}
\begin{equation}
K_{1}:= \mbox{\em the upper-semi-continuous envelope of}\,\,\limsup_{\ell\rightarrow\infty}\sqrt[\ell ]{(\ell !)^{-n}K_{\ell,1}}
\end{equation}
exists and  
\[
h_{1}:= \frac{(2\pi)^{n}}{n!}\,K_{1}^{-1}
\]
is an AZD of $a(K_{X}+D)$.
 If we set 
\begin{equation}
\omega_{1}:= \sqrt{-1}\Theta_{h_{1}},
\end{equation}
then $\omega_{1}$ is a closed positive current on $X$ and satisfies the equation:
\begin{equation}\label{omega1}
-\mbox{\em Ric}_{\omega_{1}} + \frac{a-1}{a}\omega_{0} + 2\pi[D] = \omega_{1}. 
\end{equation}
on $X$ and 
\[
h_{1} = \left(\frac{\omega_{1,abc}^{n}}{n!}\right)^{-1}\cdot (h_{P}\cdot h_{\sigma_{N}})^{a-1}\cdot h_{\sigma_{D}}
\]
holds on $X$.  
\fbox{}
\end{lemma}
Next replacing $h_{P},\omega_{0}$ by $h_{1},\omega_{1}$ respectively. 
We obtain the dynamical system of Bergman kernels 
$\{ K_{\ell,2}\}_{\ell=1}^{\infty}$ and 
\begin{equation}
K_{2}:= \limsup_{\ell\rightarrow\infty}\sqrt[\ell ]{(\ell !)^{-n}K_{\ell,2}}
\end{equation}
exists and  $h_{2}:= K_{2}^{-1}$ is a $C^{\infty}$-hermitian metric on 
$a(K_{X}+D)|U$ and is a singular hermitian metric on $a(K_{X}+D)$.  If we set 
\begin{equation}
\omega_{2}:= \sqrt{-1}\,\Theta_{h_{2}},
\end{equation}
then $\omega_{2}$ is a K\"{a}hler form on $X$ and satisfies the equation:
\begin{equation}
-\mbox{Ric}_{\omega_{2}} + \frac{a-1}{a}\omega_{1} + 2\pi[D]
= \omega_{2}. 
\end{equation}
Inductively, for every positive integer $m$, we define the dynamical system of Bergman kernels $\{ K_{\ell,m}\}_{\ell=1}^{\infty}$.
Continueing this process, we obtain the following theorem.

\begin{theorem}\label{doubleinduction2}
Let $\{ K_{\ell,m}\}_{\ell=1}^{\infty}$ be the system of Bergman kernels 
defined as above.  If we define 
\begin{equation}\label{km}
K_{m}:= \mbox{\em the upper-semi-continuous envelope of}\,\,\, 
\limsup_{\ell\rightarrow\infty}\sqrt[\ell]{(\ell !)^{-n}K_{\ell,m}},
\end{equation}
then 
\begin{equation}
h_{m} := \frac{n!}{(2\pi)^{n}}\,K_{m}^{-1}
\end{equation}
is a singular hermitian metric on $a(K_{X}+D)$ with strictly positive curvature on $U$ and is an AZD of $a(K_{X} + D)$. 
And if we define the closed positive current $\omega_{m}$ on $X$  by 
\begin{equation}
\omega_{m}:= \sqrt{-1}\,\Theta_{h_{m}},
\end{equation}
then $\{\omega_{m}\}_{m=1}^{\infty}$ is $C^{\infty}$ on $U$ and satisfy  the equations:
\begin{equation}
-\mbox{\em Ric}_{\omega_{m}} + \frac{a-1}{a}\,\omega_{m-1} + 2\pi [D]= \omega_{m} 
\end{equation}
on $U$ for every $m \geqq 1$ as {\rm (\ref{ricciind2})}. 
Hence by Theorem \ref{iteration2},
\[
\omega_{K}:= \frac{1}{a}\lim_{m\rightarrow\infty}\omega_{m}
\]
exists in the sense of current on $X$  and in the compact uniform $C^{\infty}$-topology on $U$. And  $\omega_{K}$ is the unique solution of the equation:
\[
-\mbox{\em Ric}_{\omega_{K}} + 2\pi[D] = \omega_{K},
\]
on $U$  such that 
\[
h_{K} := \left(\frac{1}{n!}\,\omega_{K,abc}^{n}\right)^{-1}\cdot h_{\sigma_{D}}
\]
is an AZD of $K_{X} + D$ (Here $h_{\sigma_{D}}$ appears just to identify a multivalued holomorphic section of 
$K_{X} + D$ with a multivalued meromorphic $(n,0)$-form with pole along $\mbox{Supp}\, D$.), i.e.
$\omega_{K}$ is the canonical K\"{a}hler-Einstein current on 
$(X,D)$. 
\fbox{}
\end{theorem} 

Now we shall consider the uniqueness of canonical K\"{a}hler-Einstein currents 
on KLT pairs.  To state the result we shall introduce an equivalence relation 
between KLT pairs. 
\begin{definition}\label{eqKLT}
Let $(X_{1},D_{1}),(X_{2},D_{2})$ be  KLT (resp. LC) pairs.   
We say that $(X_{i},D_{i}) (i= 1,2)$ are birationally equivalent, 
if there exists a KLT(resp. LC) pair $(Y,D_{Y})$ and compositions of blowing ups with smooth centers: $\mu_{i} : Y \to X_{i}(i=1,2)$
\begin{equation}
\begin{picture}(400,400)
\setsqparms[1`1`1`1;350`350]
\putAtriangle(0,00)%
[Y`X_{1} ` X_{2}; \mu_{1}`\mu_{2} `\mu_{2}\circ\mu_{1}^{-1}]
\end{picture}
\end{equation}  
 such that 
\[
K_{Y} + D_{Y}- \mu_{i}^{*}(K_{X} + D_{i})\,\,\,\,\, (i = 1,2)
\]
are effective exceptional divisors respectively. 
Here $\mu_{2}\circ\mu_{1}^{-1}$ is a birational rational map \fbox{}
\end{definition}
Using Theorem \ref{doubleinduction2}, we obtain the following uniqueness of 
canonical K\"{a}hler-Einstein currents on KLT pairs.
  
\begin{theorem}\label{unique}
Let $X$ be a smooth projective variety and let $D$ be an effective $\mathbb{Q}$-divisor on $X$ such that $(X,D)$ is a KLT pair of log general type.  
Then  canonical K\"{a}hler-Einstein current $\omega_{K}$ on the KLT pair $(X,D)$ is unique.  

Moreover if two such KLT pairs $(X_{i},D_{i}) (i= 1,2)$ are equivalent in the sense of Definition \ref{eqKLT}, then for any KLT pair $(Y,D_{Y})$ and 
the morphisms $\mu_{i} : Y\to X_{i}(i=1,2)$ as in Definition \ref{eqKLT}, 
we have that for the canonical K\"{a}hler-Einstein current $\omega_{K}$ on 
$(Y,D_{Y})$, $(\mu_{i})_{*}\omega_{K}(i = 1,2)$ are unique canonical K\"{a}hler-Einstein currents on $(X_{i},D_{i})(i=1,2)$ respectively in the sense that 
the absolutely continuous part $((\mu_{i})_{*}\omega_{K})_{abc}$ coincides that of the canonical K\"{a}hler-Einstein currents $\omega_{i}$ of $(X_{i},D_{i})$ respectively.  \fbox{} 
\end{theorem}
\noindent{\em Proof}. 
As in Section \ref{ex1}, we take a log resolution
$\mu : Y \to X$ 
of $(X,D)$ which satisfies  the followings :  
\begin{enumerate}
\item[(1)] If we write $K_{Y} = \mu^{*}(K_{X}+D) + \sum a_{i}E_{i}$, 
where $\{E_{i}\}$ are prime divisors. 
then $a_{i} > -1$ holds for every $i$. 
\item[(2)] There exists a Zariski decomposition: 
of $\mu^{*}(K_{X} +D) = P + N$  of $\mu^{*}(K_{X}+D)$ as (\ref{ZD}). 
\end{enumerate} 
We set  $I := \{ i|\, a_{i} < 0\}$. 
Then  replacing $X$ by $Y$ and $D$ by
\begin{equation}\label{DY}
D_{Y} := \sum_{i\in I} (-a_{i})E_{i}, 
\end{equation}
we obtain the new KLT pair $(Y,D_{Y})$ such that 
\begin{enumerate}
\item[(a)] $R(Y,K_{Y} + D_{Y}) \simeq R(X,K_{X} +D)$, 
\item[(b)] There exists a Zariski decomposition:  
$K_{Y} + D_{Y} = P + (N + \sum_{i\not{\in} I} a_{i}E_{i})$.
\end{enumerate}
Suppose that  Theorem \ref{unique} holds for $(Y,D_{Y})$.  
Then if there exist canonical K\"{a}hler-Einstein currents $\omega_{K,1},\omega_{K,2}$ on $(X,D)$, we see that $(\mu^{*}\omega_{K,1})_{abc} = (\mu^{*}\omega_{K,2})_{abc}$ 
holds.  Hence $\omega_{K,1} = \omega_{K,2}$ holds on $X$. 

Hence to prove Theorem \ref{unique}, we may assume that 
the Zariski decomposition $K_{X} + D = P + N$ exists on  $X$ from the beginning  and $\mbox{Supp}\, D$ and $\mbox{Supp}\, N$ are divisors with normal crossings. 

We set $n:= \dim X$.  
By the construction of $\{ K_{\ell,m}\}$ as above, we see that 
\[
K_{1}: = \limsup_{\ell\to\infty}\sqrt[\ell]{(\ell !)^{-n}K_{\ell,1}} 
\]
does not depend on $(A,h_{A})$ and if we set  
\[
h_{1}:= \mbox{the lower-semi-continuous envelope of}\,\,\, K_{1}^{-1}
\]
is an AZD of $a(K_{X} + D)$. 
On the other hand let $\omega_{1}^{\prime}$ be a solution of 
\[
-\mbox{Ric}_{\omega^{\prime}_{1}} + \frac{a-1}{a}\omega_{0} + 2\pi [D] = \omega_{1}^{\prime}. 
\]
such that 
\[
h^{\prime}_{1}:= \left(\frac{1}{n!}(\omega_{1,abc}^{\prime})^{n}\right)^{-1}\cdot (h_{P}\cdot h_{\sigma_{N}})^{a-1}\cdot h_{\sigma_{D}}
\]
is an AZD of $a(K_{X}+D)$. 
Then by the proof of Theorem \ref{doubleinduction2} above, we see that 
\begin{equation}
\omega_{1} = \sqrt{-1}\Theta_{h_{1}} = \sqrt{-1}\Theta_{h^{\prime}_{1}}
\end{equation}
holds. Hence (\ref{omega1}) has a unique solution.  
If $a = 1$, then we have nothing to prove.  
If $a  > 1$, then by the equation (\ref{ricciind2}),  we have that 
$\omega_{K}$ is unique up to the choice of $h_{P}$.  

Now we shall prove that $\omega_{K}$ does not depend on the choice of 
$h_{P}$.  This follows from the contraction property of the Ricci iterations.
 
Let us recall the equation (\ref{MAind3}).  Let 
$\{ u_{m,\delta}\}_{m=0}^{\infty}, \{ u_{m,\delta}^{\prime}\}_{m=0}^{\infty}$ be the systems of solutions 
corresponding to the metrics $h_{P}$ and $h^{\prime}_{P}$ respectively.  
Then by (\ref{MAind3}), we have that 
\begin{equation}\label{succ3}
\log \frac{(a\omega_{m,\delta,0} +\sqrt{-1}\partial\bar{\partial}u_{m,\delta})^{n}}
{(a\omega_{m-1,\delta,0} +\sqrt{-1}\partial\bar{\partial}u^{\prime}_{m-1,\delta})^{n}}  = \ (u_{m,\delta} - u^{\prime}_{m-1,\delta}) - \frac{a - 1}{a}(u_{m-1,\delta}-u^{\prime}_{m-2,\delta}) 
\end{equation}
hold for $m\geqq 1$. 
Let $\epsilon$ be a sufficiently small positive number and we define  $u_{m,\delta}(\epsilon),u^{\prime}_{m,\delta}(\epsilon)$ as (\ref{umdeltaepsilon}).  
Then repeating the same argument as in Section 3.2 (cf. (\ref{indlower})) by the maximum principle, 
we see that there exists a positive constants $C_{(-)}$ independent of $m$ and $\delta$
such that 
\begin{equation}\label{inf5}
u_{m,\delta}(\epsilon) - u^{\prime}_{m-1,\delta}(\epsilon)
-\alpha_{m}(\delta +\epsilon)\log\parallel\sigma_{E}\parallel^{2}_{h_{E}} 
\geqq   -C_{(-)}\left(\frac{a-1}{a}\right)^{m-1}
\end{equation}
holds  for every $m\geqq 1$, where $\alpha_{m}$ is as (\ref{alpham}).  
Switching $u_{m,\delta}(\epsilon)$ and $u_{m,\delta}^{\prime}(\epsilon)$, 
we see that there exists a positive constant $C_{+}$ independent of $m$ and $\delta$ such that 
\begin{equation}\label{sup5}
u_{m-1,\delta}(\epsilon) -  u^{\prime}_{m,\delta}(\epsilon) + \alpha_{m}(\delta +\epsilon)\log\parallel\sigma_{E}\parallel^{2}_{h_{E}} \leqq C_{(+)}\left(\frac{a-1}{a}\right)^{m-1}
\end{equation}
holds.  In Section 3.2, we have already observed that 
$\{ u_{m,\delta}\}_{m=0}^{\infty}$ and $\{ u_{m,\delta}^{\prime}\}_{m=0}^{\infty}$
 converge.  Hence  by the estimates (\ref{inf5}) and (\ref{sup5})
\begin{equation}
\lim_{m\to\infty} (u_{m,\delta}-u_{m,\delta}^{\prime}) = 0
\end{equation}
holds on $U$. 
  Now we note that the estimates (\ref{inf5}) and (\ref{sup5}) are 
independent of $\delta$.   
Then letting $\delta$ tend to $0$, by (\ref{um}) we have that 
$\omega_{K} = \lim_{m\rightarrow\infty}\omega_{m}$ is independent of the choice of $h_{P}$.
 
We note that if we take $h_{P}$ to be $h_{K}$, then the Ricci iteration is 
stable, i.e., $\omega_{m} = a\omega_{K}$ holds for every $m\geqq 1$\footnote{Here $\omega_{K}$ is not a pull back of the Fubini-Study K\"{a}hler form up to a constant multiple like $\omega_{P} = \sqrt{-1}\Theta_{h_{P}}$, but the Ricci iteration above and the above argument also work in this case}.

Hence this completes the proof of Theorem \ref{unique}. 
q.e.d.

\subsection{Weak semistability}\label{ws}
First we shall recall  several definitions  in \cite{v}. 
To measure the positivity of coherent sheaves, we shall introduce the 
following notion.
\begin{definition}\label{weakpositivity}
Let $Y$ be a quasi-projective reduced scheme, $Y_{0}\subseteq Y$ an open dense subscheme and let $\mathcal{G}$ be locally free sheaf on $Y$, of finite constant  rank.  Then $\mathcal{G}$ is {\bf weakly positive} over $Y_{0}$, if 
for an ample invertible sheaf $\mathcal{H}$ on $Y$ and for a given number 
$\alpha > 0$ there exists some $\beta > 0$ such that 
$S^{\alpha\cdot\beta}(\mathcal{G})\otimes \mathcal{H}^{\beta}$ 
is globally generated over $Y_{0}$.  \fbox{}
\end{definition}  
The notion of weak positivity is a natural generalization of the notion of nefness of line bundles.  Roughly speaking, the weak semipositivity of $\mathcal{G}$ over $Y_{0}$ means that 
$\mathcal{G}\otimes \mathcal{H}^{\varepsilon}$ is $\mathbb{Q}$-globally generated over $Y_{0}$ for every $\varepsilon > 0$. 
\begin{definition}\label{succeq}
Let $\mathcal{F}$ be a locally free sheaf and let $\mathcal{A}$ be an invertible sheaf, both on a quasi-projective reduced scheme $Y$.  
We denote 
\begin{equation}
\mathcal{F} \succeq \frac{b}{a}\,\,\mathcal{A},
\end{equation}
if $S^{a}(\mathcal{F})\otimes \mathcal{A}^{-b}$ is weakly positive over $Y$,
where $a,b$ are positive integers.  
\fbox{}
\end{definition}
Let $X$ be a normal variety.  We define the canonical sheaf $\omega_{X}$ of $X$  by 
\begin{equation}
\omega_{X} := i_{*}\mathcal{O}_{X_{reg}}(K_{X_{reg}}), 
\end{equation}
where $i : X_{reg}\hookrightarrow X$ is the canonical injection. 
The following notion introduced by Viehweg is closely related to the 
notion of logcanonical thresholds. 
\begin{definition}
Let $(X,\Gamma)$ be a pair of normal variety $X$ and an effective Cartier divisor $\Gamma$.  Let $\pi : X^{\prime}\to X$ be a log resolution of $(X,\Gamma)$ and let $\Gamma^{\prime}:= \pi^{*}\Gamma$. 
For a positive integer $N$ we define
\begin{equation}
\omega_{X}\left\{\frac{-\Gamma}{N}\right\}
= \pi_{*}\left(\omega_{X^{\prime}}\left(-\left\lfloor\frac{\Gamma^{\prime}}{N}\right\rfloor\right)\right)
\end{equation}
and 
\begin{equation}
\mathcal{C}_{X}(\Gamma,N)= \mbox{\em Coker}\left\{\omega_{X}\left\{\frac{-\Gamma}{N}\right\}\to \omega_{X}\right\}.
\end{equation}
If $X$ has at most rational singularities, one defines :
\begin{equation}
e(\Gamma)= \min\{N>0\,|\,\mathcal{C}_{X}(\Gamma,N) = 0\}.
\end{equation}
If $\mathcal{L}$ is an invertible sheaf,  $X$ is proper with at most rational  singularities and  $H^{0}(X,\mathcal{L})\neq 0$, then one defines
\begin{equation}
e(\mathcal{L}) = \sup\left\{ e(\Gamma)|\Gamma :\,\mbox{effective Cartier divisor 
with $\mathcal{O}_{X}(\Gamma)\simeq \mathcal{L}$}\right\}.
\end{equation}
\fbox{}
\end{definition}
Now we state the result of E. Viehweg. 

\begin{theorem}(\cite[p.191,Theorem 6.22]{v})\label{viehweg}
Let $f : X\to Y$ be a flat surjective projective Gorenstein morphism of reduced connected quai-projective schemes.  Assume that the  relative canonical sheaf $\omega_{X/Y}:= \omega_{X}\otimes f^{*}\omega_{Y}^{-1}$  
is $f$-semi-ample and that the fibers $X_{y} = f^{-1}(y)$ are reduced normal varieties with at most rational singularities.  Then one has :
\begin{enumerate}
\item[\em (1)]{\bf Functoriality}: For $m > 0$ the sheaf $f_{*}\omega_{X/Y}^{m}$ is locally free of 
rank $r(m)$ and it commutes with arbitrary base change.
\item[\em (2)]{\bf Weak semipositivity}: For $m > 0$ the sheaf $f_{*}\omega_{X/Y}^{m}$ is weakly positive over $Y$.
\item[\em (3)] {\bf Weak semistability}: Let $m > 1, e > 0$ and $\nu > 0$ be chosen so that 
$f_{*}\omega_{X/Y}^{m}\neq 0$ and 
\begin{equation}
e \geqq \sup\left\{\frac{k}{m-1}, e(\omega_{X_{y}}^{k})\,\,;\,\,\mbox{for}\,\,y\in Y\right\}
\end{equation}
hold. 
Then 
\begin{equation}\label{wst}
f_{*}\omega_{X/Y}^{m} \succeq \frac{1}{e\cdot r(k)}\,\det (f_{*}\omega_{X/Y}^{k})
\end{equation}
holds. \fbox{}
\end{enumerate}
\end{theorem}

\noindent The weak semistability of $f_{*}\omega_{X/Y}^{m}$ is very important in the application of the weak semipositivity.   For example, if $\dim Y = 1$ 
and $\deg \det (f_{*}\omega_{X/Y}^{m}) > 0$ hold, then for a sufficiently large
 $r$  $S^{r}(f_{*}\omega_{X/Y}^{m})$ is globally generated. 
Moreover by the finite generation of canonical rings (\cite{b-c-h-m}),
for every sufficiently large $m$, $f_{*}\omega_{X/Y}^{m!}$ is globally generated on $Y$.
\subsection{Weak semistability of the direct image of the relative log canonical bundle of a family of KLT pairs}
Now we state our generalization of Theorem \ref{viehweg}. 
\begin{theorem}\label{logmain}
Let $f : X \to Y$ be an algebraic fiber space and let
$D$ be an effective $\mathbb{Q}$ divisor on $X$ 
such that $(X,D)$ is KLT.  Let $Y^{\circ}$ denote the complement of the 
discriminant locus of $f$. 
We set  
\begin{equation}\label{y00}
Y_{0}:= \{y\in Y|y\in Y^{\circ}, (X_{y},D_{y})\,\,\mbox{is a KLT pair}\}.  
\end{equation}
\begin{enumerate}
\item[\em (1)]{\bf Weak semistability 1}:  Let $r$ denote $\mbox{\em rank}\,f_{*}\mathcal{O}_{X}(\lfloor m(K_{X/Y}+D)\rfloor)$. 
Let $X^{r}:= X\times_{Y}X\times_{Y}\cdots\times_{Y}X$ be the $r$-times fiber product over $Y$ and let $f^{r}: X^{r}\to Y$ be the natural morphism. 
And let $D^{r}$ denote the divior on $X^{r}$ defined by 
$D^{r} = \sum_{i=1}^{r} \pi_{i}^{*}D$,  
where $\pi_{i} : X^{r} \longrightarrow X$ denotes the projection: 
$X^{r}\ni (x_{1},\cdots ,x_{n}) \mapsto x_{i} \in X$.  

There exists a canonically defined effective divisor $\Gamma$ (depending on $m$) on $X^{r}$ which does not conatin any fiber $X^{r}_{y} (y\in Y^{\circ})$ such that if we we define the number $\varepsilon_{0}$ by 
\begin{equation}\label{delta00}
\varepsilon_{0} := \sup \{\varepsilon\, |\, (X^{r}_{y},D^{r}_{y}+\varepsilon\Gamma_{y})\,\, \mbox{is KLT for all $y\in Y^{\circ}$}\}, 
\end{equation}
then for every $0 < \varepsilon < \varepsilon_{0}$,
  there exists a singular hermitian metric $h_{m,\varepsilon}$ on 
\begin{equation}
\left(\frac{1}{m}+\varepsilon\right)r(\lfloor m(K_{X/Y}+D)\rfloor) - \varepsilon\cdot f^{*}\det \left(f_{*}\mathcal{O}_{X}(\lfloor m(K_{X/Y}+D)\rfloor)\right)^{**}
\end{equation}
such that 
\begin{enumerate}
\item $\sqrt{-1}\,\Theta_{h_{m,\varepsilon}}\geqq 0$ holds on $X$ in the sense 
of current.
\item  
 For every $y\in Y_{0}$, $h_{m,\varepsilon}|X_{y}$ is well defined 
and is an AZD of
\[ 
\left(\frac{1}{m}+\varepsilon\right)r(\lfloor m(K_{X/Y}+D)\rfloor) - \varepsilon\cdot f^{*}\det \left(f_{*}\mathcal{O}_{X}(\lfloor m(K_{X/Y}+D)\rfloor)\right)^{**}|X_{y}.
\]
\end{enumerate} 
\item[(2)]{\bf Weak semistability 2}: In addition, suppose that $K_{X/Y}+D$ is $\mathbb{Q}$-linear equivalent to a Cartier divisor $G$ on $X$. Then for every $m\geqq 1$ 
and a rational number $0 < \varepsilon < \varepsilon_{0}$, there exists an
ample line bundel $A$ on $Y$ such that for every 
sufficiently large positive integer $\ell$ with $\ell\varepsilon\in \mathbb{Z}_{>0}$,  
\[
f_{*}\mathcal{O}_{X}(\ell mG)\otimes \left(\det \left(f_{*}\mathcal{O}_{X}(mG)\right)^{**}\right)^{-\ell\varepsilon}\otimes\mathcal{O}_{Y}(A)
\]  
is globally generated over $Y_{0}$. And for every suffficiently large 
$\ell$ 
\begin{equation}\label{WS}
f_{*}\mathcal{O}_{X}(\ell mG)\succeq \frac{\ell m\varepsilon}{1+ (1+m\varepsilon)\ell r}\,\left(\det \left(f_{*}\mathcal{O}_{X}(mG)\right)^{**}\right)^{\otimes\ell}
\end{equation}
holds over $Y_{0}$.
\fbox{} 
\end{enumerate}
\end{theorem}
The main advantage of Theorem \ref{logmain} to  Theorem \ref{viehweg} is that 
we do not assume the $f$-semiampleness of the relative log canonical bundles.

\subsection{Proof of Theorem \ref{logmain}}
Let us prove Theorem \ref{logmain}. 
The proof follows closely the one of Theorem \ref{viehweg} in \cite{v}. 
But we replace the use of branched coverings  in \cite{v} by the use of Theorem \ref{variation4}.  This enables us to get rid of the assumption that 
$K_{X/Y}$ is $f$-semiample.  

Let us start the proof.  
Let $f : X \longrightarrow Y$ be an algebraic fiber space.
Let $Y_{0}$ be the Zariski open subset of $Y$ defined as (\ref{y00}) above. 
For $r = \mbox{rank}\,f_{*}\mathcal{O}_{X}(\lfloor m(K_{X/Y}+D)\rfloor)$
we set  $X^{r}:= X\times_{Y}X\times_{Y}\cdots\times_{Y}X$ be the $r$-times fiber product over $Y$ and let $f^{r}: X^{r}\to Y$ be the natural morphism. 
Then we have the natural morpshim 
\begin{equation}
\det \left(f_{*}\mathcal{O}_{X}(\lfloor m(K_{X/Y}+D)\rfloor)\right)^{**} \rightarrow 
\otimes^{r}f_{*}\mathcal{O}_{X}(\lfloor m(K_{X/Y}+D)\rfloor) = f^{r}_{*}\mathcal{O}_{X^{r}}(\lfloor m(K_{X^{r}/Y} + D^{r})\rfloor).    
\end{equation}
Hence we have the canonical global section:  
\begin{equation}
\gamma \in H^{0}\left(X^{r},f^{r*}(\det \left(f_{*}\mathcal{O}_{X}(\lfloor m(K_{X/Y}+D)\rfloor\right)^{**}))^{-1}
\otimes \mathcal{O}_{X^{r}}(\lfloor m(K_{X^{r}/Y}+D^{r})\rfloor)\right). 
\end{equation}
Let $\Gamma$ denote the zero divisor of $\gamma$.  
It is clear the $\Gamma$ does not contain any fiber over $Y_{0}$. 
Now we define $\delta_{0}$ as (\ref{delta00}). 
Let us take a positive rational number $\varepsilon < \delta_{0}$.
We set 
\begin{equation}\label{theta}
\Theta :=  \frac{1}{m}\lfloor mD^{r}\rfloor + \varepsilon\Gamma.
\end{equation}
Then there exists the relative  canonical measure
$d\mu_{can,(X^{r},\Theta)}$ on 
$f: (X^{r},\Theta)\longrightarrow Y$ as in Theorems \ref{main2} and \ref{variation4}. 
By the logarithmic pluri-subharmonicity  of the relative canonical measure (Theorem \ref{variation4}),
we see that 
\begin{equation}
\sqrt{-1}\partial\bar{\partial}\log d\mu_{can,(X^{r},\Theta)/Y} \geqq 0
\end{equation}
holds on $X$ in the sense of current.
We set 
\begin{equation}
H_{m,\varepsilon} := d\mu_{can,(X^{r},\Theta)/Y}^{-1}
\end{equation}
Then $H_{m,\varepsilon}$ is a singular hermitian metric on 
\begin{equation}\label{me}
\left(\frac{1}{m} + \varepsilon\right)(\lfloor m(K_{X^{r}/Y}+D^{r})\rfloor) - \varepsilon\cdot f^{r*}\det \left(f_{*}\mathcal{O}_{X}(\lfloor m(K_{X/Y}+D)\rfloor\right)^{**})
\end{equation}
with semipositive curvature current by Theorem \ref{variation4} 
and $H_{m,\varepsilon}|X^{r}_{y}$ is an AZD of 
\begin{equation}
\left(\frac{1}{m} + \varepsilon\right)(\lfloor m(K_{X^{r}/Y}+D^{r})\rfloor)|X^{r}_{y} - \varepsilon\cdot f^{r*}\det \left(f_{*}\mathcal{O}_{X}(\lfloor m(K_{X/Y}+D)\rfloor)\right)^{**}|X^{r}_{y}
\end{equation}
for every $y\in Y_{0}$.  
Let $h_{m,\varepsilon}$ be the restriction of $H_{m,\varepsilon}$ to 
the diagonal $\Delta(X^{r})$ of $X^{r}$.  
Then since the restriction of (\ref{me}) to $\Delta(X^{r})\simeq X$
is isomorphic to 
\begin{equation}\label{resdiagonal}
L:= \left(\frac{1}{m} + \varepsilon\right)r(\lfloor m(K_{X/Y}+D)\rfloor) - \varepsilon\cdot f^{*}\det \left(f_{*}\mathcal{O}_{X}(\lfloor m(K_{X/Y}+D)\rfloor)\right)^{**}, 
\end{equation} 
this implies that $h_{m,\varepsilon}$ is considered to be  a  
singular hermitian metric on  $L$ with semipositive curvature 
in the sense of current and for every $y\in Y_{0}$, 
$h_{m,\varepsilon}|X_{y}$ is an AZD of the restriction of $L|X_{y}$.

Now we shall prove  the assertion (2) in Theorem \ref{logmain}. 
Let $G$ be a Cartier divisor on $X$ which is  $\mathbb{Q}$-linear equivalent to 
$K_{X/Y}+D$.   
For a positive rational number $t$, we consider 
\begin{equation}
\Theta_{t}:= (1+t)D^{r} + \varepsilon\Gamma
\end{equation}
instead of $\Theta$ in (\ref{theta}) above, 
where we have considered $D$ (resp.  $D^{r}$) 
as the Cartier divisor $G - K_{X/Y}$ (resp. $G^{r}- K_{X^{r}/Y}$). 
Then as above, for the $\mathbb{Q}$-line bundle: 
\begin{equation}\label{lt}
L(t) : = rK_{X/Y} +(1+t)rD + \varepsilon rmG- \varepsilon\cdot f^{*}\det \left(f_{*}\mathcal{O}_{X}(mG)\right)^{**} 
\end{equation}
by using the relative canonical measure $d\mu_{can,(X^{r},\Theta_{t})/Y}$,
we obtain the  singular hermitian metric $h_{m,\varepsilon}(t)$ on $L(t)$ 
with semipositive curvature (in the sense of current) such that 
for every $y\in Y_{0}$, $h_{m,\varepsilon}(t)|X_{y}$ is an AZD of $L(t)|X_{y}$. 
Let $\ell$ be a positive integer such that 
$\ell \varepsilon$ is a positive integer. 
And we set $t = 1/\ell r$.  We note that if we have taken $\ell$ sufficiently large, then $t$ is enough small so that $(X^{r},\Theta_{t})$ is still a KLT pair. 
Then   
\begin{equation}
K_{X/Y} + \ell L(t)
= \left((1+\ell r) + \ell\varepsilon mr\right) G -\ell\varepsilon\cdot f^{*}\det \left(f_{*}\mathcal{O}_{X}(mG)\right)^{**}
\end{equation}
holds.  Let us fix  a $C^{\infty}$-volume form $dV_{Y}$ on $Y$.  

Let $(A,h_{A})$ be a sufficiently ample hermitian line bundle on $Y$  such that  for every point $y\in Y$, the followings hold:
\begin{enumerate}
\item[(1)] There exists a local coordinate 
$(U,t)$  with biholomorphic to a open unit polydisk $\Delta^{n}$ in $\mathbb{C}^{n}$
with center $O$ and $t(y) = 0$,  
\item[(2)] For $y\in Y$ and $(U,t)$, there exists a $\rho\in C^{\infty}(Y)$ such that 
$\mbox{Supp}\,\rho \subset\subset t^{-1}(\Delta^{n})$ and 
$\rho \equiv 1$ on a neighbourhood of $y$. 
\item[(3)] $\mbox{Ric}\, dV_{Y} +(2\dim Y)\sqrt{-1}\partial\bar{\partial}(\rho\log |t|)
+ \sqrt{-1}\Theta_{A}$ dominates a K\"{a}hler form on $Y$, 
where $|t| = \sqrt{\sum_{j =1}^{k}|t_{j}|^{2}}, t = (t_{1},\cdots, t_{k}) (k=\dim Y)$.
\end{enumerate}
Such a hermitian line bundle $(A,h_{A})$ certainly exists. 
Then since $h_{m,\varepsilon}(t)$ has semipositive curvature in the sense of current, by the choice of $(A,h_{A})$ and the $L^{2}$-extension theorem (\cite{o,o-t}),  
we see that for every $y\in Y_{0}$ and 
any element of 
\[
H^{0}(X_{y},\mathcal{O}_{X_{y}}((1+ (1+m\varepsilon)\ell r)G)\otimes \mathcal{I}(h_{m,\varepsilon}(t)^{\ell}|X_{y}))
\]
extends to an element of 
\[
H^{0}(X,\mathcal{O}_{X}((1+ (1+m\varepsilon)\ell r)G + f^{*}A - \ell\varepsilon f^{*}\det f_{*}\mathcal{O}_{X}(mG))\otimes \mathcal{I}(h_{m,\varepsilon}(t)^{\ell}). 
\] 
We note that $L(t)|X_{y} - (1+m\varepsilon)\ell r)G|X_{y} = D|X_{y}$ is effective. 
Then since $h_{m,\varepsilon}(t)|X_{y}$ is an AZD of $L(t)|X_{y}$, 
\begin{equation}
H^{0}(X_{y},\mathcal{O}_{X_{y}}((1+ (1+m\varepsilon)\ell r)G))= H^{0}(X_{y},\mathcal{O}_{X_{y}}((1+ (1+m\varepsilon)\ell r)G)\otimes \mathcal{I}(h_{m,\varepsilon}(t)^{\ell}|X_{y})) 
\end{equation}
holds. 
Hence we see that for every sufficiently large $\ell$ with $\ell\varepsilon 
\in \mathbb{Z}_{> 0}$, 
\[
f_{*}\mathcal{O}(\ell mG)\otimes (\det \left(f_{*}\mathcal{O}_{X}(mG)\right)^{**})^{-\ell \varepsilon}\otimes \mathcal{O}_{Y}(A)
\]
is globally generated over $Y_{0}$. 
Then (\ref{WS}) follows from the finite generation of the relative 
log canonical bundles (\cite{b-c-h-m}). 
  This completes the proof of Theorem \ref{logmain}. q.e.d.

Author's address\\
Hajime Tsuji\\
Department of Mathematics\\
Sophia University\\
7-1 Kioicho, Chiyoda-ku 102-8554\\
Japan \\
e-mail address: tsuji@mm.sophia.ac.jp  or h-tsuji@h03.itscom.net

\end{document}